\documentclass{amsart}

\usepackage[left=3cm,right=3cm,
    top=3cm,bottom=3cm,bindingoffset=0cm]{geometry}
\usepackage{graphicx}
\usepackage{hyperref}
\hypersetup{colorlinks,allcolors=black}

\usepackage{amsmath}
\usepackage{amsfonts}
\usepackage{amssymb}
\usepackage{amsthm}
\usepackage{mathtools}
\usepackage{tikz-cd}
\usepackage{fdsymbol}
\usepackage{dynkin-diagrams}

\usepackage{makecell}

\usepackage{biblatex}
\newcommand{\ie}{{i.\,e.\,}}
\newcommand{\eg}{{e.\,g.\,}}
\newcommand{\cf}{{cf.\,}}

\newcommand{\A}{\mathbb{A}}
\newcommand{\Ad}{\operatorname{Ad}}
\newcommand{\ad}{\operatorname{ad}}
\newcommand{\Aut}{\operatorname{Aut}}
\newcommand{\be}{\mathfrak{b}}
\newcommand{\B}{\mathcal{B}}
\newcommand{\C}{\mathbb{C}}

\newcommand{\const}{\operatorname{const}}
\newcommand{\Kir}{\mathcal{C}}
\newcommand{\coker}{\operatorname{coker}}

\newcommand{\diag}{\operatorname{diag}}
\newcommand{\D}{\mathcal{D}}

\newcommand{\E}{\mathcal{E}}
\newcommand{\End}{\operatorname{End}}
\newcommand{\ev}{\operatorname{ev}}
\newcommand{\G}{\mathbb{G}}
\newcommand{\g}{\mathfrak{g}}
\newcommand{\GL}{\mathrm{GL}}
\newcommand{\gl}{\mathfrak{gl}}
\newcommand{\Gr}{\operatorname{Gr}}
\newcommand{\h}{\mathfrak{h}}
\newcommand{\Ha}{\mathbb{Ha}}
\newcommand{\Hom}{\operatorname{Hom}}

\newcommand{\id}{\operatorname{id}}
\newcommand{\im}{\operatorname{im}}
\newcommand{\Int}{\operatorname{Int}}
\newcommand{\jay}{\mathfrak{j}}
\newcommand{\Lie}{\operatorname{Lie}}
\newcommand{\M}{\mathcal{M}}

\newcommand{\Min}{\operatorname{Min}}
\newcommand{\Mor}{\operatorname{Mor}}
\newcommand{\N}{\mathbb{N}}

\newcommand{\kay}{\mathfrak{k}}

\newcommand{\m}{\mathfrak{m}}
\newcommand{\Out}{\operatorname{Out}}
\newcommand{\Pe}{\mathbb{P}}
\newcommand{\p}{\mathfrak{p}}
\newcommand{\PGL}{\mathrm{PGL}}
\newcommand{\PSO}{\mathrm{PSO}}
\newcommand{\PSp}{\mathrm{PSp}}
\newcommand{\pt}{\mathrm{pt}}
\newcommand{\Q}{\mathbb{Q}}
\newcommand{\R}{\mathbb{R}}

\newcommand{\res}{\operatorname{res}}
\newcommand{\s}{\mathfrak{s}}

\newcommand{\SL}{\mathrm{SL}}
\newcommand{\esl}{\mathfrak{sl}}
\newcommand{\SO}{\operatorname{SO}}
\newcommand{\so}{\mathfrak{so}}

\newcommand{\esp}{\mathfrak{sp}}
\newcommand{\Spec}{\operatorname{Spec}}

\newcommand{\Sym}{\operatorname{Sym}}
\newcommand{\te}{\mathfrak{t}}
\newcommand{\tr}{\operatorname{tr}}
\newcommand{\wt}{\operatorname{wt}}
\newcommand{\Z}{\mathbb{Z}}

\newcommand{\inv}{^{-1}}

\theoremstyle{definition}
\newtheorem{defn}{Definition}[section]
\newtheorem{notn}[defn]{Notation}
\newtheorem{defprop}[defn]{Definition and Proposition}

\theoremstyle{plain}
\newtheorem{thm}[defn]{Theorem}
\newtheorem{prop}[defn]{Proposition}
\newtheorem{lem}[defn]{Lemma}
\newtheorem{cor}[defn]{Corollary}
\newtheorem{conj}[defn]{Conjecture}

\newtheorem{innerintrothm}{Theorem}
\newenvironment{introthm}[1]
{\renewcommand\theinnerintrothm{#1}\innerintrothm}
{\endinnerintrothm}
\newtheorem{innerintroconj}{Conjecture}
\newenvironment{introconj}[1]
{\renewcommand\theinnerintroconj{#1}\innerintroconj}
{\endinnerintroconj}

\theoremstyle{remark}
\newtheorem{rmk}[defn]{Remark}
\newtheorem{expl}[defn]{Example}

\DeclareRobustCommand{\SkipTocEntry}[5]{}

\title{Automorphisms of medium algebras and quasi-anisotropic real affine Grassmannians}
\author{Mischa Elkner}

\begin{document}
\begin{abstract}
    We study automorphisms of Hausel's medium algebras by identifying them with equivariant cohomology algebras of affine Schubert varieties of a complex reductive group $G$. The automorphisms in question are then given by the action of either complex automorphisms or real structures of $G$ on affine Schubert varieties. We focus on describing the fixed points of those automorphisms on the spectra of medium algebras, and prove that the torsion-free quotients of the corresponding coinvariant algebras are isomorphic to equivariant cohomology rings of fixed loci in affine Schubert varieties. In the case of real structures, this naturally leads to and involves the study of affine Grassmannians for quasi-anisotropic (also called quasi-compact) real reductive groups, closely related to recent work of O'Brien. We also discuss relations to automorphisms of Higgs bundle moduli and to symmetric pair variants of medium algebras. 
\end{abstract}

\maketitle

\setcounter{tocdepth}{1}
\tableofcontents

\section{Introduction}
\subsection{Medium algebra automorphisms and coinvariants}
Medium algebras are commutative algebras associated to dominant integral weights of complex reductive groups, appearing in several contexts. They were introduced by Hausel in \cite{Hausel2024_avatars} as subalgebras of Kirillov algebras\footnote{originally called classical family algebras} \cite{Kirillov2001_introduction}, and are in fact believed to be their centres. Hausel encountered those algebras in a different setting -- the study of Higgs bundle moduli -- where the spectra of medium algebras are proposed to describe Lagrangian closures of certain upward flows \cite{Hausel2023_ICM_mirror}. Finally, medium algebras can be identified \cite{Hausel2024_avatars} with the equivariant cohomology algebras of affine Schubert varieties, connecting them to geometric representation theory. The purpose of this article is to describe several natural automorphisms of medium algebras using this last point of view.

To be more precise now, let $G$ be a connected complex reductive group, $\g$ its Lie algebra, and $V^\lambda$ a finite-dimensional irreducible representation of highest weight $\lambda$. To this data, one assigns the \emph{Kirillov algebra} $\Kir^\lambda(\g) \coloneqq (\C[\g]\otimes \End(V^\lambda))^G$, a graded algebra over its subring $\C[\g]^G$, whose elements can be thought of as $G$-equivariant morphisms (of affine schemes) from $\g$ to $\End(V^\lambda)$. A family of such morphisms is obtained from the $G$-equivariant morphisms $\g\to \g$ by postcomposition with the representation map $\g\to \End(V^\lambda)$, and the unital $\C[\g]^G$-subalgebra generated by this family is the \emph{medium algebra} $\M^\lambda(\g)$. Medium algebras are commutative, as well as finite and free as modules over $\C[\g]^G$, and we may think of them in terms of the branched covering $\Spec \M^\lambda(\g) \to \g/\!/G$. A general fibre of this branched covering can be identified with the set of weights $\wt(\lambda)$ of $V^\lambda$. 

Two natural sources of automorphisms of medium algebras are the grading, corresponding to a $\C^\times$-action, and automorphisms of $G$. For $\eta\in \Aut(G)$, one has to assume that the twisted representation $\eta^*V^\lambda$ is isomorphic to $V^\lambda$; we shall then call $\lambda$ and $\eta$ compatible. In this case, $\eta$ induces an automorphism of $\M^\lambda(\g)$, and this automorphism turns out to depend only on the inner class of $\eta$, that is, its image in $\Out(G) = \Aut(G)/\Int(G)$. Writing $\Out(G)_\lambda$ for classes of automorphisms compatible with $\lambda$, we altogether have a map
\begin{equation}
    \Out(G)_\lambda\times \C^\times \to \Aut(\M^\lambda(\g)),
    \quad
    (a,\zeta)
    \mapsto
    \iota_{a,\zeta}.
\end{equation}

Aside from their natural construction, our interest in the automorphisms $\iota_{a,\zeta}$ comes from the following observations. Firstly, the (maximal) parameter group $\Out(G)\times \C^\times$ is reminiscent of the parameters for automorphisms of Higgs bundle moduli spaces studied by García-Prada and Ramanan \cite{GPRa2019_involutions}; given the aforementioned connection to Higgs bundles, this is more than a coincidence (see \cite[Sec.\,4]{GH2024_even_very_stable} and \S\ref{sec:Higgs}). Secondly, they realise symmetries of $\wt(\lambda)$ on general fibres of $\M^\lambda(\g)$ -- for instance, $\iota_{1,-1}$ turns out to correspond to the action of the longest Weyl group element on the set of weights (see \S\ref{subsec:MA_weight_action}). We will mostly be interested in the finite order automorphisms of the form $\iota_{a,1}$ or $\iota_{a,-1}$, which we shall refer to as \emph{type $+1$} and \emph{type $-1$} automorphisms, respectively.

Another perspective, central to this article, is that of equivariant cohomology. Namely, as observed in \cite{Hausel2024_avatars}, there is an isomorphism of algebras
\begin{equation}
\label{eq:intro_MA_eq_coh}
    \M^\lambda(\g) \cong H^*_{\Check{G}}(\Gr_{\Check{G}}^{\le \lambda},\C),
\end{equation}
where $\Check{G}$ is the Langlands dual group, $\Gr_{\Check{G}}$ its affine Grassmannian, and $\Gr_{\Check{G}}^{\le \lambda}$ the affine Schubert variety labelled by $\lambda$. Note that this is compatible with the isomorphism $\C[\g]^G \cong \C[\check{\g}]^{\Check{G}} \cong H^*_{\Check{G}}(\pt,\C)$ of the base rings. This setting also provides a natural source of automorphisms: any $\sigma\in \Aut(\Check{G})$ induces an automorphism of $\Gr_{\Check{G}}$ (also denoted $\sigma$), and as before we say that $\sigma$ is compatible with $\lambda$ if $\sigma(\Gr_{\Check{G}}^{\le\lambda}) = \Gr_{\Check{G}}^{\le\lambda}$. In this case, functoriality of equivariant cohomology produces an automorphism $\sigma^*$ of $H^*_{\Check{G}}(\Gr_{\Check{G}}^{\le \lambda},\C)$, and $\sigma^*$ turns out to depend only on the inner class of $\sigma$. Moreover, all of this applies also if we instead take $\sigma$ to be a real structure of $\Check{G}$ (that is, the conjugation corresponding to a real form). The starting point for this work is the following:
\begin{introthm}{A}(See Theorem \ref{thm:MA_autos_geometry}.)
    The isomorphism \eqref{eq:intro_MA_eq_coh} matches:
    \begin{itemize}
        \item type $+1$ automorphisms with those induced by automorphisms of $\Check{G}$, and
        \item type $-1$ automorphisms $\iota_{a, -1}$ where $a^2 = 1$ with those induced by real structures of $\Check{G}$.
    \end{itemize}
\end{introthm}

With these perspectives in mind, we now come to our main objective. To explain it, observe first that the automorphisms $\iota$ under consideration restrict to automorphisms of the base ring $\C[\g]^G$. In geometric terms, we have an automorphism $\iota$ of the branched covering $\Spec \M^\lambda \to \g/\!/G$, and we will be interested in its fixed points. In general, if $A$ is a ring and $\phi\in \Aut(A)$ an automorphism, then the fixed point scheme of $\phi$ on $\Spec A$ is given by $\Spec A_\phi$, where $A_\phi \coloneqq A/(a-\phi(a) \mid a\in A)$ is the \emph{coinvariant ring}. In our setting, $\M^\lambda(\g)_\iota$ becomes an algebra over $\C[\g]^G_\iota$, and we call this the \emph{coinvariant algebra}.\footnote{The name \emph{coinvariant algebra} is also sometimes used to refer to a quotient by positive-degree invariants, this is a different concept.} The setting is summarised by the following commutative diagram:
\begin{equation*}
\begin{tikzcd}
    \Spec (\M^\lambda(\g)_\iota)
        \arrow[r, equal]
        \arrow[d]
    & (\Spec \M^\lambda(\g))^\iota
        \arrow[r, hookrightarrow]
        \arrow[d]
    & \Spec \M^\lambda(\g)
        \arrow[r, "\iota"]
        \arrow[d]
    & \Spec \M^\lambda(\g)
        \arrow[d]
    \\
    \Spec (\C[\g]^G_\iota)
        \arrow[r, equal]
    & (\g/\!/G)^\iota
        \arrow[r, hookrightarrow]
    & \g/\!/G
        \arrow[r, "\iota"]
    & \g/\!/G.
\end{tikzcd}
\end{equation*}
The main goal of this article is to describe the coinvariant algebras of type $\pm 1$ automorphisms using the setting of equivariant cohomology. 

Similar coinvariant algebras have previously been considered in three related papers. Firstly, there is the work of Zveryk \cite{Zveryk2024_dynkin}, which in fact concerns the so-called \emph{big algebras}, also introduced in \cite{Hausel2024_avatars}. Every medium algebra is contained in a big algebra, and the automorphisms studied here extend to automorphisms of big algebras. Zveryk studies the coinvariant algebras of (extensions of) type $+1$ automorphisms, and relates the coinvariants of $G$-big algebras to big algebras for the orbit group $(\Check{G}^{\theta\circ})^\vee$. Here $\theta$ is a well-chosen automorphism $\Check{G}$, and $\Check{G}^{\theta\circ}$ denotes the identity component of its fixed point subgroup. The apparently simpler case of medium algebras is not covered in \cite{Zveryk2024_dynkin}.

Secondly, coinvariants for type $-1$ automorphisms appear in the Higgs bundle setting in \cite{GH2024_even_very_stable}, in the special case where $\lambda$ is minuscule. This special case was further studied in the author's previous work \cite{Elkner2026_minuscule}.\footnote{In this case, the medium, big, and Kirillov algebra all coincide, and \cite{Elkner2026_minuscule} is written in terms of Kirillov algebras.} The main result of \cite{Elkner2026_minuscule} describes coinvariant algebras as equivariant cohomology rings of real affine Schubert varieties (or partial flag varieties, in this case) for well-chosen real structures of $\Check{G}$. 

A common lesson of those previous works is as follows: since the equivariant cohomology action $\sigma^*$ of a complex automorphism or real structure $\sigma$ depends only on its inner class, one can and should choose an optimal representative of that inner class. In the complex automorphism setting, this role is played by \emph{pinning-preserving automorphisms} (called \emph{Dynkin automorphisms} in \cite{Zveryk2024_dynkin}), and the optimal real structures turn out to be the \emph{quasi-anisotropic} ones (called \emph{quasi-compact} in \cite{Elkner2026_minuscule}). To briefly define these terms, an automorphism $\sigma$ of $\Check{G}$ is called pinning-preserving if it preserves some pinning of $\Check{G}$, and a real structure is called quasi-anisotropic if its Cartan involution is pinning-preserving. (For more details, see \S\ref{subsubsec:autos_setup} and \S\ref{sec:ppqa}.) The reason why these are the preferred representatives has to do with the topology of their fixed loci, as we shall see below.

\subsection{Main results}
We now come to our approach to coinvariants. Although we treat complex automorphisms and real structures simultaneously later on, let us start with a complex automorphism $\sigma\in \Aut(\Check{G})$ for simplicity. In view of \eqref{eq:intro_MA_eq_coh}, we wish to find the coinvariant ring $H^*_{\Check{G}}(\Gr_{\Check{G}}^{\le \lambda},\C)_{\sigma^*}$ as an algebra over $H^*_{\Check{G}}(\pt, \C)_{\sigma^*}$, and it seems reasonable to first treat the global variant $H^*(\Gr_{\Check{G}},\C)_{\sigma^*}$. A simple guess is that this should be related to the fixed point subspace $(\Gr_{\Check{G}})^\sigma$, which can in fact be identified with $\Gr_{(\Check{G})^\sigma}$. Now, from the inclusions $\Check{G}^{\sigma\circ}\hookrightarrow \Check{G}$ and $\Gr_{\Check{G}^\sigma}\hookrightarrow \Gr_{\Check{G}}$ we obtain at least homomorphisms
\begin{align}
    &H^*_{\Check{G}}(\pt, \C)_{\sigma^*}
    \to
    H^*_{\Check{G}^{\sigma\circ}}(\pt, \C)
    \label{eq:intro:cplx_global}
    \\
    &H^*_{\Check{G}}(\Gr_{\Check{G}},\C)_{\sigma^*}
    \to 
    H^*_{\Check{G}^{\sigma\circ}}(\Gr_{\Check{G}^\sigma},\C),
    \label{eq:intro:cplx_global_base}
\end{align}
and the question is whether they are isomorphisms. The answer will depend on $\sigma$, and in fact already asking just for \eqref{eq:intro:cplx_global_base} to be surjective essentially forces us to assume that $\sigma$ is pinning-preserving.\footnote{See Remark \ref{rmk:res_surjectivity_higher_order} for the single exception.} Even then, there is an additional difficulty arising from the fact that $\Gr_{\Check{G}^\sigma}\hookrightarrow \Gr_{\Check{G}}$ need not be injective on the sets of connected components, which obstructs surjectivity in \eqref{eq:intro:cplx_global}. Since the contributions of different components to equivariant cohomology of $\Gr_{\Check{G}^\sigma}$ can be identified, one can suppress this component issue by restricting the image appropriately. We show:
\begin{introthm}{B.1}(See Theorem \ref{thm:coinv_global}.)
\label{introthm:cplx_global}
    If $\sigma$ is a pinning-preserving automorphism of $\Check{G}$, then the inclusions $\Check{G}^{\sigma\circ} \hookrightarrow \Check{G}$ and $\Gr_{\Check{G}^{\sigma}}\to \Gr_{\Check{G}}$ induce isomorphisms
    \begin{equation*}
        H^*_{\Check{G}}(\pt, \C)_{\sigma^*}
        \overset{\cong}{\longrightarrow}
        H^*_{\Check{G}^{\sigma\circ}}(\pt, \C)
        \quad
        \text{and}
        \quad
        H^*_{\Check{G}}(\Gr_{\Check{G}},\C)_{\sigma^*}
        \overset{\cong}{\longrightarrow}
        H^*_{\Check{G}^{\sigma\circ}}(\Gr_{\Check{G}^\sigma},\C)^{S(\kappa_\sigma)},
    \end{equation*}
    where $S(\kappa_\sigma)$ is a finite symmetric group acting on equivariant cohomology by permuting factors corresponding to connected components of $\Gr_{\Check{G}^\sigma}$ lying in a common connected component of $\Gr_{\Check{G}}$. 
\end{introthm}

In algebraic terms, Theorem \ref{introthm:cplx_global} is rather simple: the equivariant cohomology rings involved are direct products of polynomial rings, and the homomorphisms induced by restrictions correspond to deleting certain variables. The main work is to show that such a simple algebraic picture really does describe the equivariant topology. 

Returning to affine Schubert varieties, the situation becomes more complicated. Firstly, the issue of noninjectivity of $\pi_0(\Gr_{\Check{G}^\sigma})\to \pi_0(\Gr_{\Check{G}})$ translates to the fact that $(\Gr_{\Check{G}}^{\le\lambda})^\sigma$ need not be connected. This is now more serious than on the level of full affine Grassmannians, and we resort to avoiding those $G$ where this happens. Secondly, the coinvariant ring $H^*_{\Check{G}}(\Gr_{\Check{G}}^{\le \lambda},\C)_{\sigma^*}$ can have torsion elements as a module over the coinvariant base ring $H^*_{\Check{G}}(\pt, \C)_{\sigma^*}$, and we will have to quotient out all torsion elements. For any module (or algebra) $M$ over an integral domain, we write $M_\mathrm{tf}$ for the quotient by its torsion submodule (or ideal). 
\begin{introthm}{C.1}(See Theorem \ref{thm:coinv_main}.)
\label{introthm:cplx_main}
    Let $\lambda$ be a dominant weight of $G$, and $\sigma$ a pinning-preserving automorphism of $\Check{G}$ compatible with $\lambda$. Assume that $\Lie(G)$ does not contain a copy of $\esl_{2n+1}(\C)$, $n\in \Z_{>0}$, on which some power of $\sigma$ restricts to a nontrivial automorphism. Let $\iota\in \Aut(\M^\lambda(\g))$ be the type $+1$ automorphism corresponding to $\sigma^*$. Then we have compatible isomorphisms
    \begin{equation*}
        \C[\g]^G_\iota 
        \cong 
        H^*_{\Check{G}}(\pt, \C)_{\sigma^*}
        \cong
        H^*_{\Check{G}^{\sigma\circ}}(\pt, \C)
    \end{equation*}
    and
    \begin{equation*}
        (\M^\lambda(\g)_\iota)_\mathrm{tf}
        \cong
        (H^*_{\Check{G}}(\Gr_{\Check{G}}^{\le \lambda},\C)_{\sigma^*})_\mathrm{tf}
        \cong 
        H^*_{\Check{G}^{\sigma\circ}}(\Gr_{\Check{G}^\sigma}^{\le \lambda},\C).
    \end{equation*}
\end{introthm}

If we now take $\sigma$ to be a real structure of $\Check{G}$, then $\Gr_{\Check{G}}^\sigma$ is identified with the real affine Grassmannian $\Gr_{\Check{G}^\sigma}$ of the real reductive group $\Check{G}$. More precisely, we have so far tacitly identified affine Grassmannians of complex reductive groups with their complex points; by $\Gr_{\Check{G}^\sigma}$ we now mean the real points $\Gr_{\Check{G}^\sigma}(\R)$. This entails that $\Gr_{\Check{G}^\sigma}$ can consist of only one point, namely if $\Check{G}^\sigma$ is anisotropic (in other words, compact as a Lie group). Nevertheless, real affine Grassmannians generally have a more complicated geometry than their complex counterparts. For instance, due to the existence of complex affine pavings, the cohomology of complex affine Grassmannians is concentrated in even degrees, which need not be the case for real affine Grassmannians. 

Despite this, O'Brien \cite{obrien2024_splitting_rank} has recently obtained descriptions of the (equivariant) homology of certain real affine Grassmannians similar to earlier results of Yun--Zhu \cite{YZ2011_integral} in the complex setting. O'Brien's results apply to \emph{splitting-rank} groups\footnote{The name comes from a property of corresponding symmetric spaces.}, in particular the following simple groups of adjoint type:
\begin{itemize}
    \item the quaternionic projective general linear group $\PGL_{n,\Ha}$, $n\in \Z_{>0}$,
    \item the projective Lorentz groups $\PSO_{2n-1,1,\R}$, $n\in \Z_{>0}$, and
    \item a certain real form of $E_6$.
\end{itemize}

All of these are quasi-anisotropic (or quasi-compact) groups, as used in \cite{Elkner2026_minuscule}. However, the quasi-anisotropic groups include also $\PGL_{2n+1,\R}$ for $n\in \Z_{>0}$, which may appear surprising, given the general expectation that split forms of reductive groups should in a sense be the opposite of anisotropic. This additional family indeed causes some difficulties\footnote{The roots of these difficulties are the two adjacent nodes of a type $A_{2n}$ Dynkin diagram interchanged by its nontrivial diagram automorphism. This is the only instance of a Dynkin diagram automorphism interchanging adjacent nodes, and many constructions related to such automorphisms become more complicated when this happens.} but has to be included to obtain a complete description of coinvariants. One of these difficulties is that $\PGL_{2n+1,\R}$ has odd-dimensional affine Schubert varieties, which are not seen by Nadler's real geometric Satake equivalence \cite{Nadler2005_PSoRLG}. Since the work of Yun--Zhu, and its generalisation by O'Brien use the geometric Satake framework, we cannot easily adapt their methods to all quasi-anisotropic groups. 

Moreover, we again need to take care with connected components: this time, $\pi_0(\Gr_{\Check{G}^\sigma})\to \pi_0(\Gr_{\Check{G}})$ is injective, but there can be components of $\Gr_{\Check{G}}$ stable under $\sigma$ which nevertheless have no real points (see also \S\ref{subsec:qa_pi1}). Those components which do intersect the real locus will be called \emph{admissible}. We prove the following result analogous to Theorem \ref{introthm:cplx_global}:
\begin{introthm}{B.2}(See Theorem \ref{thm:coinv_global}.)
\label{introthm_real_global}
    Let $\sigma$ be a quasi-anisotropic real structure of $\Check{G}$, and write $(\Gr_{\Check{G}})_\mathrm{adm}$ for the union of $\sigma$-admissible connected components of $\Gr_{\Check{G}}$. The inclusions $\Check{G}^{\sigma\circ} \hookrightarrow \Check{G}$ and $\Gr_{\Check{G}^{\sigma}}\to \Gr_{\Check{G}}$ induce isomorphisms
    \begin{equation*}
        H^*_{\Check{G}}(\pt, \C)_{\sigma^*}
        \overset{\cong}{\longrightarrow}
        H^*_{\Check{G}^{\sigma\circ}}(\pt, \C)
        \quad
        \text{and}
        \quad
        H^*_{\Check{G}}((\Gr_{\Check{G}})_\mathrm{adm},\C)_{\sigma^*}
        \overset{\cong}{\longrightarrow}
        H^*_{\Check{G}^{\sigma\circ}}(\Gr_{\Check{G}^\sigma},\C).
    \end{equation*}
\end{introthm}
For the result analogous to Theorem \ref{introthm:cplx_main}, we have to exclude the family $\PGL_{2n+1,\R}$, $n\in \Z_{>0}$, and its covers. For the remaining cases, we prove:
\begin{introthm}{C.2}(See Theorem \ref{thm:coinv_main}.)
\label{introthm:real_main}
    Let $\lambda$ be a dominant weight of $G$, and $\sigma$ a quasi-anisotropic real structure of $\Check{G}$ compatible with $\lambda$, such that $\Lie(\Check{G}^\sigma)$ does not contain a copy of $\esl_{2n+1}(\R)$ for any $n\in \Z_{>0}$. Let $\iota\in \Aut(\M^\lambda(\g))$ be the type $-1$ automorphism corresponding to $\sigma^*$. Then we have compatible isomorphisms
    \begin{equation*}
        \C[\g]^G_\iota 
        \cong 
        H^*_{\Check{G}}(\pt, \C)_{\sigma^*}
        \cong
        H^*_{\Check{G}^{\sigma\circ}}(\pt, \C)
    \end{equation*}
    and
    \begin{equation*}
        (\M^\lambda(\g)_\iota)_\mathrm{tf}
        \cong
        (H^*_{\Check{G}}(\Gr_{\Check{G}}^{\le \lambda},\C)_{\sigma^*})_\mathrm{tf}
        \cong 
        H^*_{\Check{G}^{\sigma\circ}}(\Gr_{\Check{G}^\sigma}^{\le \lambda},\C).
    \end{equation*}
\end{introthm}
It is worth noting that taking torsion-free quotients really is necessary in general. In fact, already in the minuscule cases studied in \cite{Elkner2026_minuscule}, not all coinvariant algebras are free as modules. It follows from Theorem \ref{introthm:real_main} that the torsion-free coinvariant algebras are already free, so any nonfree coinvariant algebra must have torsion.

Let us also mention that we expect Theorem \ref{introthm:real_main} to hold without excluding copies of $\esl_{2n+1}(\R)$. The currently missing topological input for this is as follows: 
\begin{introconj}{D}(= Conjecture \ref{conj:typeA2_Gr_cohomology})
\label{introconj:typeA2_Gr_cohomology}
    The canonical maps $H^*(\Gr_{\PGL_{2n+1}}(\R),\Q)\to H^*(\Gr_{\PGL_{2n+1}}^{\le \lambda}(\R),\Q)$ are surjective for all dominant coweights $\lambda$. Moreover, the dimension of $H^*(\Gr_{\PGL_{2n+1}}^{\le \lambda}(\R),\Q)$ equals the number of Iwahori cells in $\Gr_{\PGL_{2n+1}}^{\le \lambda}(\R)$ whose labelling coweights are fixed by the longest element of the Weyl group. 
\end{introconj}
We expect that this conjecture can be proved using results of Kocherlakota \cite{Kocherlakota1995_IntegralHomology} and some combinatorial work. Such a proof has been carried out for $n=1$ (in a rather lengthy computation not included here), and we plan to address the general case in future work.

\subsection{Comments on methods}
Much of this article is built on the following fact: if $\sigma$ is a pinning-preserving automorphism of a complex reductive Lie algebra $\g$, then $\C[\g]^\g_{\sigma^*} \cong \C[\g^\sigma]^{\g^\sigma}$ (where $\sigma^*$ denotes precomposition with $\sigma\inv$). This can be checked using Kostant sections, and it easily implies our statements regarding coinvariants of the base rings (as observed previously in \cite{Elkner2026_minuscule}). 

What is less straightforward is how to apply this fact to the equivariant cohomology of affine Grassmannians. To do so, we first reduce (via standard equivariant formality arguments) to nonequivariant cohomology. We then seek to algebraically model the map $i^* \colon H^*(\Gr_{G},\C)\to H^*(\Gr_{G^\sigma}, \C)$ for connected complex reductive $G$ (playing the role of $\Check{G}$ above) with a pinning-preserving automorphism or quasi-anisotropic real structure $\sigma$, and $i$ denoting the inclusion $\Gr_{G^\sigma}\hookrightarrow \Gr_G$.\footnote{Especially in the real case, using algebraic methods seems necessary. While both real and complex affine Grassmannians have useful CW decompositions, these generally do not match nicely.} 

In a first step, we use the homotopy-equivalence of $\Gr_G$ with the based loop space $\Omega K$ of a maximal compact subgroup $K\le G$. If $\sigma$ is a complex automorphism, then we similarly have $\Gr_{G^\sigma}\simeq \Omega K^\sigma$, and $i$ corresponds to the inclusion $\Omega K^\sigma\to \Omega K$. For a (quasi-anisotropic) real structure, we instead have $\Gr_{G^\sigma}\cong \Omega (K/K^\sigma)$. In this case, $i$ corresponds to $\Omega C$, where $C\colon K/K^\sigma\to K$ is the Cartan embedding.

Now, the homology of a based loop space $\Omega X$ with complex coefficients is well-understood for reasonable spaces $X$. Precise assumptions will be given below, but let us also assume, for simplicity, that $X$ is connected and simply connected. The topological Milnor--Moore theorem then yields $H_*(\Omega X, \C)\cong U(\pi_*(\Omega X)\otimes_\Z \C)$ where $\pi_*(\Omega X)\otimes_\Z \C$ is the complexification of the direct sum of homotopy groups of $X$. This is a graded Lie algebra, and $U(\pi_*(\Omega X)\otimes_\Z \C)$ denotes its universal enveloping algebra. Since $\pi_i(\Omega X) = \pi_{i+1}(X)$, it suffices to understand $\pi_*(X)\otimes_\Z \C$. We want cohomology rather than homology, and this amounts to taking the dual $(\pi_*(X)\otimes_\Z \C)^*$.

To access this data, we use rational homotopy theory, or more precisely $\C$-homotopy theory. Namely, the space $X$ admits a \emph{$\C$-minimal Sullivan model} $(\Lambda(V_X),d)$; this is a differential graded algebra over $\C$ determined uniquely by certain properties, one of which is that its homology is $H^*(X,\C)$. The key property is that this data includes the graded vector space $V_X$, which is known to be isomorphic to $(\pi_*(X)\otimes_\Z \C)^*$. 

Returning now to our problem, we have to understand the minimal models of $K$, $K^\sigma$ and $K/K^\sigma$, as well as the maps between them induced by $K^\sigma\to K$ and $C\colon K/K^\sigma\to K$.\footnote{Here we assume that $K$ is simply connected; the general case requires some additional care.} While such questions can be difficult in general, it is fairly well known that the minimal models of compact Lie groups are given by their cohomology algebras, that is, $(\Lambda(V_K), d) = (H^*(K,\C),0)$. The reason for this is that $K$ is a formal space with $H^*(K,\C)$ a free graded (exterior) algebra. Of course, the same reasoning applies to $K^\sigma$. For $K/K^\sigma$ we derive the analogous properties using a classical Theorem of Samelson \cite{Samelson1941_Beitraege} and the quasi-anisotropy assumption. 

At this point, we may use a standard spectral sequence argument to pass from $K$ and $K^\sigma$ to the classifying spaces $BK$ and $BK^\sigma$, and this turns the problem into the key fact mentioned at the start of this subsection. The case of $K/K^\sigma$ is handled similarly.

So far, this summary explains the proof of Theorems \ref{introthm:cplx_global} and \ref{introthm_real_global}, concerning coinvariant algebras on the level of full affine Grassmannians. Intuitively, our main results about (medium algebras and) affine Schubert varieties should then follow, but this appears not to be quite so simple. Nevertheless, we may conclude that the maps $H^*_{\Check{G}}(\Gr_{\Check{G}}^{\le \lambda},\C)_{\sigma^*}\to H^*_{\Check{G}^{\sigma\circ}}((\Gr_{\Check{G}}^{\le\lambda})^\sigma,\C)$ underlying Theorems \ref{introthm:cplx_main} and \ref{introthm:real_main} are surjective. The theorems are then derived by showing that these maps induce isomorphisms on a general fibre.

\subsection{Further topics, remarks and discussion}
\subsubsection{Torsion in coinvariant algebras}
Regarding the appearance of torsion in coinvariant algebras, we have the following conjecture generalising \cite[Thm.\,1.2]{Elkner2026_minuscule}:
\begin{introconj}{E}(See Conjecture \ref{conj:fixed_lambda_1}.)
    Let $\lambda$ be a dominant weight, defined in terms of a maximal torus $T$, and $\iota$ a type $-1$ automorphism of $\M^\lambda(\g)$. Suppose that a splitting-rank\footnote{Explicitly, we mean that $\sigma$ should be quasi-anisotropic and the corresponding real form of the Lie algebra should not contain any summands of the form $\esl_{2n+1}(\R)$, $n>0$.} real structure $\sigma$ of $\Check{G}$ with $\sigma^*$ corresponding to $\iota$ can be chosen such that $\sigma(\lambda) = \lambda$. Then $\M^\lambda(\g)_\iota$ is torsion-free.
\end{introconj}
\begin{rmk}\leavevmode
\begin{enumerate}
    \item One can assume that the real structure $\sigma$ preserves the maximal torus of $\Check{G}$ implicit in the setup, and then compatibility translates to $\sigma(\lambda) \in W\cdot \lambda$, where $W$ is the Weyl group. However, it is not always possible to improve this to $\sigma(\lambda) = \lambda$, which would explain the occurrence of torsion.
    \item It is also worth noting that for splitting-rank groups all coweights fixed by the real structure are of the form $\sigma(\mu) + \mu$ for some coweight $\mu$. These are the coweights labelling the affine Schubert varieties relevant to Nadler's real geometric Satake equivalence, which could be helpful in proving the conjecture.
\end{enumerate}
\end{rmk}

\subsubsection{Real affine Grassmannians}
Theorem \ref{introthm_real_global} provides an apparently novel point of view on the cohomology of real affine Grassmannians (for quasi-anisotropic groups). Upon dualisation, this also applies to homology, where one could use it to reprove some results of O'Brien \cite{obrien2024_splitting_rank}. We illustrate the connection in \S\ref{sec:coh_context}. However, we work with complex coefficients, and the results referred to hold over the integers localised at certain (finitely many) primes. We briefly discuss in \S\ref{subsec:coefficients} below how our methods could be adapted accordingly, but leave the details for future work.

That being said, even with complex coefficients, the real general linear group $\SL_{2n+1,\R}$ is accessible through our setting but not covered by \cite{obrien2024_splitting_rank}. The affine Grassmannians for splitting-rank groups admit cohomology descriptions in terms of their Nadler groups (that is, the associated groups from \cite{Nadler2005_PSoRLG}). For $\SL_{2n+1,\R}$ there is a similar-looking description, but in terms of a different group, see \S\ref{subsubsec:A2n-1_case}.

On a more conceptual level, it also seems interesting to compare our rather concrete setting to the usual real geometric Satake setup initiated in \cite{Nadler2005_PSoRLG}, in particular the derived versions in \cite{CMNO2025_quaternionic, CO2024_lorentzian_octonionic}. The latter pass from complex to real points using a well-chosen nearby cycles functor; the passage here corresponds to pullback along the inclusion map.

\subsubsection{Moduli of Higgs bundles}
As mentioned at the start of the introduction, Hausel's introduction of medium algebras was in relation to Higgs bundles, more precisely to conjectures on mirror symmetry \cite{Hausel2023_ICM_mirror}.\footnote{We should note that a more central role in these conjectures is played by big algebras.} Although the precise relation has (to the author's knowledge) not been worked out, this seems to explain the similarities between our automorphisms and those for Higgs bundle moduli described in \cite{GPRa2019_involutions}. In \S\ref{sec:Higgs}, we consider a known case of the medium-algebra-Higgs-bundles relation and show that the automorphisms indeed match.

\subsubsection{Symmetric pairs}
There is a well-known generalisation of Lie algebras to symmetric pairs, where one views a Lie algebra $\g$ as the symmetric pair $(\g\oplus \g,\g)$ obtained from the entry swap automorphism $(x,y)\mapsto (y,x)$ of $\g\oplus \g$. Along these lines, one can generalise medium algebras to \emph{symmetric pair medium algebras}, which we briefly discuss in \S\ref{sec:sym_pairs}. The connection to our main topic stems from the fact that certain symmetric pair medium algebras are also candidates for coinvariant algebras of type $-1$ automorphisms. (For type $1$, the corresponding candidates are usual medium algebras for fixed point subalgebras.) A more precise formulation of this is given in Conjecture \ref{conj:principal_coinv}.  

This is interesting for two reasons. Firstly, the above discussion of Higgs bundles suggests that type $-1$ coinvariant medium algebras should be related to Higgs bundles fixed by certain involutions, but our quasi-anisotropic groups do not easily fit the fixed points given in \cite{GPRa2019_involutions}. Symmetric pairs seem to match these more closely. Secondly, if we forget the coinvariants entirely, the conjecture combines with our main result to suggest that quasi-split (or principal) symmetric pair medium algebras are isomorphic to the equivariant cohomology of quasi-anisotropic affine Schubert varieties. This would extend the known relation between usual medium algebras and complex affine Schubert varieties. 

\subsection{Structure of the article}
\S \ref{sec:prelim} reviews background material related to algebraic topology and reductive groups, thereby also setting up notation. 

\S \ref{sec:MA} is about medium algebras and their automorphisms. It contains details on the isomorphism with equivariant cohomology and the weight set interpretation of general fibres. After giving those, we introduce and discuss the automorphisms studied in this article. In particular, we prove Theorem \ref{thm:MA_autos_geometry} relating automorphisms in the algebraic and equivariant-topological settings, as well as Proposition \ref{prop:MA_weight_action} which describes the action on well-behaved fibres. Here there is some overlap with \cite{Elkner2026_minuscule}. \S\ref{subsec:examples} contains some explicit examples of medium algebras and their automorphisms.

In \S \ref{sec:ppqa} we review the properties of pinning-preserving complex automorphisms and quasi-anisotropic real structures to be used afterwards. None of this should be novel, but for quasi-anisotropic real structures, such a collection of properties was not to be found elsewhere.

\S \ref{sec:coinv} contains our main results on coinvariant algebras. It begins with coinvariants of equivariant cohomology of full affine Grassmannians, where we prove Theorem \ref{thm:coinv_global}. As part of the proof, we also obtain a nonequivariant version, Theorem \ref{thm:noneq_coinv}. We then move to affine Schubert varieties (hence medium algebras) where we briefly digress to discuss their fixed loci. After that, we prove our main result on coinvariant medium algebras, Theorem \ref{thm:coinv_main}. We then discuss how the proof differs from what may be the expected proof, and how this relates to torsion.

We then address several related topics, starting in \S\ref{sec:Higgs} with moduli spaces of Higgs bundles. There, we prove Theorem \ref{thm:HB_MA_autos}, which matches our automorphisms to Higgs bundle automorphisms in a special case. 

In \S\ref{sec:coh_context} we discuss the relation of this article to other recent work on real affine Grassmannians, and show how our approach relates to familiar descriptions in terms of Nadler groups. This also includes a less familiar instance, where a different group appears. 

Finally, in \S\ref{sec:sym_pairs} we comment on medium algebras for symmetric pairs and formulate Conjecture \ref{conj:principal_coinv}, connecting them to the coinvariant algebras studied here.

Appendix \ref{appx:case_proofs} contains case-by-case proofs of some statements about quasi-anisotropic real structures.

\subsection{Acknowledgements}
I thank Tamás Hausel for his support and guidance. I also thank Miguel González, Daniel Holmes, Jakub Löwit, and Vladyslav Zveryk for helpful discussions. I also thank Jakub Löwit for useful comments on a draft of this article.

This research project was carried out during the author's PhD studies. The initial planning was done from June to September 2025, with the main work carried out February to September 2026. During the first time period, the author was supported by the Austrian Science Fund (FWF) 10.55776/P35847, and during the second time period by the ERC grant ViaFiPoS, with grant ID 101199663.

\subsubsection{Statement on the use of generative AI}
Generative AI (ChatGPT versions 5.3 to 5.6) was used for literature search and example computations (not contained in the article). None was used for writing the article and the design of its mathematical arguments. After finishing the first version, an AI-generated review offered by Tamás Hausel (using ChatGPT version 6.1 Sol Ultra) found some errors, which the author then corrected independently. 

\subsection{Conventions}
All rings are by convention commutative and unital. The cohomology $H^*(X,R)$ of a space (with coefficients in a ring $R$) is by convention the \emph{direct sum} $\bigoplus_{n=0}^\infty H^n(X,R)$ of its cohomology groups, rather than the direct product. Thus, for instance, the cohomology rings of affine Grassmannians with complex coefficients are polynomial rings. The direct product variant of cohomology then arises from ours by a completion. We also note that the cohomology of the empty topological space is $0$.

Centralisers of elements in groups or Lie algebras are denoted using subscripts, for instance, $G_x$ is the centraliser of $x$, which depending on context could be an element of $G$ or of a space acted on by $G$. Nevertheless, fixed points are denoted with superscripts, so for instance $G^\sigma$ is the subgroup of fixed points of $\sigma\colon G\to G$. The identity component of a group $G$ is denoted $G^\circ$.

To avoid any confusion that may arise from simultaneously using real and complex points, we will be uncommonly careful in distinguishing linear algebraic groups from their Lie groups of real or complex points.

Lastly, we abuse notation in the following ways:
\begin{itemize}
    \item For a homomorphism of Lie/linear algebraic groups $\phi\colon G\to H$, we often denote its derivative $\Lie(G)\to \Lie(H)$ also by $\phi$.
    \item Similarly, if $f\colon A\to B$ is a ring homomorphism, then the resulting map $\Spec B\to \Spec A$ will also be denoted $f$.
\end{itemize}

\section{Preliminaries and setting}
\label{sec:prelim}
\subsection{Graded algebras}
\subsubsection{Generalities}
\label{subsubsec:generalities}
In the following, $k$ denotes a field of characteristic zero, all vector spaces are over $k$, and all graded vector spaces are understood to be graded by the natural numbers $\N$, which contain $0$. The ground field $k$ is thought of as a graded vector space concentrated in degree zero. The tensor product of graded vector spaces $V,W$ becomes a graded vector space via $(V\otimes_k W)^n \coloneqq \bigoplus_{p+q=n}V^p\otimes W^q$. We write $T_{V,W}\colon V\otimes_k W\to W\otimes_k V$ for the \emph{swap automorphism} defined by $T(v\otimes w) = (-1)^{pq} w\otimes v$ for $v\in V^p$ and $w\in W^q$. The (graded) dual of a graded vector space $V$ is the graded vector space $V^*\coloneqq \bigoplus_{n\ge 0}(V^n)^*$, where $(V^n)^*$ denotes the usual dual space. A graded vector space $V$ is said to be of \emph{finite type} if each component $V^n$ is finite-dimensional; for such $V$ we have $(V^*)^*\cong V$. Given a graded vector space $V$ and $i\in \Z$, we denote by $V[i]$ the graded vector space with $(V[i])^n = V^{n+i}$, $n\in \N$.

By a \emph{graded algebra} we always mean an associative unital $k$-algebra with grading $A = \bigoplus_{n\ge 0} A^n$ such that $A^m\cdot A^n \subseteq A^{m+n}$ for all $m,n\ge 0$. We may think of the multiplication and unit of $A$ as grading-preserving linear maps $\mu\colon A\otimes_k A\to A$ and $\eta\colon k\to A$. We call the graded algebra $A$ \emph{commutative} if $\mu = \mu\circ T_{A,A}$. Dualising the axioms of a (commutative) graded algebra, one obtains the notion of a (cocommutative) graded coalgebra. If $A$ is a graded algebra (coalgebra) of finite type, then $A^*$ is naturally a graded coalgebra (algebra), also of finite type.

A \emph{free commutative graded algebra} is one of the form $\Lambda(V) = S(V_0)\otimes_k E(V_1)$ -- here $V$ is a graded vector space, $V=V_0\oplus V_1$ its decomposition into graded subspaces concentrated in even, resp.\,odd degrees, $S(V_0)$ is the symmetric algebra of $V_0$ and $E(V_1)$ the exterior algebra of $V_1$, both with gradings inherited from the tensor grading.  

An \emph{augmentation} of a graded algebra $A$ is a grading preserving algebra homomorphism $A\to k$. The kernel of $\epsilon$ is denoted by $I(A)$. The graded vector space $I(A)/(I(A)\cdot I(A))$ is called the space of \emph{indecomposable elements} and denoted by $Q(A)$. The image of any section of $I(A)\twoheadrightarrow Q(A)$ generates $A$ as a unital algebra.\footnote{In general, $I(A)$ and $Q(A)$ depend on the choice of $\epsilon$.} We call a graded algebra \emph{connected} if $\dim A^0 = 1$, equivalently $A^0 = \eta(k)$; in this case $A$ admits a unique augmentation. 

\subsubsection{Graded Hopf algebras}
Here we follow \cite[ch.\,20--22]{MayPontoMCAT}. 
A \emph{graded bialgebra} $(A,\mu,\eta, \Delta, \epsilon)$ (usually denoted only by $A$) consists of a graded vector space $A=\bigoplus_{n\ge 0} A^n$ and grading-preserving linear maps $\mu\colon A\otimes_k A\to A$, $\eta\colon k\to A$, $\Delta\colon A\to A\otimes A$, $\epsilon\colon A\to k$ such that
\begin{itemize}
    \item $(A,\mu, \eta)$ is a graded algebra,
    \item $(A,\Delta, \epsilon)$ is a graded coalgebra, and
    \item $\Delta$ and $\epsilon$ are algebra homomorphisms.
\end{itemize}
Let us observe that $\epsilon$ is then a canonical augmentation of the algebra $(A,\mu, \eta)$; $I(A)$ and $Q(A)$ (see \ref{subsubsec:generalities}) will always be defined via this augmentation.

A morphism of graded bialgebras is by definition a grading-preserving linear map compatible with all structure maps. If $A$ is of finite type, then its \emph{dual} $(A^*, \Delta^*, \epsilon^*, \mu^*, \eta^*)$ is also a graded bialgebra. A grading-preserving $k$-linear map $S\colon A\to A$ is called an \emph{antipode} if
\begin{equation*}
    \mu\circ (S\otimes \id)\circ \Delta 
    = \eta\circ \epsilon
    = \mu\circ (\id\otimes S)\circ \Delta.
\end{equation*}
A graded bialgebra $A$ can admit at most one antipode; if it does, one calls $A$ a \emph{graded Hopf algebra}. The dual of a finite type graded Hopf algebra is also a graded Hopf algebra. It is straightforward to equip the tensor product of two graded bialgebras (or graded Hopf algebras) with the structure of a graded bialgebra (resp.\,graded Hopf algebra). 

By a \emph{graded Lie algebra} we mean a graded vector space $L = \bigoplus_{n\ge 0}L^n$ together with a grading-preserving linear map $[\cdot, \cdot] \colon L\otimes_k L\to L$ such that for all $x\in L^m$, $y\in L^n$ and $z\in L^p$ we have
\begin{itemize}
    \item $[x,y] = (-1)^{mn+1}[y,x]$, and
    \item $(-1)^{mp}[x, [y,z]] + (-1)^{mn}[y,[z,x]] + (-1)^{np}[z,[x,y]] = 0$.
\end{itemize}
The \emph{universal enveloping algebra} $U(L)$ of a graded Lie algebra $L$ can be constructed similarly as in the nongraded case, that is, as the quotient of the graded tensor algebra $T(L)$ by the two-sided ideal generated by $x\otimes y - \Pi(x\otimes y) - [x,y]$ for all homogeneous $x,y\in L$. It is a graded Hopf algebra with the algebra structure inherited from $T(L)$ and coproduct induced by $\Delta(x) = x\otimes 1 + 1\otimes x$ for $x\in L$; if moreover $L^0=0$, then $U(L)$ is connected. Any graded $k$-vector space $V$ can be viewed as a graded Lie algebra with trivial bracket; in this case $U(V)$ agrees as an algebra with the free commutative graded algebra $\Lambda(V)$.

The \emph{primitive elements} $P(A)$ of a graded bialgebra $A$ are given by
\begin{equation*}
    P(A)\coloneqq \{x\in A\mid \Delta(x) = 1\otimes x + x\otimes 1\};
\end{equation*}
they form a graded $k$-linear subspace closed under the commutator of $(A,\mu)$, hence a graded Lie algebra. If $x\in P(A)$, the counit property implies $\epsilon(x) = 0$, thus $P(A)\subseteq I(A)$. It follows that any graded bialgebra has a canonical map $\pi_A\colon P(A)\to Q(A)$; this map fulfils the following:
\begin{prop}[{\cf\cite[Cor.\,22.3.3]{MayPontoMCAT}}]
\label{prop:HA_P-to-Q}
    Let $A$ be a connected graded Hopf algebra. Then $\pi_A$ is injective if and only if $A$ is commutative and surjective if and only if $A$ is cocommutative. 
\end{prop}
One can also define $P(A)$ as follows: since $\Delta$ sends $1$ to $1\otimes 1$, it descends to a map $\underline{\Delta} \coker \eta \to \coker \eta \otimes \coker \eta$, and $P(A)$ is mapped isomorphically onto its image in $\coker \eta$. Identifying $P(A)$ with that image, one has $P(A) = \ker \underline{\Delta}$. This description is precisely dual to that of $Q(A)$, so if $A$ is of finite type we have $P(A) \cong Q(A)^*$.

The following classical structure theorems will be useful:
\begin{thm}[{\cf\cite[Thm.\,22.4.1 and Cor.\,22.4.3]{MayPontoMCAT}}]
\label{thm:HA_structure_1}
    Let $A$ be a connected commutative graded Hopf algebra. Then any choice of section $Q(A)\to I(A)$ yields a graded algebra isomorphism $A\cong \Lambda(Q(A))$. If $Q(A)$ is concentrated in odd degrees, then $A$ is cocommutative and $P(A)\cong Q(A)$. 
\end{thm}
\begin{thm}[{cf.\,\cite[Thm.\,22.3.1, Cor.\,22.3.3]{MayPontoMCAT}}]
\label{thm:HA_structure_2}
    Let $A$ be a connected cocommutative graded Hopf algebra. Then $A$ is canonically isomorphic to $U(P(A))$.
\end{thm}

\subsubsection{Differential graded algebras}
A \emph{differential graded algebra} $(A,d)$ (often just written as $A$) is a graded algebra $A$ together with a graded derivation $d\colon A\to A$ of degree one: $d$ maps each $A^n$ to $A^{n+1}$ and fulfils $d(ab) = d(a)b + (-1)^{pq}ad(b)$ for all $a\in A^p$, $b\in A^q$. We use the common abbreviations \emph{dga} for ``differential graded algebra'' and \emph{cdga} for ``commutative differential graded algebra''. The cohomology algebra of a dga $A$ is $H(A)\coloneqq \ker(d)/\im(d)$, it is again a graded algebra. A morphism of dgas $f\colon A\to B$ is a graded algebra morphism such that $d\circ f = f\circ d$; it induces a morphism $H(f)\colon H(A)\to H(B)$. If $H(f)$ is an isomorphism, one calls $f$ a \emph{quasi-isomorphism}. The equivalence relation between dgas generated by quasi-isomorphisms is called \emph{weak equivalence}. On any graded algebra $A$, the zero map is a graded derivation making $A$ a dga; one calls a dga that is weakly equivalent to one with $d=0$ \emph{formal}. Equivalently, a dga $A$ is formal if and only if it is weakly equivalent to $(H(A),0)$.
   
\subsection{Rational homotopy theory}
\label{subsec:prelim_RHT}
We follow \cite[\S\S10--16]{FHT_RHT}, which may be consulted for proofs and further details. Throughout this subsection, $X$ denotes a connected, simply connected, compactly generated Hausdorff space, and $k$ a field of characteristic zero. One can define a $k$-cdga $A_\mathrm{PL}(X,k)$ of piecewise linear differential forms on $X$; any quasi-isomorphism $m\colon A\to A_\mathrm{PL}$ from a cdga $A$ is called a \emph{($k$-)model} of $X$. There are preferred models defined as follows. A \emph{Sullivan algebra} is a cdga $(\Lambda(V),d)$ with free underlying commutative graded algebra $\Lambda(V)$, such that $V^0 = 0$, and such that $V$ is the increasing union of graded subspaces $V(n)$, $n\in \N$ with $d(V(n)) \subseteq V(n-1)$ for all $n\ge 1$ and $d(V(0)) = 0$. Observe that $\Lambda(V)$ is then connected. A Sullivan algebra $(\Lambda(V), d)$ is called \emph{minimal} if $d(V) \subseteq I(\Lambda(V))\cdot I(\Lambda(V))$. A model $m\colon (\Lambda(V), d)\to A_\mathrm{PL}(X,k)$ by a minimal Sullivan algebra is called a \emph{($k$-)minimal model} of $X$. We collect some fundamental facts in the following:
\begin{prop}
\label{prop:RHT_basics}
    Let $X$ and $k$ be as above. Then:
    \begin{itemize}
        \item[(i)] $X$ admits a minimal model.
        \item[(ii)] Any two minimal Sullivan algebras weakly equivalent to $A_\mathrm{PL}(X,k)$ are isomorphic.
        \item[(iii)] If $k=\R$ and $X$ is a smooth manifold, then $A_\mathrm{PL}(X,k)$ is weakly equivalent to the de Rham cdga $A_\mathrm{dR}(X)$ of differential forms on $X$.
        \item[(iv)] $A_\mathrm{PL}(X,k)$ is weakly equivalent to the dga $C^*(X,k)$ of singular cochains on $X$ with coefficients in $k$.  
    \end{itemize}
\end{prop}

By parts (i) and (ii) of the above proposition, we may assign to $X$ an essentially unique minimal model, which will be denoted $\Min(X,k)$, or $\Min(X)$ if the choice of $k$ is clear. It is easy to see that $\Min(X,k)\otimes_k k'\cong \Min(X,k')$ whenever $k\to k'$ is a field extension. Observe further that a (minimal) Sullivan algebra $(\Lambda(V), d)$ is fully determined by $V$ and $d|_V\colon V\to \Lambda(V)$. If $(\Lambda V, d) = \Min(X)$, we write $V_X\coloneqq V$ and note that $V_X \cong Q(\Min(X))$. Various constructions involved here are functorial. In particular, if $Y$ is another topological space fulfilling the same assumptions as $X$, and $f\colon X\to Y$ a continuous map, then $f$ gives rise to a morphism of dgas $\Min(f,k)\colon \Min(Y,k)\to \Min(X,k)$ (called its Sullivan representative) inducing $f^*$ in singular cohomology, and further to a morphism $Q(f)\colon V_Y\to V_X$ of graded vector spaces.

By part (iv) of Proposition \ref{prop:RHT_basics}, the cohomology algebra of $\Min(X,k)$ is isomorphic to singular cohomology. In particularly favourable situations, the singular cohomology is itself a minimal model: let us call $X$ \emph{formal} if $A_\mathrm{PL}(X,k)$ (equivalently $C^*(X,k)$ or, if $k=\R$, $A_\mathrm{dR}(X,k)$) is a formal dga. It then suffices to note that any free graded commutative algebra with zero differential is a minimal Sullivan algebra to see the following: 
\begin{lem}
\label{lem:formal_min_model}
    If $X$ is formal and $H^*(X,k)$ is a free graded commutative algebra, then $\Min(X,k) \cong (H^*(X,k), 0)$. 
\end{lem}

The topological space $X$ is said to be of \emph{finite type} (with respect to $k$) if $H^*(X,k)$ is of finite type, \ie if all cohomology groups are finite-dimensional; we now assume this. A key advantage of minimal models is then that they give access to homotopy groups, as follows. We consider $\pi_*(X)\otimes_\Z k$ as a graded vector space equipped with the \emph{Whitehead product}, whose definition we briefly recall. For all $m,n\ge 1$, one defines a so-called universal Whitehead product $a_{m,n}\in \pi_{m+n-1}(S^m\vee S^n)$ by viewing $S^{m+n-1}$ as the boundary of a cube $I^m\times I^n$, and mapping $I^m$, $I^n$ to $S^m$ resp. $S^n$ by contracting their boundaries. Precomposition with $a_{m,n}$ yields a map
\begin{equation}
\label{eq:Whitehead_product}
    [\cdot, \cdot]_W\colon \pi_m(X)\times \pi_n(X)\to \pi_{m+n-1}(X),
    \quad
    [f,g]_W \coloneqq \big(S^{m+n-1} \overset{a_{m,n}}{\longrightarrow} S^m\vee S^n \overset{f\vee g}{\longrightarrow} X \big).
\end{equation}
Extending coefficients to $k$ and shifting degrees, the Whitehead product $[\cdot,\cdot]_W$ becomes a bilinear map $(\pi_*(X)\otimes_Z k)[1]\otimes (\pi_*(X)\otimes_X k)[1] \to (\pi(X)\otimes_Z k)[1]$.\footnote{The bracket $[\cdot, \cdot]_W$ is in general similar to but not quite the same as a Lie bracket. See also \eqref{eq:Whitehead_Samelson} below.} 

Now, if $\Min(X) = (\Lambda V_X, d)$, note first that $d(V_X)\subseteq \Lambda^{\ge 2}(V_X)$, so that the projection $\Lambda(V_X)\to V_X$ induces a map $\zeta_X\colon H^+(\Lambda(V_X)\to V_X$. We further let $d_1\colon V_X\to \Lambda^2 V_X$ be the $\Lambda^2V_X$-component of $d|_{V_X}$; it is called the \emph{quadratic part} of $d$. 
\begin{thm}[{\cf\cite[\S13(c), Prop.\,13.16, Thm.\,15.11, and Cor.\,on p.\,210]{FHT_RHT}}]
\label{thm:RHT_main}
    There is a perfect pairing $\langle \cdot, \cdot \rangle \colon V_X\times (\pi_*(X)\otimes_\Z k)\to k$. Thus, $\pi_*(X)\otimes_\Z k$ is canonically isomorphic to $V_X^*$. Moreover:
    \begin{itemize}
        \item[(i)]
        The isomorphisms $H^*(\Lambda(V_X)) \cong H^*(X,k) \cong (H_*(X,k))^*$ and $V_X\cong (\pi_*(X)\otimes_\Z k)^*$ identify $\zeta_X\colon H^+(\Lambda(V_X))\to V_X$ with the dual of the Hurewicz homomorphism.
        \item[(ii)]
        If we define the trilinear map $\langle \cdot; \cdot, \cdot\rangle\colon \Lambda^2V_X\times (\pi_*(X)\otimes_\Z k)\times (\pi_*(X)\otimes_\Z k) \to k$ by setting
    \begin{equation*}
        \langle v\wedge w; f, g\rangle 
        \coloneqq
        \langle v, f\rangle \langle w, g\rangle +(-1)^{mn}\langle v, g\rangle\langle w, f\rangle
    \end{equation*}
    for all $v,w\in V_X$, $f\in \pi_m(X)\otimes_Z k$ and $g\in \pi_n(X)\otimes_\Z k$, then (for all such $v,f,g$):
    \begin{equation}
    \label{eq:d1_and_Whitehead}
        \langle v, [f,g]_W\rangle 
        =
        (-1)^{m+n-1}\langle d_1(v); f,g\rangle.
    \end{equation}
    \end{itemize}
\end{thm}

\subsection{Cohomology of H-spaces}
\label{subsec:prelim_Hsp}
Our main reference here is \cite[\S9]{MayPontoMCAT}.
\subsubsection{Definitions}
By an \emph{H-space} we mean a well-pointed topological space $(X, e)$ equipped with a map $m\colon X\times X\to X$ such that the maps $m(e,\cdot), m(\cdot,e)\colon X\to X$ are both homotopic to the identity. It is called \emph{homotopy-associative}, or an \emph{H-monoid}, if $m\circ (m\times \id)$ is homotopic to $m\circ (\id \times m)$. An H-monoid is called an \emph{H-group} if there exists a map $\nu\colon X\to X$ such that $x\mapsto m(x,\nu(x))$ and $x\mapsto m(\nu(x), x)$ are homotopic to the constant map $\const_e$ with value $e$. We call $X$ \emph{homotopy-commutative} if the maps $X\times X\to X$ given by $m$ and $(x,y)\mapsto m(y,x)$ are homotopic.

\subsubsection{The Samelson product}
\label{subsubsec:Samelson_product}
If $X$ is an $H$-group, we may define a map $\operatorname{comm}\colon X\times X\to X$ by
\begin{equation*}
    \operatorname{comm}(x,y) \coloneqq \mu(\mu(x,y), \mu(\nu(x), \nu(y))); 
\end{equation*}
this is the analogue of a the commutator map of a topological group. After passing to a homotopic map if necessary, we may assume that $\operatorname{comm}$ descends to the smash product $X\wedge X$. The \emph{Samelson product} of $X$ is then
\begin{equation}
    [\cdot, \cdot]_S\colon \pi_m(X)\times \pi_n(X)\to \pi_{m+n}(X),
    \quad
    [f,g]_S \coloneqq \big(S^{m+n} \cong S^m\wedge S^n\overset{f\wedge g}{\longrightarrow} X\wedge X\overset{\operatorname{comm}}{\longrightarrow} X\big).
\end{equation}
\begin{prop}[{\cf\cite[ch.\,X, Thms.\,5.1 and 5.4]{WhiteheadEHT}}]
    If $X$ is path-connected, the Samelson product makes $\pi_*(X)$ into a graded Lie algebra over $\Z$.
\end{prop}
For a field $k$, we may extend coefficients to $k$ to obtain $\pi_*(X)\otimes_\Z k$ as a graded Lie algebra over $k$ (whose bracket we shall also call Samelson product).

\subsubsection{Cohomology of connected H-spaces}
The multiplication map $m$ of an $H$-monoid of finite type induces $m^*\colon H^*(X,k)\to H^*(X\times X,k) \cong H^*(X,k)\otimes H^*(X,k)$; this makes $(H^*(X,k),\smile,\eta,m^*,\epsilon)$ a graded bialgebra, where $\smile$ is the cup product, $\eta$ is induced by $\const_e$, and $\epsilon$ is the usual augmentation. If $X$ is an $H$-group with inversion map $\nu$, then $H^*(X,k)$ is a Hopf algebra with antipode $\nu^*$. One can similarly make $H_*(X,k)$ into a graded bialgebra (or Hopf algebra) and this agrees with the dual of $H^*(X,k)$. We note that these graded bialgebras are connected if and only if $X$ is path-connected.

We now let $X$ be a path-connected H-group of finite type. Consider the Hurewicz homomorphisms $h_n\colon \pi_n(X,e)\otimes_\Z k \to H_n(X,k)$, $n\ge 1$. It follows easily from the constructions that any image of any $h_n$ is primitive, and in fact the following is true:
\begin{thm}[{\cf\cite[appendix]{MM1965_OtSoHA} or \cite[Prop.\,9.2.4]{MayPontoMCAT}}]
\label{thm:Hurewicz_primitives}
    If $X$ is of finite type and path-connected, the Hurewicz homomorphisms define an isomorphism of graded Lie algebras $\pi_*(X)\otimes_\Z k \to P(H_*(X,k))$.
\end{thm}

Being mainly interested in cohomology, we note that $Q(H^*(X,k)) \cong P(H_*(X,k))^*$. Moreover, $H^*(X,k)$ is always commutative (as a property of the cup product). If it is also cocommutative -- for instance when $X$ is homotopy-commutative -- then Proposition \ref{prop:HA_P-to-Q} yields an identification $P(H^*(X,k)) \cong Q(H^*(X,k)) \cong P(H_*(X,k))^* \cong (\pi_*(X)\otimes_\Z k)^*$. We write $P_X\coloneqq P(H^*(X,k))$. The structure theorems \ref{thm:HA_structure_1} and \ref{thm:HA_structure_2} now imply the following:
\begin{thm}
\label{thm:H-group_coh_structure}
    Let $X$ be a path-connected $H$-group of finite type.
    \begin{itemize}
        \item[(i)] If $H^*(X,k)$ is finite-dimensional (that is, $H^n(X,k) = 0$ for large $n$), then $P_X\cong (\pi_*(X), k)^*$ is concentrated in odd degrees and there is a Hopf algebra isomorphism $H^*(X,k) \cong E(P_X)$.
        \item[(ii)] If $H^*(X,k)$ is cocommutative, \eg when $X$ is homotopy-commutative, then we have a Hopf algebra isomorphism $H^*(X,k) \cong \Lambda(P_X) \cong \Lambda((\pi_*(X), k)^*)$.
    \end{itemize}
\end{thm}

\subsubsection{Disconnected $H$-spaces}
Finally, we address the case where $X$ is of finite type but not path-connected; note that $\pi_0(X)$ is nevertheless finite due to the finite type assumption. The $H$-group structure induces a group structure on $\pi_0(X)$, and $H^0(X,k)$ can be identified with $k[\pi_0(X)]$, the functions on $\pi_0(X)$. Since $H^0(X,k)$ is stable under all structure maps, this makes $k[\pi_0(X)]$ into a Hopf algebra (with trivial grading); in fact this is simply the usual Hopf algebra structure for $\pi_0(X)$ viewed as an affine algebraic group. For notational convenience, identify $\pi_0(X)$ as a set with chosen representatives (including $e$) of the components. The path component $X_a$ of $a\in \pi_0(X)$ is homotopy-equivalent to $X_e$ through $x\mapsto m(a,x)$, hence $H^*(X,k)\cong H^*(X_e\times \pi_0(X),k) \cong H^*(X_e,k) \otimes_k k[\pi_0(A)]$ as a graded algebra. To describe the comultiplication in this picture, consider the map $\pi_0(X)\times X_e\to X_e$ given by $(a,x)\mapsto m(a, m(x, \nu(a)))$. It is a group action up to homotopy, hence equips $H_*(X_e,k)$ with the structure of a module over the group algebra of $\pi_0(X)$; dualising yields a comodule structure $\rho\colon H^*(X_e,k)\to k[\pi_0(X)]\otimes_k H^*(X_e,k)$. One can now verify that $H^*(X,k)$ is the resulting smash coproduct (see \cite{Molnar1977_SDProducts}):
\begin{prop}
\label{prop:coh_of_disconnected_H-group}
    Let $(X,e)$ be an $H$-group of finite type and $X_e$ the path component of $e$. Then the graded Hopf algebra $H^*(X,k)$ has underlying graded $k$-algebra $H^*(X_e,k)\otimes k[\pi_0(X)]$. Denoting the (co)multiplications of $H^*(X_e,k)$ and $k[\pi_0(X)]$ by $\Delta', \mu'$ resp.\,$\Delta'', \mu''$, and the swap automorphism for tensor products by $T$, the comultiplication of $H^*(X,k)$ is given by
    \begin{equation*}
        (\id\otimes \mu'' \otimes \id \otimes \id) \circ (\id\otimes \id \otimes T\otimes \id) \circ (\id\otimes \rho\otimes \Delta'')\circ (\Delta'\otimes \id).
    \end{equation*}
\end{prop}
\begin{rmk}(Infinite component groups)
    If $X_e$ is of finite type but $\pi_0(X)$ is infinite, then $H^0(X,k)$ should be identified with $k^{\pi_0(X)}$ rather than $k[\pi_0(X)]$. This is still almost a Hopf algebra, but the identification of $H^0(X\times X, k)$ with $H^0(X,k)\otimes H^0(X,k)$ implicit in the coproduct is problematic. (For reasons like this, homology is often preferred.) Nevertheless, the same reasoning as above still yields an isomorphism of \emph{graded algebras} $H^*(X,k)\cong H^*(X_e,k)\otimes k^{\pi_0(X)}$. 
\end{rmk}

\subsection{Compact Lie groups and symmetric spaces}
Throughout this subsection, $K$ denotes a connected compact Lie group, and $\kay$ its Lie algebra.

\subsubsection{Universal bundle and transgression}
Since $K$ is in particular an H-group with finite-dimensional cohomology, Theorem \ref{thm:H-group_coh_structure}(i) applies, showing that $H^*(K,\R)$ is isomorphic to an exterior algebra $E(P_K)$. In principle, the space of primitive elements $P_K$ is then determined by the homotopy groups of $K$, but a more convenient description is as follows. Let $EK\to BK$ denote the universal $K$-bundle, thus $EK$ is a contractible CW complex with free $K$-action (unique up to $K$-homotopy) and $BK=EK/K$ the classifying space of $K$.
\begin{thm}[{\cf\cite[Sec.\,19--20]{Borel1953_esp_fib_et_hom}}]
\label{thm:transgression_EG}
    For the Leray--Serre spectral sequence associated to the fibration $K\to EK\to BK$, the transgressive elements of $H^*(K,\R)$ are precisely the primitive elements $P_K$. Moreover, the transgression takes $P_K$ isomorphically to the space of indecomposables $Q_{BK}\coloneqq Q(H^*(BG,\R))$, which is concentrated in even degrees, and any choice of section $Q_{BK}\to H^*(BG,\R)$ yields an isomorphism $S(Q_{BK}) \cong H^*(BG,\R)$.
\end{thm}

This is helpful because $H^*(BG,\R)$ is well known to coincide with the ring $\R[\kay]^K$ of $K$-invariant polynomials on $\kay$. In turn, if $T\le K$ is a maximal torus, $W=N_K(T)/T$ the Weyl group, and $\te$ the Lie algebra of $T$, then the Chevalley restriction theorem gives $\R[\kay]^K\cong \R[\te]^W$, yielding an easy description of $P_{BK}$ in terms of the degrees of the Weyl group. It is moreover worth pointing out that the assignments $K\mapsto EK, BK, P_K, P_{BK}$ are functorial in $K$ (the latter two contravariantly), and that the above identifications are natural. Thus, if $\phi\colon H\to K$ is a homomorphism of compact Lie groups, then the following diagram commutes for every $n\ge 0$:
\begin{equation}
\label{eq:transgression_naturality_diagram}
\begin{tikzcd}
    H^n(K,\R)
        \arrow[d, "\phi^*"]
    & P_K^n
        \arrow[l, hook']
        \arrow[r, "\tau_K"]
        \arrow[d]
    & Q^{n+1}_{BK}
        \arrow[d]
    & H^{n+1}(BK,\R)
        \arrow[l, twoheadrightarrow]
        \arrow[r, "\cong"]
        \arrow[d, "(B\phi)^*"]
    & \R[\kay]^K
        \arrow[d, "(\phi')^*"]
    \\
    H^n(H,\R)
    & P_H^n
        \arrow[l, hook']
        \arrow[r, "\tau_H"]
    & Q^{n+1}_{BH}
    & H^{n+1}(BH,\R)
        \arrow[l, twoheadrightarrow]
        \arrow[r, "\cong"]
    & \R[\h]^H
\end{tikzcd}
\end{equation}
Here, $\h$ denotes the Lie algebra of $H$, $\tau_K$, $\tau_H$ the appropriate transgressions, and $\phi'\colon \h\to \kay$ the derivative of $\phi$. 

\subsubsection{Compact symmetric spaces}
\label{subsubsec:cpct_sym_spaces}
Recall that a \emph{(compact, globally) symmetric space} is a quotient $K/H$ where $H\le K$ is a closed subgroup such that $K^{\sigma\circ} \subseteq H \subseteq K^\sigma$ for some involution $\sigma$ of $K$. For the symmetric spaces discussed in this article, we usually have $H=K^\sigma$. Such symmetric spaces can be identified with closed subspaces of $K$ via the \emph{Cartan embedding}
\begin{equation}
\label{eq:Cartan_embedding}
    C\colon K/H \to K,
    \quad
    xK^\sigma \mapsto x\sigma(x\inv).
\end{equation}
More generally, $K/H$ as above is a finite covering space of $K/K^\sigma$, and the map $C$ as above is well-defined but need not be injective. If $K$ is simply connected (or more generally $\pi_1((K^\sigma)^\circ) \to \pi_1(K)$ surjective), then $K/(K^\sigma)^\circ$ is the universal covering of $K/K^\sigma$.

It is useful to note that a compact Lie group $K$ is itself a symmetric space. For this one considers the swap involution $\sigma\colon (x,y)\mapsto (y,x)$ on $K\times K$; the map $K\times K\to K$ sending $(x,y)$ to $xy\inv$ then exhibits $K$ as $(K\times K)/(K\times K)^\sigma$. Thus, results about compact symmetric spaces may generalise results about compact Lie groups.

\subsubsection{Action of \texorpdfstring{$\sigma$}{σ} and formality}
\label{subsubsec:theta_and_formality}
As before, let $K/H$ be a compact symmetric space corresponding to the involution $\sigma$ of $K$. Since $\sigma$ fixes $H$ pointwise, it descends to an involution of $K/H$ (which is really the essence of the symmetric space). Now, using classical ideas of É.\,Cartan, one can model the real cohomology of $K/H$ by $K$-invariant differential forms. That is, one identifies $H^*(K/H,\R)$ with de Rham cohomology and then shows that the inclusion map $\Omega^*(K/H)^K\to \Omega^*(K/H)$ is a quasi-isomorphism (\cf\eg\cite[p.\,150, Thm.\,1]{Onishchik1994_tttg}). A key advantage is that $K$-invariant differential forms on $K/H$ are fully determined by their values at the distinguished tangent space $T_{eH}(K/H) \cong \kay/\h$. This allows for an easy computation of $\sigma^*\colon \Omega^*(K/H)^K\to \Omega^*(K/H)^K$: since the derivative $\sigma'$ equals $-1$ on $\kay/\h$, $\sigma^*$ takes a $k$-form $\omega$ to $(-1)^k\omega$. Since we also have $d\circ \sigma^* = \sigma^*\circ d$, it follows that $d=0$ on $\Omega^*(K/H)$. We conclude:
\begin{prop}[{\cf\cite[p.\,167, Thm.\,7]{Onishchik1994_tttg}}]
\label{prop:theta_and_formality}
    The cdga $\Omega^*(K/H)^K$ is formal and isomorphic to $H^*(K/H,\R)$. The involution $\sigma^*$ on $H^*(K/H,\R)$ equals $-\id$.
\end{prop}
\begin{cor}
\label{cor:formal_min_model_group}
    The cdga $(H^*(K,\R),0)$ is a minimal model of $K$.
\end{cor}
\begin{proof}
    Indeed, $K$ is itself a compact symmetric space as indicated in \S\ref{subsubsec:cpct_sym_spaces}, hence formal by Proposition \ref{prop:theta_and_formality}. It is moreover an exterior algebra by Theorem \ref{thm:H-group_coh_structure}(i), so the claim follows from Lemma \ref{lem:formal_min_model}.
\end{proof}

\subsubsection{Induced morphisms in cohomology}
\label{subsubsec:prelim_maps_in_coh}
We keep the notation of \ref{subsubsec:cpct_sym_spaces} and moreover write $p\colon K\to K/H$ for the canonical projection and $i\colon H\to K$ for the inclusion. Here we discuss the morphisms $i^*$, $p^*\colon H^*(K/H,\R)\to H^*(K,\R)$ and $C^*\colon H^*(K,\R)\to H^*(K/H,\R)$. To begin with, it is clear that the graded Hopf algebra homomorphism $i^*$ sends the primitive elements $P_K$ to $P_H$. Since $P_K$ generates $H^*(K,\R)$, we have:
\begin{prop}
    The kernel of $i^*\colon H^*(K,\R)\to H^*(H,\R)$ is generated by $\ker(i^*|_{P_K}\colon P_K\to P_H)$.
\end{prop}
Next, one can use the left action of $K$ on $K/H$ (more precisely the resulting comodule structure in cohomology) to show:
\begin{prop}[{dual of \cite[Satz V]{Samelson1941_Beitraege}, \cf\cite[\S9, Prop.\,16]{Onishchik1994_tttg}}]
\label{prop:induced_morphism_p}
    The image of $p^*\colon H^*(K/H,\R)\to H^*(K,\R)$ is generated as a unital subalgebra by its intersection with $P_K$. 
\end{prop}
The two preceding propositions feature subspaces of $P_K$. To describe them, note that the involution $\sigma^*\colon H^*(K,\R)\to H^*(K,\R)$ takes $P_K$ to itself; we thus have an eigenspace decomposition $P_K = P_K^+ \oplus P_K^-$ into $\sigma^*$-eigenspaces for $1$ resp.\,$-1$. By Proposition \ref{prop:theta_and_formality} and the fact that $\sigma$ is trivial on $H$, we have
\begin{equation}
\label{eq:P_K-}
    \im(p^*)\cap P_K \subseteq P_K^- \subseteq \ker(i^*)\cap P_K.
\end{equation}

To say more, we require an additional assumption, which is expressible in several ways:
\begin{thm}
\label{thm:Samelson}
    Suppose $H$ is connected. Then the following conditions are equivalent:
    \begin{itemize}
        \item[(i)] The map $i^*\colon H^*(K,\R)\to H^*(H,\R)$ is surjective.
        \item[(ii)] The fundamental homology class $[H]$ of $H$ is nonzero in $H_*(K,\R)$.
        \item[(iii)] The map $(Bi)^*\colon H^*(BK,\R)\to H^*(BH,\R)$ is surjective.
        \item[(iv)] The restriction map $\R[\kay]^K\to \R[\h]^H$ is surjective.
    \end{itemize}
    If the equivalent conditions (i)--(iv) hold, then $p^*\colon H^*(K/H,\R)\to H^*(K,\R)$ is injective with image $E(P_K^-)$.
\end{thm}
\begin{proof}
    Item (i) is equivalent to injectivity of $i_*\colon H_*(H,\R)\to H_*(K,\R)$. The homology Hopf algebras of these compact Lie groups are dual to their cohomology Hopf algebras and therefore also exterior algebras. Since $[H]$ spans the top degree part of the exterior algebra $H_*(H,\R)$ and $i_*$ respects the algebra structure, $i_*$ is injective if and only if $i_*([H])\ne 0$; thus (i) is equivalent to (ii). The equivalence of (i) with (iii) and (iv) follows from Theorem \ref{thm:transgression_EG} using diagram \eqref{eq:transgression_naturality_diagram}. We now assume that (i)--(iv) hold and prove the remaining assertion, which is a classical theorem of Samelson \cite[Satz VI]{Samelson1941_Beitraege}. Using (i) and the Leray--Serre theorem, we have 
    \begin{equation}
    \label{eq:LeraySerre_K/H}
        H^*(K,\R)\cong H^*(K/H,\R)\otimes H^*(H,\R)    
    \end{equation}
    as $H^*(K/H,\R)$-modules, so in particular $p^*$ is injective. By Proposition \ref{prop:induced_morphism_p}, its image is an exterior subalgebra of $H^*(K,\R) = E(P_K)$. By choosing preimages in $P_K$ of a basis of $P_H$, and then using their products for the preimages of a basis of $H^*(H,\R)$ implicit in the Leray--Serre theorem, we in fact obtain \eqref{eq:LeraySerre_K/H} as an isomorphism of graded algebras. It follows that the inclusions in \eqref{eq:P_K-} must be equalities, and in particular $\im(p^*) = E(P_K^-)$ by Proposition \ref{prop:induced_morphism_p}.
\end{proof}
\begin{defn}
\label{def:nhz}
    Following Samelson \cite{Samelson1941_Beitraege}, we call a connected closed subgroup $H\le K$ \emph{nonhomologous to zero} if the equivalent conditions (i)--(iv) are fulfilled.\footnote{The name of course comes from condition (ii), which we have included for context.}
\end{defn}

A first application is a generalisation of Corollary \ref{cor:formal_min_model_group} (which admits essentially the same proof):
\begin{cor}
\label{cor:nhz_formal}
    If $H\le K$ is nonhomologous to zero then $(H^*(K/H,\R), 0)$ is a minimal model of $K/H$.
\end{cor}
Moreover, in the situation of Theorem \ref{thm:Samelson} we can give the following description of $C^*$:
\begin{prop}
\label{prop:nhz_Cartan}
    Suppose $K/H$ is a symmetric space with $H\le K$ nonhomologous to zero. Then $C^*\colon H^*(K,\R)\to H^*(K/H,\R)$ is surjective, takes $P_K^+$ to zero and $P_K^-$ isomorphically to $p^*(P_K^-)$.
\end{prop}
\begin{proof}
    From the definition of $C$ (see \eqref{eq:Cartan_embedding}) we find
    \begin{equation*}
        p^*\circ C^*
        =
        \mu \circ (\id\otimes (S\circ \sigma^*)) \circ \Delta
    \end{equation*}
    where $\mu$, $S$ and $\Delta$ denote the multiplication, antipode and comultiplication of the graded Hopf algebra $H^*(K,\R)$. The antipode $S$ sends any primitive element $x\in P_K$ to $-x$, so $S\circ \sigma^*$ is equal to $S$ on $P_K^+$ and to the identity on $P_K^-$. Thus,
    \begin{equation*}
        (p^*\circ C^*)(x) = \begin{cases}
            0 & \text{ for } x\in P_K^+\\
            2x & \text{ for } x\in P_K^-
        \end{cases}
    \end{equation*}
    from which the claims about $P_K^+$ follows immediately. Surjectivity and the claim about $P_K^-$ also follow, via Theorem \ref{thm:Samelson}.
\end{proof}

\subsection{Based loop spaces}
\label{subsec:prelim_blp}
Given a connected well-pointed topological space $(X,x)$, its \emph{based loop space} is the space $\Omega X$ of continuous maps $(S^1,1)\to (X,x)$ with the compact-open topology. The usual notions of concatenation and reversal of loops equip $\Omega X$ with the structure of an $H$-group. Thus, provided $\Omega X$ is of finite type, its cohomology is naturally a graded Hopf algebra. 

The based loop space functor is right adjoint to the reduced suspension $\Sigma$, that is, to the wedge product with $S^1$; this in particular yields isomorphisms 
\begin{equation}
\label{eq:loop_space_homotopy_groups}
    \tau_n\colon \pi_n(\Omega X)\overset{\cong}{\rightarrow} \pi_{n+1}(X),
    \quad
    n\ge 0.
\end{equation}
As worked out in \cite[ch.\,X, (7.9) and (7.10)]{WhiteheadEHT}, for a straightforward choice of $\tau_i$ one obtains the following relation between the Samelson product (see \ref{subsubsec:Samelson_product}) of $\Omega X$ and the Whitehead product (see \eqref{eq:Whitehead_product}) of $X$: for $f\in \pi_{m+1}(X)$ and $g\in \pi_{n+1}(X)$, we have
\begin{equation}
\label{eq:Whitehead_Samelson}
    \tau_{m+n}([f,g]_W) = (-1)^m [\tau_m(f), \tau_n(g)]_S.
\end{equation}

There is also the \emph{path space} $PX$ of continuous maps $([0,1], 0)\to (X,x)$, which carries the endpoint map $PX\to X$ sending $\gamma\in PX$ to $\gamma(1)$; this is a fibration with fibre $\Omega X$. We note that $PX$ is contractible. Now, if $X$ is simply connected, then $\Omega X$ is path-connected, so the convergence of the Serre spectral sequence for $\Omega X\to PX \to X$ shows that $\Omega X$ is of finite type provided $X$ is. Thus, in this case $H^*(\Omega X, k)$ and $H^*(\Omega X, k)$ are (dual) connected graded Hopf algebras for any coefficient ring $k$; as before we will assume $k$ to be a field of characteristic zero. Since $H_*(\Omega X, k)$ is cocommutative, Theorems \ref{thm:HA_structure_2} and \ref{thm:Hurewicz_primitives} yield $H_*(\Omega X, k) \cong U(\pi_*(\Omega X)\otimes_\Z k)$, where $\pi_*(\Omega X)\otimes_Z k$ carries the Samelson product. Via \eqref{eq:loop_space_homotopy_groups}, we may further identify $\pi_*(\Omega X)\otimes_\Z k$ with $(\pi_*(X)\otimes_\Z k)[1]$, and \eqref{eq:Whitehead_Samelson} shows that this takes the Samelson product to the Whitehead product modified by a sign. Finally, using Theorem \ref{thm:RHT_main} we arrive at
\begin{equation}
\label{eq:loop_space_Q}
    Q(H^*(\Omega X,k))
    \cong
    P(H_*(\Omega X,k))^*
    \cong 
    ((\pi_*(X)\otimes_\Z k)[1])^*
    \cong
    V_X[1],
\end{equation}
We may conclude the following: 
\begin{prop}
\label{prop:loop_space_cohomology_1}
    Suppose $X$ is simply connected and of finite type. Then $H^*(\Omega X, k)$ and $\Lambda(V_X[1])$ are isomorphic as graded algebras. If the differential of $\Min(X)$ has trivial quadratic part $d_1$ -- for instance, if $X$ is formal with free cohomology algebra -- then they are isomorphic as graded Hopf algebras.  
\end{prop}
\begin{proof}
    The first claim follows directly from Theorem \ref{thm:HA_structure_1} and \eqref{eq:loop_space_Q}. If $d_1=0$, then Theorem \ref{thm:RHT_main}, specifically identity \eqref{eq:d1_and_Whitehead}, shows that all Whitehead products on $X$ vanish. With the above discussion, it follows that all Samelson products on $\Omega X$ vanish, so that the graded Lie algebra $P_*(H_*(\Omega X,k))$ has trivial Lie bracket. Thus, $H^*(\Omega X, k)$ is cocommutative and the claim follows from Theorem \ref{thm:H-group_coh_structure}(ii).
\end{proof}

One can obtain a variant of this result in a manner parallel to Theorem \ref{thm:transgression_EG}; this version was pointed out by Bott \cite[\S7]{Bott1958_loops} in the case where $X$ is a compact Lie group:
\begin{prop}
\label{prop:Bott_SpSeq_variant}
    Suppose $X$ is simply connected and that $H^*(X,k)=\Lambda(V)$ is a finitely generated free graded commutative algebra. Then the primitive elements of $H^*(\Omega X, k)$ coincide with the transgressive elements for the Leray--Serre spectral sequence associated to the fibration $\Omega X\to PX \to X$, and the transgression maps them isomorphically to $V$.
\end{prop}
\begin{proof}
    This comes, via standard arguments (see e.\,g.\,\cite[Thm.\,3.27]{McClearyUGSS} for a variant), from an analysis of the spectral sequence and the fact that $PX$ is contractible. One chooses a homogeneous basis $\{x_1,x_2,\ldots x_n\}$ of $V$ ordered by nondecreasing degree and first proves inductively that each $x_i$ is in the image of the transgression, say with preimage $y_i$. One then combines each pair $(x_i,y_i)$ into an elementary spectral sequence with second page $F(y_i,x_i)$ converging to $k$ (in bidegree $(0,0)$) in which $y_i$ is again transgressive and mapped to $x_i$. The tensor product of these elementary spectral sequences maps into the Leray--Serre spectral sequence and agrees with it on the $\infty$-page as well as the fibre terms on the second page, so Zeeman's comparison theorem (\cf\cite[Thm.\,3.26]{McClearyUGSS}) yields $H^*(\Omega X, k) \cong F(y_1,\ldots, y_n)$ as graded algebras. It is then straightforward to check that the $y_i$ span the spaces of primitive and of transgressive elements, which finishes the proof.  
\end{proof}
    
Lastly, we treat the case where $X$ is path-connected and has finite but possibly nontrivial fundamental group. Then \eqref{eq:loop_space_homotopy_groups} yields $\pi_0(\Omega X)\cong \pi_1(X)$, and the component of the constant loop can be identified with $\Omega \Tilde{X}$ where $\Tilde{X}$ is the universal covering space of $X$. Moreover, for $n\ge 1$, the bijections in \eqref{eq:loop_space_homotopy_groups} identify the action of $\pi_0(\Omega X)$ on $\pi_n(\Omega X)$ with the usual action of $\pi_1(X)$ on $\pi_{n+1}(X) \cong \pi_{n+1}(\Tilde{X})$. We recall that $X$ is called \emph{homotopy-simple} if $\pi_1(X)$ acts trivially on the higher homotopy groups. Via Theorem \ref{thm:Hurewicz_primitives} and Proposition \ref{prop:coh_of_disconnected_H-group}, we obtain:
\begin{prop}
\label{prop:coh_of_disconnected_loop_space}
    If $X$ is homotopy-simple, then there are Hopf algebra isomorphisms
    \begin{equation*}
        H^*(\Omega X,k)
        \cong 
        H^*(\Omega \Tilde{X}, k)\otimes_k k[\pi_1(X)]
        \cong
        U(V_{\Tilde{X}}[1]) \otimes_k k[\pi_1(X)].
    \end{equation*}
\end{prop}

\subsection{Equivariant cohomology}
Let $G$ be a topological group. By a $G$-space we mean a topological space $X$ with a continuous group action $G\times X\to X$. As before, we write $EG\to BG$ for the universal $G$-bundle. The \emph{equivariant cohomology} of a $G$-space $X$ with coefficients in a ring $R$ is defined as
\begin{equation*}
    H^*_G(X, R)\coloneqq H^*(X_G, R),
\end{equation*}
where $X_G\coloneqq X\times^G EG = (X\times EG)/G$. The definition of $X_G$ yields a canonical fibration $X_G\to BG$ with fibre $X$. In particular we can (and always will) view $H^*_G(X,R)$ as a graded algebra over $H^*(BG,R)$; in turn, $H^*(BG,R)$ can be viewed as $H^*_G(\pt,R)$, the equivariant cohomology of a singleton space with trivial $G$-action. Moreover, the fibre inclusion $X\to X_G$ yields (after choosing a base point) a canonical homomorphism $H^*_G(X,R)\to H^*(X,R)$. 

For simplicity, we take $R$ to be a field $k$ for the remainder of this subsection. We say that the $G$-space $X$ is \emph{equivariantly formal} (over $k$) if the canonical homomorphism $H^*_G(X,k)\to H^*(X,k)$ is surjective. Equivalently, $X$ is equivariantly formal if and only if the Leray--Serre spectral sequence associated to the fibration $X\to X_G\to BG$ collapses on the second page. This second characterisation immediately yields the following facts:
\begin{prop}
\label{prop:eq_formality}
    Assume that $G$ is connected.
    \begin{itemize}
        \item[(i)] If $X$ is equivariantly formal, then $H^*_G(X,k)\cong H^*_G(\pt, k)\otimes_k H^*(X,k)$ (as graded modules).
        \item[(ii)] If both $H^*(X,k)$ and $H^*(BG,k)$ are concentrated in even degrees, then $X$ is equivariantly formal.
    \end{itemize}
\end{prop}
If $X$ is equivariantly formal, the first characterisation entails that $H^*(X,k)$ is a quotient of the equivariant cohomology, and more precisely one has:
\begin{prop}
\label{prop:eq_formal_sp_fibre}
    If $X$ is equivariantly formal, then $H^*(X,k) = H_G^*(X,k)/(H_G^+(\pt,k))$. Here $H_G^+(\pt,k)$ denotes the ideal generated by the image of the positively graded part of $H_G^*(\pt,k)$ in $H^*_G(X,k)$.
\end{prop}

Equivariant cohomology is functorial with respect to spaces with group action. More precisely, there is a category whose objects are pairs $(G,X)$ where $G$ is a topological group with a continuous group action $a_X\colon G\times X\to X$, and a morphism $(G,X)\to (H,Y)$ consists of a group homomorphism $\varphi\colon G\to H$ and a continuous map $f\colon X\to Y$ such that $f\circ a_X = a_Y\circ (\varphi\times f)$. Equivariant cohomology (with coefficients in $k$) is then a contravariant functor from this category to the category of algebras of graded $k$-algebras, whose objects are graded $k$-algebra homomorphisms and whose morphisms are commutative squares between them. It assigns $(G,X)\mapsto H^*_G(X,k)$ and takes $(\varphi, f)$ to the homomorphism $H^*_H(Y,k)\to H^*_G(X,k)$ induced by
\begin{equation*}
    f\times^\varphi E\varphi\colon 
    X\times^G EG\to Y\times^H EH,
    \quad
    [(x, u)]\mapsto [(f(x), E\varphi(u))],
\end{equation*}
lying over the homomorphism $H^*(BH,k)\to H^*(BG,k)$ induced by $B\varphi$. We thus have the following commutative diagram whose left square depicts the image of the morphism $(\varphi,f)$:
\begin{equation}
\label{eq:eq_coh_funct}
\begin{tikzcd}
    H^*_H(\pt,k)
        \arrow[r]
        \arrow[d, "(B\varphi)^*"]
    & H^*_H(Y,k)
        \arrow[r]
        \arrow[d, "(f\times^\varphi E\varphi)^*"]
    & H^*(Y,k)
        \arrow[d, "f^*"]
    \\
    H^*_G(\pt,k)
        \arrow[r]
    & H^*_G(X,k)
        \arrow[r]
    & H^*(X,k)
\end{tikzcd}
\end{equation}

In this context, the following lemma will be useful:
\begin{lem}
\label{lem:eq_coh_restriction_surjective}
    Let $X$, $Y$, $G$, $H$, $\varphi$ and $f$ be as above. Suppose that $X$ and $Y$ are equivariantly formal and both $(B\phi)^*\colon H^*(BH,k)\to H^*$ and $f^*\colon H^*(Y,k)\to H^*(X,k)$ are surjective. Then the homomorphism $H^*_H(Y,k)\to H^*_G(X,k)$ induced by $(\varphi,f)$ is also surjective.
\end{lem}
\begin{proof}
    Under the stated assumptions, the induced morphism of Leray--Serre spectral sequences is surjective on the second page. By equivariant formality, this is sufficient.
\end{proof}

Lastly, we recall that if $G$ is a compact or complex Lie group with maximal torus $T$ and Weyl group $W$, then $W$ acts naturally on $H^*_T(X,k)$, and we have $H^*_G(X,k)\cong H^*_T(X,k)^W$. (See \eg\cite[ch.III, p.\,38, Prop.\,1]{Hsiang_CToTTG} for a proof.)

\subsection{Complex and real reductive groups}
By a \emph{complex reductive group}, we mean a linear algebraic group $G$ over $\C$ whose unipotent radical is trivial. Similarly, a \emph{real reductive group} is a linear algebraic group $G$ over $\R$ with trivial unipotent radical; this is equivalent to $G_\C$ being a complex reductive group. We may define $G$ to be \emph{semisimple} if its Lie algebra $\g$ is semisimple, and \emph{almost simple} if $\g$ is simple. When considering the complex or (where applicable) real points of $G$ as above, we always take them with their natural Lie group structure. For our purposes, the groups in question should be (Zariski-)connected, so we will usually assume this; note however that the Lie group $G(\R)$ of a connected real reductive group $G$ need not be connected.

If $G$ is any (complex or real) reductive group, with (complex or real) Lie algebra $\g$, we write $G_\mathrm{ad}\le \GL(\g)$ for the image of its adjoint representation. We have $G_\mathrm{ad} \cong G/Z(G)$ where $Z(G)$ is the centre of $G$, and $G_\mathrm{ad}$ is semisimple. We say that $G$ is of \emph{adjoint type} if $G\cong G_\mathrm{ad}$. The \emph{derived subgroup} of $G$ is denoted $G_\mathrm{der}\coloneqq [G,G]$, and its Lie algebra $\g_\mathrm{der}$; we recall also that $\g = \g_\mathrm{der} \oplus Z(\g)$. An \emph{isogeny} of connected complex reductive groups is a surjective morphism with finite kernel; equivalently, $\phi\colon G\to H$ is an isogeny if and only if its derivative $\phi'\colon \operatorname{Lie}(G)\to \operatorname{Lie}(H)$ is a Lie algebra isomorphism. 

Below, we recall further standard definitions specific to the complex or real case, as well as facts about automorphisms and real structures.

\subsubsection{Setup for connected complex reductive groups}
\label{subsubsec:cplx_red_setup}
Here we largely follow \cite{Springer1979_RG}. Let $G$ be a connected complex reductive group. Choose a maximal torus $T\le G$. The \emph{root datum} of $(G,T)$ is the quadruple $\Psi = \Psi_{G,T} = (\Lambda, \Phi, \Lambda^\vee, \Phi^\vee)$ where $\Lambda = \Hom(T, \G_m)$ is the \emph{weight lattice}, $\Lambda^\vee \coloneqq \Hom(\G_m,T)$ the \emph{coweight lattice}, $\Phi \subseteq \Lambda$ the set of \emph{roots}, that is, the nontrivial characters of $T$ for the adjoint representation, and $\Phi^\vee \subseteq \Lambda^\vee$ the set of coroots. Composition defines a perfect pairing $\langle \cdot, \cdot \rangle \colon \Lambda\times \Lambda^\vee \to \Aut(\G_m) \cong \Z$. The set of roots $\Phi$ has the structure of a root system (in $\Phi\otimes_\Z \R$) whose Weyl group $W$ is canonically isomorphic to $N_G(T)/T$. Moreover, to each root $\alpha$ there corresponds a connected, one-dimensional, unipotent subgroup $U_\alpha \le G$ called its \emph{root group}.

We further define the \emph{root lattice} $Q\coloneqq \Z\Phi\subseteq \Lambda$ and \emph{coroot lattice} $Q^\vee \coloneqq \Z\Phi^\vee \subseteq \Lambda^\vee$. The quotient group $\Lambda^\vee/Q^\vee$ is known to be isomorphic to the fundamental group $\pi_1(G) = \pi_1(G(\C))$. In particular, $G$ is simply connected if and only if $Q^\vee = \Lambda^\vee$; note that this implies that $G$ is semisimple (which is equivalent to $\operatorname{rank}(Q^\vee) = \operatorname{rank}(\Lambda^\vee)$).

The notions just defined depend on the choice of $T$, and we will usually suppress the dependence and work with some fixed choice. However, it is sometimes useful to avoid choices, which can be done because all maximal tori of $G$ are conjugate. We may therefore define an abstract version $\Psi_{G,\mathrm{abs}}$ of the root datum, as the (co)limit of all $\Psi_{G,T}$ along the isomorphisms $\Psi_{G,T}\to \Psi_{G,T'}$ induced by conjugations. (Note that all conjugations sending $T$ to $T'$ induce the same isomorphism.)

There is a category of (abstract, reduced) root data such that the assignment $(G,T)\mapsto \Psi_{G,T}$ defines a contravariant functor from the category of pairs $(G,T)$ with isogenies to that category, and this turns out to be an equivalence of categories. Thus, for any root datum $\Psi$ whose associated root system is reduced, there exist a connected complex reductive group $G$ and a maximal torus $T\le G$, unique up to isomorphism (of pairs), such that $\Psi \cong \Psi_{G,T}$. (In choice-free terms, $\Psi\cong \Psi_{G,\mathrm{abs}}$.)

For $\Psi_{G,T}$ as above, one can check that $\Check{\Psi} \coloneqq (\Lambda^\vee, \Phi^\vee, \Lambda, \Phi)$ also fulfils the root datum axioms, so this yields a pair $(\Check{G}, \Check{T})$, unique up to isomorphism, such that $\Psi_{\Check{G}, \Check{T}} \cong \Check{\Psi}$. The connected complex reductive group $\Check{G}$ is called the \emph{Langlands dual} of $G$. By definition, the coweights (coroots, weights, roots) of $\Check{G}$ are identified with the weights (roots, coweights, coroots) of $G$. 

We now further choose a Borel subgroup $B$ of $G$ containing $T$. This gives rise to the following sets: the \emph{simple roots} $\Pi \subseteq \Phi$, \emph{positive roots} $\Phi_+ \subseteq \Phi$, \emph{simple coroots} $\Pi^\vee\subseteq \Phi^\vee$, \emph{positive coroots} $\Phi^\vee_+\subseteq \Phi^\vee$, \emph{dominant weights} $\Lambda_+\subseteq \Lambda$ and \emph{dominant coweights} $\Lambda^\vee_+\subseteq \Lambda^\vee$. Conversely, $B$ is recovered from $\Pi$ as the subgroup generated by $T$ and the root groups of simple roots. Thus, the choice of $B\ge T$ is equivalent to a choice of simple roots, and in turn also to such choices for $\Check{G}$. In more abstract terms, this refines the root datum $\Psi$ to a so-called \emph{based root datum}, which contains also $\Pi$ and $\Pi^\vee$. One can again take a colimit of the based root data for all $B,T$ and define an abstract version of the based root datum.

The irreducible finite-dimensional (complex) representations of $G$ are in bijection with $\Lambda_+$; the representation corresponding to $\lambda\in \Lambda_+$ is denoted $V^\lambda$ and called the \emph{representation of highest weight $\lambda$}. This is justified by the fact that $\lambda$ indeed arises as a $T$-weight of $V^\lambda$, and it is maximal in the set $\wt(\lambda)$ of such $T$-weights. The maximality here is with respect to the partial ordering on $\Lambda$, defined by setting $\mu \le \lambda$ if and only if $\lambda-\mu$ is a sum of simple roots. There is an analogous partial ordering $\le$ on $\Lambda^\vee$ defined in terms of simple coroots.

\subsubsection{Automorphisms}
\label{subsubsec:autos_setup}
Let $G$ be a connected complex reductive group, $\g$ its Lie algebra. We denote its automorphism group by $\Aut(G)$, and the normal subgroup of inner automorphisms (those given by conjugation by group elements) by $\Int(G)$. The quotient group is $\Out(G)\coloneqq \Aut(G)/\Int(G)$, so the following is a short exact sequence by definition:
\begin{equation}
\label{eq:aut_SES}
    1 \rightarrow 
    \Int(G) \rightarrow
    \Aut(G) \rightarrow
    \Out(G) \rightarrow 1.
\end{equation}

Let $T\le G$ be a maximal torus and $\Psi = (\Lambda,\Phi, \Lambda^\vee, \Phi^\vee)$ the corresponding root datum. If $\sigma\in \Aut(G)$ maps $T$ to itself, then it defines an automorphism of $\Psi$, which we also denote by $\sigma$. Explicitly, the action on weights is $\sigma(\lambda) \coloneqq \lambda\circ \sigma\inv$ for $\lambda\in \Lambda = \Hom(T,\G_m)$, and the action on coweights is $\sigma(\mu) \coloneqq \sigma \circ \mu$ for $\mu\in \Lambda^\vee = \Hom(\G_m,T)$. A refinement of this is as follows:
\begin{defn}\leavevmode
\begin{itemize}
    \item[(i)]
    A \emph{pinning} of $G$ is a tuple $(B, T, \{u_\alpha\}_{\alpha \in \Pi})$, where $B\le G$ is a Borel subgroup, $T\le B$ a maximal torus, $\Pi$ denotes the set of simple roots of $G$ with respect to $(B,T)$, and, for each $\alpha\in \Pi$, $u_\alpha\colon \G_a\to U_\alpha$ is an isomorphism.
    \item[(ii)]
    An automorphism $\sigma\in \Aut(G)$ is said to \emph{preserve} a pinning $(B,T,\{u_\alpha\}_{\alpha \in \Pi})$ of $G$ if $\sigma(B) = B$, $\sigma(T) = T$, and postcomposition by $\sigma$ permutes the $u_\alpha$. 
    \item[(iii)] 
    An automorphism $\sigma\in \Aut(G)$ is called \emph{pinning-preserving} if it preserves some pinning of $G$.
\end{itemize}
\end{defn}
\begin{rmk}\leavevmode
\begin{enumerate}
    \item 
    The Borel subgroup of a pinning is redundant in the sense that it is generated by $T$ and the images of the $u_\alpha$.
    \item 
    There is an analogous notion of pinning for the Lie algebra $\g$ -- a choice of Borel subalgebra, Cartan subalgebra and nonzero simple root vectors -- such that pinnings of $G$ are in bijection with pinnings of $\g$. An automorphism of $G$ then preserves a pinning of $G$ if and only if its derivative preserves a pinning of $\g$.
\end{enumerate}
     
\end{rmk}
\begin{prop}[{\cf\cite[\S2]{Springer1979_RG}}]
\label{prop:aut_SES_splitting}
    Let $(B,T,\{u_\alpha\}_{\alpha \in \Pi})$ be a pinning of $G$. For any $\sigma\in \Aut(G)$, there exists a unique $\phi\in \Int(G)$ such that $\sigma\phi$ preserves $(B,T,\{u_\alpha\}_{\alpha \in \Pi})$. The assignment of $[\sigma]\in \Out(G)$ to $\sigma\phi$ defines a splitting of \eqref{eq:aut_SES} which identifies $\Out(G)$ with the group of automorphisms of $G$ which preserve $(B,T,\{u_\alpha\}_{\alpha \in \Pi})$.
\end{prop}
\begin{rmk}
\label{rmk:Out_LD}
    In the same way, we may identify $\Out(G)$ with the automorphism group of the abstract based root datum. This yields in particular a canonical isomorphism $\Out(G)\cong \Out(\Check{G})$, which can be made concrete by choosing compatible pinnings of $G$ and $\Check{G}$.
\end{rmk}

The choice of pinning in Proposition \ref{prop:aut_SES_splitting} is not important as all pinnings are conjugate in the following sense:
\begin{prop}[{\cf\cite[Exp.\,XXIV, Lem.\,1.5]{SGA3_III}}]
\label{prop:pinnings_conjugate}
    Let $(B,T,\{u_\alpha\}_{\alpha \in \Pi})$ and $(B',T',\{u'_\beta\}_{\beta \in \Pi'})$ be pinnings of $G$. Then there exists a unique $\phi \in \Int(G)$ such that $\phi(B) = B'$, $\phi(T) = T'$, precomposition with $\phi$ takes $\Pi'$ to $\Pi$, and $u_{\beta} = \phi\circ u_{\beta\circ \phi}$ for all $\beta\in \Pi'$.
\end{prop}
\begin{cor}
\label{cor:pi_pr_auto_uniqueness}
    Suppose $\sigma_1,\sigma_2\in \Aut(G)$ are pinning-preserving automorphisms such that $\sigma_2 \in \sigma_1 \Int(G)$. Then there exists $\phi\in \Int(G)$ such that $\sigma_2 = \phi \sigma_1 \phi\inv$.
\end{cor}
\begin{proof}
    By Proposition \ref{prop:pinnings_conjugate}, there exists $\phi\in \Int(G)$ such that $\sigma_2$ and $\phi \sigma_1 \phi\inv$ preserve the same pinning. If $\phi$ is conjugation by $g\in G$, then $\phi\sigma_1\phi\inv = \sigma \psi$, where $\psi \in \Int(G)$ is conjugation by $\sigma(g)g\inv$. Thus, $\sigma_2 \in \sigma_1\Int(G) = (\phi\sigma_1\phi\inv)\Int(G)$, so Proposition \ref{prop:aut_SES_splitting} implies $\sigma_2 = \phi\sigma_1\phi\inv$.    
\end{proof}

\begin{rmk}[Diagram automorphisms]
\label{rmk:diagram_autos}
    Since $G$ is assumed connected, there is a natural inclusion of $\Aut(G)$ into the automorphism group $\Aut(\g)$ of the Lie algebra. It is an isomorphism if $G$ is simply connected or of adjoint type, but need not be surjective in general. A discussion of this discrepancy an be found in \cite[\S1.5]{ConradReductiveGroupSchemes}. In any case, one can describe all pinning-preserving automorphisms of $G$ by first describing all pinning-preserving automorphisms of $\g$ and, if necessary, studying which of them lift to $\Aut(G)$. If $\g$ is semisimple, the automorphisms of $\g$ preserving a fixed pinning are in bijection with automorphisms of the Dynkin diagram (as follows from the group case above, but see also \cite[\S4]{Onishchik2004_real}), which provides a convenient way of classifying them.
\end{rmk}

\subsubsection{Real structures and Cartan involutions}
Here we mainly follow \cite{AT_2018_Galois_and_Cartan}; another useful source is \cite[Sec.\,2]{GPRa2019_involutions}. For this subsection, we fix a complex reductive group $G$. A \emph{real form} of $G$ is then a real reductive group $G_0$ such that $(G_0)_\C \cong G$. A \emph{real structure} of $G$ is a multiplicative, involutive automorphism of the $\R$-scheme $G$ lying over the complex conjugation automorphism of the $\R$-scheme $\Spec \C$. By Galois descent, isomorphism classes of real forms are in a natural bijection with real structures up to conjugation by automorphisms of $G$. We denote the real form corresponding to a real structure $\sigma$ by $G^\sigma$. It is often convenient to consider two further equivalence relations on real structures.\footnote{There does not seem to be one standard terminology for all of these.}

\begin{defn}
\label{def:real_str_eqrels}
    Let $\sigma_1, \sigma_2$ be real structures of $G$.
    \begin{itemize}
        \item[(i)] We call $\sigma_1$ \emph{isomorphic} to $\sigma_2$ if their corresponding real forms are isomorphic, that is, if there exists $\phi \in \Aut(G)$ such that $\sigma_1 = \phi \sigma_2 \phi\inv$.
        \item[(ii)] We call $\sigma_1$ \emph{inner-isomorphic} to $\sigma_2$ if there exists $\phi \in \Int(G)$ such that $\sigma_1 = \phi \sigma_2 \phi\inv$.
        \item[(iii)] We call $\sigma_1$ \emph{inner to} $\sigma_2$ if they have the same image in $\Out(G)$, that is, if there exists $\phi\in \Int(G)$ such that $\sigma_1 = \phi \sigma_2$. The resulting equivalence classes are called \emph{inner classes}.
    \end{itemize}
\end{defn}
We note that all of these are equivalence relations, and that the second refines both the first and the third. In general, $G$ will have several nonisomorphic real forms; the following two kinds are particularly important:
\begin{defn}
    Let $G_0$ be a real form of $G$ and $\sigma$ a real structure with $G^\sigma \cong G_0$. 
    \begin{itemize}
        \item[(i)] We call $G_0$ or $\sigma$ \emph{split} if $G_0$ has a split maximal torus.
        \item[(ii)] We call $G_0$ or $\sigma$ \emph{anisotropic} if it contains no nontrivial split torus.
    \end{itemize}
\end{defn}

For a general field extension $k\to k'$, the usual approach to $k$-forms of a $k'$-group is via split $k$-forms. The real--complex case is set apart by the closely parallel theory of complex Lie groups and their real forms. This gives a particular importance to anisotropic forms, in light of the following:
\begin{prop}[{cf.\,\cite[24.6(c)(ii)]{Borel_LAG}}]
    A real form $G_0$ of $G$ is anisotropic if and only if $G_0(\R)$ is compact in the classical topology.
\end{prop}

\begin{prop}[{cf.\,\cite[\S3]{AT_2018_Galois_and_Cartan} and references therein}]
\label{prop:cpct_real_forms}
    Any maximal compact subgroup of $G(\C)$ is the group of fixed points of an antiholomorphic involutive automorphism. Moreover, all maximal compact subgroups of $G(\C)$ are conjugate.
\end{prop}

\begin{cor}
    Any complex reductive group has a unique inner-isomorphism class of anisotropic real structures.
\end{cor}

\begin{defn}
    Let $\sigma$ be a real structure of $G$. A \emph{Cartan involution} of $\sigma$ is an involution $\theta\in \Aut(G)$ such that $\theta \sigma = \sigma \theta$ and such that the real structure $\theta\sigma$ is anisotropic.
\end{defn}
Such involutions (and more generally automorphisms) admit parallel equivalence relations to those in Definition \ref{def:real_str_eqrels}:
\begin{defn}
    Let $\theta_1, \theta_2\in \Aut(G)$.
    \begin{itemize}
        \item[(i)] We call $\theta_1$ \emph{isomorphic} to $\theta_2$ if there exists $\phi \in \Aut(G)$ such that $\theta_1 = \phi \theta_2 \phi\inv$.  
        \item[(ii)] We call $\theta_1$ \emph{inner-isomorphic} to $\theta_2$ if there exists $\phi \in \Int(G)$ such that $\theta_1 = \phi \theta_2 \phi\inv$.
        \item[(iii)] We call $\theta_1$ \emph{inner to} $\theta_2$ there exists $\phi\in \Int(G)$ such that $\theta_1 = \phi \theta_2$.
    \end{itemize} 
\end{defn}

\begin{thm}[{cf.\,\cite[Thm.\,3.13 and Cor.\,3.17]{AT_2018_Galois_and_Cartan}}]
\label{thm:Cartan_bijection}
    Every real structure of $G$ admits a Cartan involution. Sending a real structure to its Cartan involution yields a bijection
    between equivalence classes of real structures and equivalence classes of involutive automorphisms for each of the equivalence relations \emph{``isomorphic''}, \emph{``inner-isomorphic''} and \emph{``inner to''}. 
\end{thm}

\begin{rmk}[Choice of Cartan involution]
    A given real structure $\sigma$ generally admits several Cartan involutions, but a key part of Theorem \ref{thm:Cartan_bijection} is that all Cartan involutions of $\sigma$ are inner-isomorphic. We allow ourselves to speak of ``the'' Cartan involution of $\sigma$ whenever it is clear from Theorem \ref{thm:Cartan_bijection} that the precise choice does not matter. On the other hand, the flexibility is also an advantage. For instance, given a maximal torus $T$ preserved by $\sigma$, it is useful to know that $\theta$ can be chosen so that it also preserves $T$. Such facts are perhaps most easily seen by inspecting the construction of Cartan involutions underlying Theorem \ref{thm:Cartan_bijection}. 
\end{rmk}

An advantage of the Cartan involution approach is the following topological fact:
\begin{prop}
\label{prop:max_cpct_htpy}
    Let $\sigma$ be a real structure of $G$ with Cartan involution $\theta$, and $\tau \coloneqq \theta\sigma$. Then $G(\C)^\tau$ and $G(\R)^\tau = G(\R)^\theta$ are maximal compact subgroups of the Lie groups $G(\C)$, $G(\R)$, respectively, and the inclusions $G(\C)^{\tau}\to G(\C)$ and $G(\R)^\tau\to G(\R)$ are homotopy equivalences.  
\end{prop}
\begin{proof}
    The first claim holds for $\C$-points by Proposition \ref{prop:cpct_real_forms}; for $\R$-points it follows from \cite[Lem.\,2.2]{Mostow1955_Selfadjoint}. The second claim follows from the existence of Cartan decompositions -- see \cite[Lem.\,2.1 and 2.2]{Mostow1955_Selfadjoint}.
\end{proof}

\subsubsection{The compact symmetric space corresponding to a real structure}
\label{subsubsec:corresponding_sym_sp}
Let $G$ be a complex connected reductive group with a real structure $\sigma$. If $\theta$ is a Cartan involution for $\sigma$, then $\tau \coloneqq \theta\sigma$ is by definition anisotropic. Thus, $K\coloneqq G(\C)^\tau$ is a connected compact Lie group. The involution $\theta$ preserves $K$ (and agrees there with $\sigma$), and we put $H\coloneqq K^\theta = K^\sigma$. Then $K/H$ is a connected compact symmetric space (see \ref{subsubsec:cpct_sym_spaces}). One can show that a different choice of Cartan involution $\theta$ leads to a symmetric space isomorphic to $K/H$, so this space is naturally associated to the real structure.

\subsubsection{Setup for connected real reductive groups}
\label{subsubsec:real_red_setup}
Let $G$ be a connected real reductive group, $\g$ its Lie algebra, and $\sigma$ the corresponding real structure on $G_\C$. We choose a \emph{maximally split maximal torus} $T\le G$, that is, a maximal torus whose maximal $\R$-split subtorus $S$ is a maximal $\R$-split torus of $G$. We note that all maximally split maximal tori of $G$ are conjugate (through $G(\R)$). 

The real structure $\sigma$ acts on the root datum $(\Lambda,\Phi, \Lambda^\vee, \Phi^\vee)$ of $(G_\C, T_\C)$ as follows: for $\lambda\in \Lambda = \Hom(T_\C, \G_{m,\C})$, we set $\sigma(\lambda) \coloneqq c\circ \lambda \circ \sigma$, where $c\colon \G_{m,\C}\to \G_{m,\C}$ is complex conjugation. For $\mu\in \Lambda^\vee = \Hom(\G_{m,\C}, T_\C)$, we define $\sigma(\mu) \coloneqq \sigma\circ \mu \circ c$. These involutions of $\Lambda$ and $\Lambda^\vee$ preserve $\Phi$ and $\Phi^\vee$, respectively. If $\theta$ is a Cartan involution for $\sigma$ which preserves $T_\C$, then the anisotropic real structure $\theta\sigma$ acts on the root datum as $-1$, so the actions of $\sigma$ and $\theta$ differ by a sign. 

\begin{rmk}
\label{rmk:real_str_out}
    The similarly defined action of a real structure on the abstract root datum (or abstract based root datum) depends only on its inner class. Similarly to Remark \ref{rmk:Out_LD}, this allows us to match inner classes of real structures for Langlands-dual complex reductive groups, according to their actions on the root datum.
\end{rmk}

A root $\alpha\in \Phi$ is called \emph{real} if $\sigma(\alpha) = \alpha$, \emph{imaginary} if $\sigma(\alpha) = - \alpha$ (equivalently $\theta(\alpha) = \alpha$), and \emph{complex} if it neither real nor imaginary. Moreover, an imaginary root $\alpha$ is called \emph{compact} if $\theta|_{\g_\alpha} = \id_{\g_\alpha}$ and \emph{noncompact} if $\theta|_{\g_\alpha} = -\id_{\g_\alpha}$.

The restriction to $S$ of a root $\alpha \in \Phi$ is trivial precisely when $\alpha$ is imaginary. The nontrivial restrictions of roots to $S$ form the set of \emph{relative roots} $\Phi_\R \subseteq \Hom(S,\G_m)$, which is a (possibly nonreduced) root system. There is a root space decomposition $\g = \m \oplus \bigoplus_{\alpha\in \Phi_\R}\g_\alpha$, where $\g_\alpha$ denotes the weight space of $\alpha$ and $\m$ is the centraliser of $S$ in $\g$. We write $n_\alpha \coloneqq \dim_\R \g_\alpha$ for the \emph{multiplicity} of $\alpha$ (which need not be one). We also define \emph{relative coweights} $\Lambda^\vee_\R \coloneqq \Hom(\G_m,S)$ and note that (via base change and the inclusion $S_\C\to T_\C$) we may identify $\Lambda^\vee_\R$ with $(\Lambda^\vee)^\sigma$. 

We now further choose a minimal (real) parabolic subgroup $P\le G$ containing $T$, and a Borel subgroup $B$ of $G_\C$ such that $T_\C\le B \le P_\C$. Using $B$ to define positive roots $\Phi_+\subseteq \Phi$ etc., the nontrivial restrictions of positive roots to $S$ define positive relative roots $\Phi_{\R,+} \subseteq \Phi_\R$. We may also define \emph{dominant relative coweights} $\Lambda^\vee_{\R,+} \subseteq \Lambda^\vee_\R$ as those whose pairings with any positive relative root is nonnegative; equivalently, $\Lambda^\vee_{\R,+} = (\Lambda^\vee_+)^\sigma$. Although there is also a notion of relative coroots, we shall use the partial ordering $\le$ on $\Lambda^\vee_\R$ obtained by restricting the partial ordering of $\Lambda^\vee$.

\subsection{Affine Grassmannians}
\label{subsec:RealGr}
Here we follow \cite{Zhu2016_introduction} and \cite{Mitchell1988_ThmOnBuildings}. For a reductive linear algebraic group over a field $k$ (or in fact more generally), we may form the \emph{(formal) loop group} $LG$, the \emph{positive loop group} $L^+G$ and the \emph{affine Grassmannian} $\Gr = \Gr_G$. The former two are the sheaves given by $LG(R) \coloneqq G(R((t)))$ and $L^+G(R) \coloneqq G(R[[t]])$ for $k$-algebras $R$, and are representable by (and identified with) an ind-scheme and a scheme over $k$, respectively. The \emph{affine Grassmannian} $\Gr = \Gr_G$ is an ind-projective ind-scheme of ind-finite type over $k$ representing the (étale or equivalently fpqc) quotient $LG/L^+G$. 

Below, we will take $k$ to be either $\R$ or $\C$ and consider the resulting spaces of $k$-points in the classical topology (that is, the direct limit of classical topologies obtained from the ind-scheme structure). In fact, since any complex reductive group $G$ has a real form $G_0$, and since $\Gr_{G_0}(\C) \cong \Gr_{(G_0)_\C}(\C) \cong \Gr_G(\C)$, the setting of a real reductive group suffices to treat both cases.\footnote{One could further identify (as topological spaces) the complex points of $\Gr_G$ with the real points of $\Gr_{G_\R}$, where $G_\R$ is the Weil restriction to $\R$, but this appears more confusing than helpful here.} It is worth pointing out that the real structure $\sigma$ of $G$ corresponding to $G_0$ naturally defines a real structure of $\Gr_G$ (with real form $\Gr_{G_0}$) such that $\Gr_{G_0}(\R) = \Gr_G(\C)^\sigma$.

\subsubsection{Real and complex points}
\label{subsubsec:Gr_points}
Let $G$ be a connected real reductive group, $\Gr = \Gr_G$ as above, and use the setup of \S\ref{subsubsec:real_red_setup}. 
The inclusions $\R[t,t\inv]\to \R((t))$ and $\C[t,t\inv]\to \C((t))$ induce compatible inclusions $\Lambda^\vee_\R\to LG(\R)$ and $\Lambda^\vee\to LG(\C)$. We denote the image of $\lambda \in \Lambda^\vee$ by $t^\lambda$. One has the double coset decompositions (see \cite[Cor.\,2.7]{IM1965_Bruhat} for the second, \cite[Prop.\,4.4.3]{BruhatTits_1} for both)
\begin{equation}
\label{eq:double_coset_decompositions}
    LG(\R) = \bigsqcup_{\lambda\in \Lambda^\vee_{\R,+}} L^+G(\R)t^\lambda L^+G(\R),\quad
    LG(\C) = \bigsqcup_{\lambda\in \Lambda^\vee_+} L^+G(\C)t^\lambda L^+G(\C).
\end{equation}
For $\lambda\in \Lambda^\vee$, denote now the induced real or complex point of $\Gr$ by $\underline{t}^\lambda$. We define the corresponding \emph{affine Schubert cell} $\Gr^\lambda$ as follows: first, note that there is an action of $L^+G$ on $\Gr$ by left multiplication. If $\lambda\in \Lambda^\vee_\R$, then we put $\Gr^\lambda\coloneqq L^+G\cdot \underline{t}^\lambda$; otherwise $\Gr^\lambda \coloneqq (L^+G)_\C\cdot \underline{t}^\lambda$. The affine Schubert cells $\Gr^\lambda$ are finite-dimensional smooth (locally closed) subschemes of $\Gr$ or $\Gr_\C$, respectively. Moreover, if $\lambda$ and $\lambda'$ are in the same Weyl group orbit, then $\Gr^\lambda = \Gr^{\lambda'}$, so that one usually assumes $\lambda\in \Lambda^\vee_+$. From \eqref{eq:double_coset_decompositions} one it follows at once that
\begin{equation}
\label{eq:Schubert_cell_partitions}
    \Gr(\R) = \bigsqcup_{\lambda\in \Lambda^\vee_{\R,+}} \Gr^\lambda(\R), \quad
    \Gr(\C) = \bigsqcup_{\lambda\in \Lambda^\vee_+} \Gr^\lambda(\C).
\end{equation}
The closures $\Gr^{\le \lambda}\coloneqq \overline{\Gr^\lambda}$ are finite-dimensional projective varieties called \emph{affine Schubert varieties}. The notation is explained by the fact that
\begin{equation}
\label{eq:aSc_closures}
    \Gr^{\le \lambda}(\C) = \bigsqcup_{\mu\in \Lambda^\vee_+,\, \mu \le \lambda} \Gr^\mu(\C), \text{ for } \lambda\in \Lambda^\vee_+
    \quad
    \text{and}
    \quad
    \Gr^{\le \lambda}(\R) = \bigsqcup_{\mu\in \Lambda^\vee_{\R,+},\, \mu \le \lambda} \Gr^\mu(\R) \text{ for } \lambda\in \Lambda^\vee_{\R,+}.
\end{equation}
(For the latter, see \cite[Prop.\,3.6.1]{Nadler2005_PSoRLG}.)

The positive loop groups $L^+G$ and $L^+G_\C$ have homomorphisms $\ev_0$ to $G$ and $G_\C$, respectively, defined by sending $t$ to $0$. The orbits of the subgroups $\Tilde{P}\coloneqq \ev_0\inv(P) \le L^+G$ and $\Tilde{B}\coloneqq \ev_0\inv(B)\le L^+G_\C$ refine the decompositions into affine Schubert cells; we call those orbits \emph{Iwahori cells}. A full set of representatives of the orbits is given by $\Lambda^\vee_\R$ in the real and $\Lambda^\vee$ in the complex case, thus, an affine Schubert cell $\Gr^\lambda$ is the union of the Iwahori cells $C_\mu$ (through $t^\mu$) with $\mu$ in the same Weyl orbit as $\lambda$. Moreover, the Iwahori cells form an affine paving, in fact a CW decomposition \cite{Mitchell1988_ThmOnBuildings}. In the real case, we have 
\begin{equation}
\label{eq:Iwahori_cell_dimension}
    \dim(C_\mu) = \sum_{\alpha\in \Phi_{\R,+}, \langle \mu,\alpha \rangle \ge 0} n_\alpha \langle \mu,\alpha \rangle + \sum_{\alpha\in \Phi_{\R,+}, \langle \mu,\alpha \rangle < 0}  n_\alpha (- \langle \mu,\alpha \rangle -1).
\end{equation}
The same formula applies to complex Iwahori cells if we either replace $\Phi_\R$ by $\Phi$ (with all multiplicities $n_\alpha$ equal to one) and view the result as complex dimension, or if we view $G_\C$ as a real reductive group (with all $n_\alpha = 2$). In any case, all Iwahori cells of a complex affine Grassmannian have even real dimension, hence yield a basis of the (co)homology. In the real case, the situation is generally more complicated (see \eg\cite{Kocherlakota1995_IntegralHomology}) but simplifies in several of the cases of interest discussed in this paper -- see Proposition \ref{prop:even_Iwahori_dimensions}. 

The connected components of $\Gr$, or equivalently $\Gr(\C)$, are in bijection with $\pi_1(G) \cong \Lambda^\vee/Q^\vee$. In particular, there are finitely many connected components when $G$ is semisimple, and $\Gr$ is connected when $G$ is simply connected. An affine Schubert variety $\Gr^{\le \lambda}$ belongs to the component labelled by the coset $\lambda + Q^\vee$. For the real points, we have $\pi_0(\Gr(\R)) \cong \Lambda^\vee_\R/(\Lambda^\vee_\R\cap Q^\vee)$ by \cite[Lem.\,3.8.1]{Nadler2005_PSoRLG}; this means explicitly that a connected component of $\Gr(\C)$ either has no real points or its real points form a single connected component of $\Gr(\R)$. 

\subsubsection{Reduction procedure to the simple case}
\label{subsubsec:Gr_to_sesi}
The assignment $G\to \Gr_G$ is functorial in the following sense: given another real reductive group $H$ and a morphism $G\to H$, one obtains a morphism $\Gr_G\to \Gr_H$. If $G\to H$ is the quotient map with respect to a finite central subgroup of $G$, then this map embeds $\Gr_G$ as a union of connected components of $\Gr_H$ \cite[Prop.\,4.5.6]{BD_Quantization}. From this, one concludes the following \cite[\nopp 4.5.7]{BD_Quantization}: writing $G_\mathrm{ab}\coloneqq G/[G,G]$ for the abelianisation, $\Gr_{G}$ is isomorphic to a union of connected components of $\Gr_{G_\mathrm{ad}}\times \Gr_{G_\mathrm{ab}}$. Moreover, $\Gr_{G_\mathrm{ab}}$ is topologically a discrete space. Lastly, the affine Grassmannian of a product $G_1\times G_2$ is canonically isomorphic to $\Gr_{G_1}\times \Gr_{G_2}$, coweights of $G_1\times G_2$ correspond to pairs of coweights, and we have $\Gr_{G_1\times G_2}^{\le (\lambda_1,\lambda_2)} \cong \Gr_{G_1}^{\le \lambda_1} \times \Gr_{G_2}^{\le \lambda_2}$. Combining these facts, one can reduce topological questions about affine Schubert varieties to the case where $G$ is simple and of adjoint type.

\subsubsection{Comparison to based loop spaces}
\label{subsubsec:Gr_based_loops}
Here we use the associated compact symmetric space $K/H$ from \S\ref{subsubsec:corresponding_sym_sp}. The following is an observation of Quillen, proven by Mitchell:
\begin{thm}[{\cite[Thms.\,4.1, 5.1 and 5.2]{Mitchell1988_ThmOnBuildings}}]
\label{thm:mitchell--quillen}
    Let $G$ be a connected real reductive group, $\Gr = \Gr_G$ its affine Grassmannian, and $K/H$ the corresponding compact symmetric space. Then:
    \begin{itemize}
        \item[(i)] $\Gr(\C)$ is homotopy-equivalent to the based loop space $\Omega K$.
        \item[(ii)] $\Gr(\R)$ is homotopy-equivalent to the based loop space $\Omega (K/H)$.\footnote{In Section 5 of \cite{Mitchell1988_ThmOnBuildings}, $G$ is assumed to be simple and simply connected. One way of obtaining the reductive case from this is to argue componentwise. Indeed, the connected components of $\Gr$ are all homeomorphic, and (see \S\ref{subsubsec:Gr_to_sesi}) can in a first step be identified with connected components of $G_\mathrm{ad}$. The same works on the side of $\Omega(K/H)$, reducing to the semisimple adjoint type case. In this setting, it is easy to take simple factors, thereby reducing to simple $G$. Finally, we may again identify connected components of $\Gr_G$ and $\Omega(K/H)$ with (the unique) connected components for the simply connected covering group, which places things in the setting of \cite{Mitchell1988_ThmOnBuildings}.}   
    \end{itemize}
\end{thm}
We will need to understand the inclusion of $\Gr(\R)$ into $\Gr(\C)$ under these homotopy equivalences, and to this end we now review selected details of the proof of Theorem \ref{thm:mitchell--quillen}. Let us first note that Mitchell uses polynomial loops, in particular the space $\Omega_\mathrm{poly}K \subseteq \Omega K$ of continuous maps $(S^1,1)\to (K,1)$ with finitely many nonzero Fourier coefficients. These are precisely the restrictions to $S^1$ of based morphisms $\C^\times \to G(\C)$, and can therefore be included into $LG(\C)$. This inclusion induces a homeomorphism $\Omega_\mathrm{poly}K\to \Gr(\C)$ (\cf\cite[Thm.\,1.6.1]{Zhu2016_introduction}). Upon identifying $\Gr(\C)$ with $\Omega_\mathrm{poly}K$, the homotopy equivalence in part (i) of the theorem is the inclusion $\Omega_\mathrm{poly}K\to \Omega K$. 
    
Now let $\sigma$ be the real structure corresponding to the real form $G$ of $G_\C$. We also denote by $\sigma$ the resulting involutions on $G(\C)$ (preserving $K$) and $\Gr(\C)$, the latter taking a point $a\colon \Spec \C((t))\to G_\C$ to $\sigma \circ a \circ c$ where $c$ denotes complex conjugation. Note that we have $\Gr(\R) = \Gr(\C)^\sigma$. A similar involution $\sigma$ can be defined on $\Omega K$, by sending a based loop $\gamma$ to $\sigma \circ \gamma \circ c$, and the homeomorphism $\Gr(\C)\to \Omega_\mathrm{poly}K$ identifies $\Gr(\R)$ with $(\Omega_\mathrm{poly}K)^\sigma$. We now observe that a loop $\gamma\in(\Omega K)^\sigma$ is determined by its restriction to the upper unit circle $\{e^{ix}: x\in [0,\pi]\}$, and that this restriction is a path from $1\in K$ to the element $\gamma(-1)\in K^\sigma = H$. In turn, composing such a path with the projection $K\to K/H$ defines a based loop in $K/H$. Altogether, we obtain a map $f\colon (\Omega K)^\sigma \to \Omega (K/H)$ which turns out to be a homotopy equivalence. 

Indeed, an inverse equivalence to $f$ is given by $g\colon \Omega(K/H)\to (\Omega K)^\sigma$, defined as
\begin{equation*}
    g(\beta)(e^{it}) =
    \begin{cases}
        \Tilde{\beta}(e^{2it}) & t\in [0,\pi],\\
        \sigma(\Tilde{\beta}(e^{2i(\pi-t)})) & t\in [\pi, 2\pi].
    \end{cases}
\end{equation*}
where $\Tilde{\beta}$ is a path starting at $1$ and lifting $\beta$ along the fibration $K\to K/H$. To ensure that $g$ is continuous, one really has to take the lifts coherently, but that is possible, see \eg\cite[Sec.\,7.2]{May1999}. The homotopy equivalence in part (ii) of Theorem \ref{thm:mitchell--quillen} is the composition $\Gr(\R) \to (\Omega_\mathrm{poly} K)^\sigma \to (\Omega K)^\sigma \to \Omega(K/H)$. (That the middle step is a homotopy equivalence is the hard part of the theorem.) 

Recall now the Cartan embedding $C\colon K/H\to K$ (see \eqref{eq:Cartan_embedding}). Applying the functor $\Omega$, we obtain a map
\begin{equation*}
    \Omega C\colon \Omega(K/H)\to \Omega K,
    \quad
    \beta \mapsto \big(S^1\ni z\mapsto \beta(z)\sigma(\beta(z))\inv \big). 
\end{equation*}
For paths in a topological group, pointwise multiplication is homotopic to path concatenation, and pointwise inversion is homotopic to path inversion. Thus, $\Omega (K/H)\overset{g}{\longrightarrow} (\Omega K)^\sigma \to \Omega K$ is homotopic to $\Omega C$. The situation is summarised in the following diagram, commutative up to homotopy, in which all vertical arrows are (at least) homotopy equivalences and all horizontal arrows are inclusions. 
\begin{equation*}
\begin{tikzcd}
    \Gr(\R) \ar[r] & \Gr(\C)\\
    (\Omega_\mathrm{poly}K)^\sigma 
        \ar[r] \ar[u] \ar[d]
    & \Omega_\mathrm{poly} K
        \ar[u] \ar[d] \\
    (\Omega K)^\sigma
        \ar[r] \ar[d]
    & \Omega K\\
    \Omega (K/H)
        \ar[ur, "\Omega C"']
\end{tikzcd}
\end{equation*}

We conclude the following:

\begin{prop}
\label{prop:inclusion_to_loop_Cartan}
    The homotopy equivalences of Theorem \ref{thm:mitchell--quillen} identify the inclusion $\Gr(\R)\to \Gr(\C)$ with $\Omega C\colon \Omega(K/H)\to \Omega K$.  
\end{prop}

\subsubsection{The Geometric Satake equivalence}
In this subsection, it is convenient to have $G$ denote a connected \emph{complex} reductive group, and $G_0$ a real form of $G$. A central result related to affine Grassmannians is the geometric Satake equivalence, which relates perverse sheaves on $\Gr_G$ to representations of the Langlands dual group $\Check{G}$. More precisely, let $\mathcal{P}_G$ be the abelian category of $L^+G$-equivariant perverse sheaves on $\Gr_G(\C)$ with coefficients in a commutative Noetherian ring $k$ of finite global dimension, for instance $k = \C$ or $k= \Z$. Then the equivalence can be phrased as follows:
\begin{thm}[{\cite[Prop.\,6.3 and Thm.\,7.3]{MV}}]
\label{thm:geomSat}
    $\mathcal{P}_G$ admits the structure of a neutral Tannakian category with fibre functor given by hypercohomology $\Ha\colon \mathcal{P}_G\to \mathsf{Mod_k}$. Through the Tannakian formalism, $\Ha$ becomes an equivalence of tensor categories $\mathcal{P}_G\to \mathsf{Rep}_{\Check{G}}$, where $\Check{G}$ is the unique split reductive group scheme over $k$ whose root datum is dual to that of $G$. In particular, the intersection cohomology $IH(\Gr_G^{\le \lambda}(\C), k)$ is canonically isomorphic to the $\Check{G}$-module $V^\lambda$ of highest weight $\lambda$.
\end{thm}

\subsubsection{Equivariant cohomology rings}
\label{subsubsec:Gr_cohomology_review}
The cohomology ring of $\Gr(\C)$ can be computed via the homotopy equivalence $\Gr(\C)\simeq \Omega K$ of \S\ref{subsubsec:Gr_based_loops} and the results reviewed in \S\ref{subsec:prelim_blp}. This shows that $H^*(\Gr(\C)_0, \C)$ is a graded polynomial ring, where $\Gr(\C)_0$ denotes the distinguished (or any) connected component, and relates the generator degrees to the topology of $G$. Ginzburg \cite{Ginzburg2000_GeomSat} observed that one can give a more Lie-theoretic description using the geometric Satake equivalence \ref{thm:geomSat}. If $G$ is simply connected, the result is an isomorphism $H^*(\Gr_G(\C),\C) \cong U(\check{\g}_e)$, where $\check{\g}$ is the Langlands dual Lie algebra and $e\in \check{\g}$ a regular nilpotent.

One can extend Ginzburg's description to equivariant cohomology, and here the choice of group action deserves comment. Firstly, $\Gr(\C)$ has a natural action of $L^+G(\C) = G(\C[[t]])$, which contains $G(\C)$ as a subgroup. In fact, the subgroup inclusion is a homotopy equivalence, so that $L^+G(\C)$-equivariant cohomology is the same as $G(\C)$-equivariant cohomology. Secondly, the affine Grassmannian admits an additional ``loop rotation'' $\G_m$-action, descended from the $\G_m$ on $LG$ by rescaling the formal Laurent series parameter. The resulting $G(\C)\times \C^\times$-equivariant cohomology contains more information and has technical advantages over $G(\C)$-equivariant cohomology.\footnote{Mainly, the torus action by $T\times \C^\times$ is a GKM-action \cite{HHH2005_GKM}.} However, loop rotation does not seem to adapt well to real points, so we prefer to avoid it here. 

Another preliminary point concerns connected components (avoided above by taking $G$ simply connected); we identify $\pi_0(\Gr(\C))$ with $\Lambda^\vee/Q^\vee$ and write $\Gr(\C)_\gamma$ for the connected component corresponding to $\gamma\in \Lambda^\vee/Q^\vee$. Although all connected components of $\Gr(\C)$ are homeomorphic, they are not equivariantly homeomorphic. Nevertheless, there are canonical isomorphisms $H^*_{G(\C)}(\Gr(\C)_\gamma,\C)\to H^*_{(G(\C))}(\Gr(\C)_0,\C)$ for all $\gamma\in \Lambda^\vee/Q^\vee$, constructed in \cite{BezrukavnikovFinkelberg2008_EquivariantSatake} (specifically the 2026 arXiv version 4, pp.\,8--9).\footnote{The construction there constructs a $W$-equivariant isomorphism of $T\times \G_m$-equivariant cohomology, using an explicit description of $H^*_{T\times \G_m}(\Gr(\C)_0,\C)$ and concludes by taking $W$-invariants. Since we do not use the $\G_m$-factor, we may either specialise it away or straightforwardly adapt the construction in \cite{BezrukavnikovFinkelberg2008_EquivariantSatake}.} We will  therefore freely identify the algebras $H^*_G(\Gr(\C)_\gamma,\C)$.

Now, the equivariant cohomology ring of $\Gr$ can be computed using localisation techniques \cite{HHH2005_GKM}. A particularly simple closed form, extending Ginzburg's, was given by \cite[Sec.\,2.6]{BezrukavnikovFinkelberg2008_EquivariantSatake} and \cite[Thm.\,6.4]{YZ2011_integral}; the latter with integral coefficients. To describe the result, let $\check{\s}\subseteq \check{\g}$ be a Kostant slice (\cf\cite{Kostant1963_LGRoPR}) and $\jay^\vee_\mathrm{der}$ the locally constant sheaf of abelian Lie algebras over $\check{\s}$ whose fibre over $x\in \check{\s}$ is the centraliser of $x$ in $\check{\g}_\mathrm{der}$.\footnote{In other words, this is a variant of the universal centraliser where the contribution of the centre is deleted from the fibres, but not from the base.} Then the aforementioned closed form can be phrased as follows:
\begin{thm}
\label{thm:Gr_eq_coh}
    $H^*_{G(\C)}(\Gr(\C), \C) \cong \Sym(\jay^\vee_\mathrm{der})\otimes_\C \C^{\Lambda^\vee/Q^\vee}$ as graded algebras over $H^*_{G(\C)}(\pt, \C) \cong \C[\check{\s}]$, where polynomials are graded by twice their usual degree. 
\end{thm}
\begin{rmk}
    It is worth reiterating that \cite{YZ2011_integral} gives such a description with integral coefficients, as well as a description of the Hopf algebra structure after inverting finitely many primes. For the latter, it is often more convenient to work with homology, which is the main setting of \cite{YZ2011_integral}. The work of O'Brien \cite{obrien2024_splitting_rank} extends this to real affine Grassmannians for a family of real groups contained (and almost equal) to the quasi-anisotropic groups discussed in this article. We shall say more about the connection in \S\ref{sec:coh_context}.   
\end{rmk}

As mentioned above, Theorem \ref{thm:Gr_eq_coh} relates to the geometric Satake equivalence. Namely, the equivariant cohomology of $\Gr(\C)$ acts naturally on the equivariant (hyper)cohomology of (perverse) sheaves on $\Gr(\C)$. In particular, it acts on the equivariant intersection cohomology of affine Schubert varieties. These are free modules over $H^*_G(\pt,\C)$, and via Theorem \ref{thm:geomSat} one can identify $IH^*(\Gr^\le(\C), \C)$ with $\C[\check{\s}]\otimes V^\lambda$.
\begin{thm}
\label{thm:eq_coh_aSv}
Let $\lambda\in \Lambda^\vee_+$ be a dominant coweight of $G$. 
    \begin{itemize}
        \item[(i)] The action of $H^*_{G(\C)}(\Gr(\C), \C)$ on $IH^*(\Gr^{\le \lambda}(\C), \C)$ corresponds isomorphically to the action of $\Sym(\jay^\vee_\mathrm{der})$ on $\C[\check{\s}]\otimes_\C V^\lambda$ obtained by the inclusion $\jay^\vee_\mathrm{der}\hookrightarrow \check{\s}\otimes \check{\g}$. (Here we mean the copy of $\Sym(\jay^\vee_\mathrm{der})$ in $\Sym(\jay^\vee_\mathrm{der})\otimes \C^{\Lambda^\vee/Q^\vee}$ corresponding to $\lambda+\Q^\vee$. The other copies act trivially.)
        \item[(ii)] Let $v_\lambda\in V^\lambda$ be a lowest weight vector, and $c_\lambda \in IH^*(\Gr^{\le\lambda}_G(\C),\C)$ the corresponding class. The map $H^*_{G(\C)}(\Gr(\C), \C) \to IH^*(\Gr^{\le \lambda}(\C), \C)$ given by $a\mapsto a\cdot c_\lambda$ factors through an injection $H^*_{G(\C)}(\Gr^{\le \lambda}(\C), \C) \to IH^*(\Gr^{\le \lambda}(\C), \C)$. This identifies $H^*_{G(\C)}(\Gr^{\le\lambda}(\C),\C)$ with $\Sym(\jay^\vee_\mathrm{der})/\operatorname{Ann}(v_\lambda)$, where $\operatorname{Ann}(v_\lambda)$ denotes the annihilator of $1\otimes v_\lambda\in \C[\check{\s}]\otimes V^\lambda$ in $\Sym(\jay^\vee_\mathrm{der})$
        \item[(iii)] We have $\operatorname{Ann}(v_\lambda) = \operatorname{Ann}(V^\lambda)$, where the right-hand side denotes the ideal of $\Sym(\jay^\vee)$ acting on $\C[\check{\s}]\otimes V^\lambda$ by zero. 
    \end{itemize}
\end{thm}
\begin{proof}
    The first statement is \cite[Lem.\,4.9]{BezrukavnikovFinkelberg2008_EquivariantSatake}. A nonequivariant version of (ii) is the main theorem of \cite{Ginzburg2008_variations}, but section 3 of that reference treats also the equivariant version.

    For part (iii) one can argue as follows (though there may well be a simpler reason eluding the author): the action of $\Sym(\jay^\vee_\mathrm{der})$ factors through that of a certain commutative subalgebra $\B\subseteq \End_{\C[\check{\s}]}(\C[\check{\s}]\otimes V^\lambda)$. More precisely, $\B$ is Hausel's big algebra \cite{Hausel2024_avatars}, or equivalently the image of an inhomogeneous Gaudin algebra. By results of \cite{FFR2010_irregular_singularity}, it is known that $v_\lambda$ is a cyclic vector for the action of the fibre $\B_x$ on $V^\lambda$ for $x$ in a nonempty Zariski-open (hence dense) subset of $\check{\s}$, see also \cite[Thm.\,2.1(iv)]{Hausel2024_avatars}.\footnote{To be precise, \cite{FFR2010_irregular_singularity} establishes that $v_\lambda$ is a cyclic vector for the fibre over the regular nilpotent in $\check{\s}$. It is also observed in \cite{FFR2010_irregular_singularity} that admitting a cyclic vector is an open condition for families of subalgebras of $\End(V^\lambda)$, and in fact the same is true for admitting the fixed vector $v_\lambda$ as cyclic vector.} It follows that an element of $\B$ that annihilates $1\otimes v_\lambda$ must annihilate all of $\C[\check{\s}]\otimes V^\lambda$. A fortiori, we have $\operatorname{Ann}(v_\lambda)\subseteq \operatorname{Ann}(V_\lambda)$, and the converse inclusion is clear. 
\end{proof}

As observed above, a key ingredient for the equivariant cohomology of affine Grassmannians is localisation. More precisely, this is the determination of $T$-equivariant (or commonly $T\times \C^\times$-equivariant) cohomology by restriction to the locus of $T$-fixed points. It is not hard to see that $(\Gr(\C))^T = \Lambda^\vee$ (included via $\lambda\mapsto t^\lambda$), and we shall mainly need the following aspect of the localisation theorem:
\begin{thm}[{\cite[Thm.\,2.3 and Sec.\,5]{HHH2005_GKM}}]
    Let $R$ be a subring of $\Q$. The canonical restriction map $\operatorname{Loc}_T\colon H^*_{G(\C)}(\Gr(\C),R)\to H^*_{T(\C)}(\Gr(\C)^{T},R) \cong H^*_{T(\C)}(\pt,R)^{\Lambda^\vee}$ is an injection.
\end{thm}
Although this will not be strictly necessary, let us note that the image of $\operatorname{Loc}_T$ consists of $W$-invariant elements for the action of the Weyl group $W$ on both $H^*_{T(\C)} \cong \C[\te]$ and $\Lambda^\vee$.

From the proof of Theorem 1 in \cite{BezrukavnikovFinkelberg2008_EquivariantSatake} (and further discussion in 2.6 there), one may extract an algebraic description of $\operatorname{Loc_T}$. To give it, let us first decompose $\C[\check{\s}]$ as $\C[\check{\s}\cap\check{\g}_\mathrm{der}]\otimes \C[Z(\g)]$, and further identify $\C[\check{\s}\cap\check{\g}_\mathrm{der}]$ with the (graded) polynomial ring $\C[p_1,\ldots,p_r]$ for basic invariant polynomials $p_1,\ldots, p_r\in \C[\check{\g}_\mathrm{der}]^{\Check{G}_\mathrm{der}}$, restricted to $\check{\s}\cap \check{\g}_\mathrm{der}$. Using the Killing form of $\check{\g}_\mathrm{der}$, we obtain an isomorphism $(\check{\g}_\mathrm{der})^*\cong \check{\g}_\mathrm{der}$, so the derivatives of the $p_i$ can be viewed as elements $p_i'\in\C[\check{\g}]\otimes \check{\g}_\mathrm{der}$ (extended from $\check{\g}_\mathrm{der}$) and restricted to sections $q_i$ of $\jay^\vee_\mathrm{der}$. In fact, the $q_i$ form a basis, so that $\Sym(\jay^\vee_\mathrm{der}) \cong \C[p_1,\ldots, p_s, q_1,\ldots, q_r]$, where $p_{r+1},\ldots, p_s$ are arbitrary linear coordinates of $\C[Z(\g)]$. Using some nondegenerate invariant extension of the Killing form of $\check{\g}_\mathrm{der}$ and the canonical isomorphism $\check{\te}\cong \te^*$, we may also identify $\C[\te]$ with $\C[\check{\te}]$. Finally, $\Check{G}$-equivariance implies that the restriction of $p_i'\colon \check{\g}\to \check{\g}$ to $\check{\te}$ has values in $\check{\te}$, so we may compose it with the derivative of a character $\mu\in \Lambda^\vee = \Hom(\Check{T}, \G_m)$ to obtain a function on $\check{\te}$. With this in mind, the result is as follows:

\begin{prop}
\label{prop:Gr_eq_coh_loc_algebraic}
    The map $\Sym(\jay^\vee)\overset{\cong}{\longrightarrow} H^*_{G(\C)}(\Gr_G(\C),\C)\overset{\operatorname{Loc}_T}{\longrightarrow} H^*_{T(\C)}(\pt, \C)^{\Lambda^\vee} \cong \C[\check{\te}]^{\Lambda^\vee}$ is given explicitly by 
    \begin{equation*}
        p_i\mapsto (p_i|_{\check{\te}})_{\mu\in \Lambda^\vee}, \, i = 1,\ldots, s
        \text{ (no dependence on } \mu\text{),}
        \quad
        \text{and}
        \quad
        q_i\mapsto (p_i'|_{\check{\te}}(\mu))_{\mu\in \Lambda^\vee},
        \,i=1,\ldots, r.
    \end{equation*}
\end{prop}

\section{Medium algebras and automorphisms}
\label{sec:MA}
In this section, we review the medium algebras introduced in \cite{Hausel2024_avatars}. Although their name stems from them being smaller than the so-called big algebras, it turns out that medium algebras connect various interesting topics, thereby earning the name \emph{medium}. After discussing the definition in \S\ref{subsec:MA_def}, we describe in \S\ref{subsec:MA_eq_coh} a relation to the equivariant cohomology of affine Schubert varieties, and in \S\ref{subsec:MA_fibres} how general fibres can be identified with sets of weights. We then introduce the automorphisms to be studied in this article, first from an algebraic perspective (\S\ref{subsec:MA_autos}), and then in more geometric terms (\S\ref{subsec:MA_autos_geometry}), and describe their actions on weight sets (\S\ref{subsec:MA_weight_action}). Finally, \S\ref{subsec:examples} contains some explicit examples.

Throughout, $G$ will denote a connected complex reductive group, with further standard setup and notation as in \S\ref{subsubsec:cplx_red_setup}.

\subsection{Definitions}
\label{subsec:MA_def}
Let $\lambda \in \Lambda^+$ be a dominant integral weight of $G$, and $V^\lambda$ the finite-dimensional irreducible representation of highest weight $\lambda$. The \emph{Kirillov algebra}, introduced by Kirillov \cite{Kirillov2000_family, Kirillov2001_introduction} under the name of \emph{classical family algebra}, is $\Kir^\lambda(\g) \coloneqq (\C[\g]\otimes \End(V^\lambda))^G$. Here $\C[\g]$ is equipped with the $G$-action arising from the adjoint representation, and $\End(V^\lambda)$ with that induced from the $G$-action on $V^\lambda$. The Kirillov algebra is naturally a graded (associative unital) algebra over its subring $\C[\g]^G\cong \C[\g]^G\otimes 1$. Let us note that, since $G$ is connected, the $G$-actions could be replaced by $\g$-actions without changing the result. Elements of $\Kir^\lambda(\g)$ can be viewed as $G$-equivariant morphisms (of affine schemes) $\g\to \End(V^\lambda)$.
\begin{rmk}
    Kirillov's definition actually used $S(\g)$ (or $\C[\g^*]$) where we have used $\C[\g]$. Using $S(\g)$ is in some ways more natural: firstly, it connects better to a similar version using $U(\g)$, and secondly it removes some instances of identifications $\g\cong\g^*$ appearing in the theory. Since we do always have such $G$-equivariant isomorphisms $\g\cong \g^*$ for reductive $G$, the resulting algebras are isomorphic, and we use $\C[\g]$ for convenience.
\end{rmk}

In \cite{Kirillov2000_family}, Kirillov introduces also certain elements $M_P\in \Kir^\lambda(\g)$ depending on invariant polynomials $P\in \C[\g]^G$. To define them, fix a $G$-equivariant isomorphism $\g^*\to \g$ and use it to view the derivative $P'$ as a $G$-equivariant morphism from $\g$ to $\g$. Composing $P'$ with the representation morphism $\rho^\lambda\colon \g\to \End(V^\lambda)$ yields $M_P\coloneqq \rho^\lambda \circ P' \in \Kir^\lambda(\g)$. It is then observed in \cite[Thm.\,M]{Kirillov2000_family} that the $M_P$ are central elements. Hausel \cite{Hausel2024_avatars} calls these elements \emph{medium operators} and defines the \emph{medium algebra} to be the unital subalgebra they generate. Moreover, Hausel proposes that this medium algebra is in fact the centre of $\Kir^\lambda(\g)$. We will use the following equivalent definition:

\begin{defn}
    Let $G$, $\g$, $\lambda$, and $V^\lambda$ be as above, and let $\rho^\lambda\colon U(\g)\to \End(V^\lambda)$ be the structure map of the $U(\g)$-module $V^\lambda$. 
    \begin{itemize}
        \item[(i)] The \emph{universal Kirillov algebra} is $\Kir(\g) \coloneqq (\C[\g]\otimes U(\g))^G$. 
        \item[(ii)] The \emph{universal medium algebra} is the unital $\C[\g]^G$-subalgebra $\M(\g)\subseteq \Kir(\g)$ generated by the subspace $(\C[\g]\otimes \g_\mathrm{der})^G$.
        \item[(iii)] The \emph{medium algebra (for $\lambda$)} is the subalgebra $\M^\lambda(\g)\coloneqq (\id\otimes \rho^\lambda)(\M) \subseteq \Kir^\lambda(\g)$.
    \end{itemize}
\end{defn}
As before, we may view elements of $\Kir(\g)$ as $G$-equivariant morphisms from $\g$ to $U(\g)$.
\begin{rmk}\leavevmode
\begin{enumerate}
    \item As reflected in the notation, these algebras depend on $\g$ rather than $G$. Due to independent interest of the related geometry (see \S\ref{subsec:MA_eq_coh}), it nevertheless seems fitting to work with $G$.
    \item One could replace $\g_\mathrm{der}$ by $\g$ in item (iii) without changing the medium algebras of item (iii).
\end{enumerate}
    
\end{rmk}

\begin{prop}[{\cf\cite{Kirillov2000_family}}]
    We have $\M(\g)\subseteq Z(\Kir(\g))$ and $\M^\lambda(\g)\subseteq Z(\Kir^\lambda(\g))$. In particular, medium algebras are commutative.
\end{prop}
\begin{proof}
    It suffices to prove the first claim. To that end, consider $G$-equivariant morphisms $u\colon \g\to U(\g)$ and $v\colon \g\to \g_\mathrm{der}$. We have to show that $u(x), v(x)\in U(\g)$ commute for every $x\in \g$; in other words, that $v(x)$ annihilates $u(x)$ under the adjoint representation. By equivariance, $u(x)$ and $v(x)$ are invariant under $G_x$, hence annihilated by $\g_x$. In particular, $v(x)$ is annihilated by $x\in \g_x$, which we reinterpret as $v(x)\in \g_x$, so that $v(x)$ indeed annihilates $u(x)$.  
\end{proof}

\subsection{Restriction to a Kostant slice and equivariant cohomology}
\label{subsec:MA_eq_coh}
Let $(e,f,h)$ be a principal $\esl_2$-triple in $\g$. Following \cite[Sec.\,1]{Hausel2024_avatars}, we can describe the above algebras by their restrictions to the Kostant slice $\s \coloneqq e + \g_f$. Recall (see \cite[p.\,382, Thm.\,8]{Kostant1963_LGRoPR}) that $\s$ is contained in the regular locus $\g^\mathrm{reg}$ and intersects each $G$-orbit in $\g^\mathrm{reg}$ precisely once. Now, since $\g^\mathrm{reg}$ is dense in $\g$, any $u\in \Kir(\g)$ is determined by its restriction to $\g^\mathrm{reg}$, and due to $G$-equivariance it is even determined by $u|_\s$. For any $x\in \g$, equivariance also implies that $u(x)\in U(\g)$ is fixed by $G_x$, so that we obtain an embedding of $\Kir(\g)$ into the algebra
\begin{equation*}
    \D\coloneqq \{\varphi\in \operatorname{Mor}(\s, U(\g))\mid \varphi(x)\in U(\g)^{G_x} \text{ for all } x\in \s\}.
\end{equation*}

Conversely, any $\varphi\in \D$ may be extended to an element of $\Kir(\g)$ as follows. Define $\Tilde{\varphi}\colon G\times\s\to U(\g)$ by $\Tilde{\varphi}(g,x)\coloneqq g\cdot \varphi(x)$. We wish to have $\Tilde{\varphi}$ descend to a morphism $\psi\colon \g^\mathrm{reg}\to U(\g)$ such that $\psi(g\cdot x) = \Tilde{\varphi}(g,x)$ for all $x\in \s$ and $g\in G$. To this end, recall the \emph{universal centraliser} $J$, the group scheme over $\s$ whose fibre over $x\in \s$ is $G_x$. (That is, $J$ is the fixed point scheme for the adjoint action, pulled back to $\s$.) The action map $G\times \s\to \g^\mathrm{reg}$ is a $J$-torsor, and $\Tilde{\varphi}$ is $J$-invariant, which implies that $\Tilde{\phi}$ indeed descends to a unique $\psi\colon \g^\mathrm{reg}\to U(\g)$, which is $G$-equivariant by construction. Since the codimension of $\g^\mathrm{reg}$ in $\g$ is known to be 3 (hence more than 2), Hartogs's Lemma implies that $\psi$ extends to a morphism $\g\to U(\g)$, which is again $G$-equivariant since this is the case on a dense subset. We have therefore constructed an inverse to the injection $\Kir(\g)\to A$, and thereby a proof of the following:
\begin{prop}[{\cite[Sec.\,1]{Hausel2024_avatars}}]
\label{prop:kir_slice_restriction}
    Restriction to a Kostant slice $\s$ defines isomorphisms 
    \begin{align*}
        &\Kir(\g) \cong \{\varphi\in \operatorname{Mor}(\s, U(\g))\mid \varphi(x)\in U(\g)^{G_x} \text{ for all } x\in \s\},
        \\
        &\Kir^\lambda(\g) \cong \{A\in \operatorname{Mor}(\s, \End(V^\lambda))\mid A(x)\in \End(V^\lambda)^{G_x} \text{ for all } x\in \s\}.
    \end{align*}
\end{prop}

We can similarly describe medium algebras, using the universal centraliser Lie algebra $\jay_\mathrm{der}\coloneqq \Lie(J) = \{(x,y)\in \s\times \g\mid y\in (\g_\mathrm{der})_x\}$. This is a the total space of a locally free sheaf of abelian Lie algebras over $\s$ (in other words, a vector bundle), and in particular we may form its relative universal enveloping algebra $U_\s(\jay_\mathrm{der}) \cong \Sym(\jay_\mathrm{der})$, a bundle of algebras whose fibre over $x\in \s$ is $U((\g_\mathrm{der})_x)$. We may embed it in the trivial algebra bundle $\s\times U(\g)$, and in particular restrict the morphism $\id \times \rho^\lambda \colon \s\times U(\g) \to \s\times\End(V^\lambda)$ to $U_\s(\jay_\mathrm{der})$. The description of $\Kir(\g)$ of Proposition \ref{prop:kir_slice_restriction} identifies the generating space $(\C[\g]\otimes \g_\mathrm{der})$ of $\M(\g)$ with the sections of $\jay_\mathrm{der}\to \s$. We conclude:
\begin{prop}
\label{prop:MA_slice_restriction}
    Restriction to a Kostant slice $\s\subseteq \g$ identifies $\M(\g)$ with the global section algebra $\Gamma(\s, U_\s(\jay_\mathrm{der}))$, and $\M^\lambda(\g)$ with the image of $\Gamma(\s, U_\s(\jay_\mathrm{der}))$ under $U_\s(\jay_\mathrm{der})\to \s\times \End(V^\lambda)$.
\end{prop}

Comparing Proposition \ref{prop:MA_slice_restriction} with Theorems \ref{thm:Gr_eq_coh} and \ref{thm:eq_coh_aSv} yields the following observation of Hausel:
\begin{prop}[{\cf\cite[Thm.\,3.1]{Hausel2024_avatars}}]
\label{prop:MA_eq_coh}
    Let $\Check{G}$ be the Langlands dual group of $G$, with affine Grassmannian $\Gr_{\Check{G}}$. Identify $\pi_0(\Gr_{\Check{G}}(\C))$ with $\Lambda/Q$, and denote the component corresponding to $\gamma\in \Lambda/Q$ by $\Gr_{\Check{G}}(\C)_\gamma$.
    \begin{itemize}
        \item[(i)] $\M(\g) \cong H^*_{\Check{G}(\C)}(\Gr_{\Check{G}}(\C)_\gamma, \C)$ for every $\gamma\in \Lambda/Q$; in fact $H^*_{\Check{G}(\C)}(\Gr_{\Check{G}}(\C), \C) \cong \M(\g)\otimes_\C \C^{\Lambda/Q}$.
        \item[(ii)] $\M^\lambda(\g) \cong H^*_{\Check{G}(\C)}(\Gr^{\le \lambda}_{\Check{G}}(\C),\C)$.
    \end{itemize}
    All of the above cohomology rings are concentrated in even degrees, and grading them by half the cohomological degree makes all isomorphisms into isomorphisms of graded algebras over the isomorphic base rings $\C[\g]^G\cong H^{2*}_{\Check{G}(\C)}(\pt, \C)$. 
\end{prop}
\begin{rmk}[Loop rotation]
    In fact, Hausel uses a further parameter $\hbar$, which geometrically corresponds to the ``loop rotation'' action of $\G_m$ on $\Gr$ induced by scaling the Laurent series parameter. We may view $\Check{G}(\C)\times \C^\times$-equivariant cohomology as a family over $\A^1\cong \Spec H^{2*}_{\C^\times}(\pt,\C)$ with coordinate $\hbar$. He defines also an $\hbar$-version of Kirillov and medium algebras, and identifies these medium algebras with the $\Check{G}(\C)\times \C^\times$-equivariant cohomology of affine Schubert varieties. The medium algebras considered here are recovered by setting $\hbar$ to zero. We choose to work without loop rotation because it does not seem to interact well with real structures. However, the results of this article concerned with complex automorphisms could be extended to apply to the $\hbar$-versions.
\end{rmk}
\begin{cor}
\label{cor:MA_free}
    $\M^\lambda(\g)$ is a finitely generated free module over $\C[\g]^G$, of rank equal to the number of weights of $V^\lambda$ (counted without multiplicity).
\end{cor}
\begin{proof}
    Indeed, the singular cohomology of complex affine Schubert varieties has a basis corresponding to Iwahori cells (see \S\ref{subsubsec:Gr_points}), so the claim follows from equivariant formality, see Proposition \ref{prop:eq_formality}. 
\end{proof}

\subsection{Fibres and weight sets}
\label{subsec:MA_fibres}
We may think of an element $A\in \M^\lambda(\g) \subseteq \Mor(\g, \End(V^\lambda))^G$ as a family of endomorphisms of $V^\lambda$. In particular, we may evaluate $A$ at a point $x\in \g$; this defines a map $\ev_x\colon \M^\lambda(\g)\to \End(V^\lambda)$. On the other hand, we may project $x$ to a point of $\g/\!/G$ and consider the functions on its fibre. For regular $x$, these operations turn out to coincide:
\begin{prop}
\label{prop:MA_ev=fibre}
    Let $x\in \g^\mathrm{reg}$ and let $\m_x\subseteq \C[\g]^G$ be the kernel of the map $\C[\g]^G\to \C$ which sends $f\mapsto f(x)$. Then 
    \begin{equation*}
        \ev_x(\M^\lambda(\g))
        \cong
        \M^\lambda(\g)\otimes_{\C[\g]^G} \C[\g]^G/\m_x,
    \end{equation*}
    and both coincide with the image of $U(\g_x)$ under $\rho^\lambda\colon U(\g)\to \End(V^\lambda)$.
\end{prop}
\begin{proof}
    It suffices to show that the dimension of $\ev_x(\M^\lambda(\g))$ is at (least the) number of weights of $V^\lambda$ for all $x$ in some Kostant slice $\s = e+\g_f$ (with notation as in \S\ref{subsec:MA_eq_coh}). At $e$, we have $\ev_e(\M^\lambda(\g)) = \rho^\lambda(U(\g_e))$, which has the required dimension by \cite[Thm.\,1.3.2]{Ginzburg2008_variations}. It follows that $\dim\ev_y(\M^\lambda(\g))$ is at least the number of weights for all $y$ in some open neighbourhood $U$ of $e$ in $\s$. For general $x\in \s$, choose $z\in \C^\times$ small enough so that $z\cdot x\in U$. Then, since the dimension of $\ev_x(\M^\lambda(\g))$ equals that of $\ev_{z\cdot x}(\M^\lambda(\g))$ (which is easily seen using a homogeneous generating set), we have that $\dim\ev_x(\M^\lambda(\g)) = \dim\ev_{z\cdot x}(\M^\lambda(\g))$ is at least the number of weights, as we needed to show.
\end{proof}
For the next proposition, recall that our setup \ref{subsubsec:cplx_red_setup} includes a choice of Cartan subalgebra $\te = \Lie(T)$ of $\g$.
\begin{prop}
\label{prop:MA_sesi_fibre}
    If $x\in \te$ is regular, there is a canonical isomorphism $\ev_x(\M^\lambda(\g))\cong \C[\wt(\lambda)]$, where $\wt(\lambda)$ denotes the set of weights of $V^\lambda$. Thus, $(\Spec \M^\lambda(\g))_{\underline{x}} \cong \wt(\lambda)$, where $\underline{x}$ denotes the image of $x$ in $\Spec \C[\g]^G$ and the set $\wt(\lambda)$ is viewed as a finite reduced scheme over $\C$.
\end{prop}
\begin{proof}
    We first observe that $\g_x = \te$, so that Proposition \ref{prop:MA_ev=fibre} yields $\ev_x(\M^\lambda(\g))\cong \rho^\lambda(U(\te))$. Since $U(\te)$ acts by the character given by the weight on each weight space, an application of the Chinese remainder theorem completes the proof. 
\end{proof}

\subsection{Automorphisms: algebraic picture}
\label{subsec:MA_autos}
Medium algebras admit a natural class of automorphisms. We begin by describing them in the algebraic setting of \S\ref{subsec:MA_def}, where they come from automorphisms of Kirillov algebras. The following is essentially a variation of \cite[Sec.\,3]{Elkner2026_minuscule}. 
\begin{defprop}
\label{defprop:autos}
    Let $\eta\in \Aut(G)$ and $\zeta\in \C^\times$. We also write $\eta$ for the automorphisms of $U(\g)$ it induces. 
    \begin{itemize}
        \item[(i)] We define $\Tilde{\iota}_{\eta,\zeta}\in \Aut(\Kir(\g))$ by $\Tilde{\iota}_{\eta,\zeta}(u)(x) \coloneqq \eta(u(\zeta\cdot\eta\inv(x)))$, for all $u\in \Kir(\g)$, viewed as $G$-equivariant morphisms $\g\to U(\g)$, and all $x\in \g$.
        \item[(ii)] $\Tilde{\iota}_{\eta,\zeta}$ restricts to an automorphism of $\M(\g)$ also denoted $\Tilde{\iota}_{\eta,\zeta}$.
        \item[(iii)] 
        Let $\lambda\in \Lambda^+$, with corresponding highest weight representation $V^\lambda$. Let $\eta^*V^\lambda$ denote the $U(\g)$-module with underlying vector space $V^\lambda$ and action twisted by $\eta$, that is, with $u$ acting through $\rho^\lambda\circ \eta$. We call $\eta$ and $\lambda$ \emph{compatible} if $\eta^*V^\lambda \cong V^\lambda$. In this case $\Tilde{\iota}_{\eta,\zeta} \in \Aut(\M(\g))$ (resp.\,$\Aut(\Kir(\g))$) descends to an automorphism $\iota_{\eta, \zeta}$ of $\M^\lambda(\g)$ (resp.\,$\Kir^\lambda(\g)$). 
    \end{itemize}
    These automorphisms all lie over the automorphism $\zeta\circ \eta$ of the base ring $\C[\g]^G$.
\end{defprop}
\begin{proof}
    For (i), one has to check that $(\zeta\circ \eta)\otimes \eta$ preserves $G$-invariance, which follows immediately from the fact that it intertwines the action of $g\in G$ with that of $\eta(g)$. It is then clear that $\Tilde{\iota}_{\eta,\zeta}$ preserves also the space $(\C[\g]\otimes \g_\mathrm{der})^G$ generating $\M(\g)$, hence all of $\M(\g)$. 
    
    Now let $I_\lambda$ denote the kernel of the surjection $\id\otimes \rho^\lambda\colon \Kir(\g)\to \Kir^\lambda(\g)$. It remains to verify that, under the assumptions of part (iii), we have $\Tilde{\iota}_{\eta,\zeta}(I_\lambda)\subseteq I_\lambda$. Explicitly, a $G$-equivariant morphism $u\colon \g\to U(\g)$ is contained in $I_\lambda$ if $u(x)$ annihilates $V^\lambda$ for all $x\in \g$. Checking $\Tilde{\iota}_{\eta,\zeta}(u) \in I_\lambda$ then amounts to showing that also $\eta(u(x))$ annihilates $V^\lambda$ for every $x\in \g$. But the action of $\eta(u(x))$ on $V^\lambda$ is that of $u(x)$ on $\eta^*V^\lambda \cong V^\lambda$, so we indeed have $\Tilde{\iota}_{\eta,\zeta}(u)\in I_\lambda$, as we needed to show. 
\end{proof}
\begin{rmk}
\label{rmk:grading}
    Kirillov and medium algebras are graded, and the parameter $\zeta$ simply implements the resulting $\C^\times$-action. The automorphisms introduced above combine this $\C^\times$-action with an action of $\Aut(G)$ or, as we shall see, $\Out(G)$. We will mainly consider the finite order automorphisms of the form $\iota_{\eta, 1}$ and $\iota_{\eta, -1}$.
\end{rmk}

The following composition rule results immediately from the definitions:
\begin{prop}
\label{prop:MA_auto_composition}
    Let $\eta_1,\eta_2\in \Aut(G)$ and $\zeta_1, \zeta_2\in \C^\times$. Then $\Tilde{\iota}_{\eta_1,\zeta_1}\circ \Tilde{\iota}_{\eta_2,\zeta_2} = \Tilde{\iota}_{\eta_1\eta_2,\zeta_1\zeta_2}$. If both $\eta_1$ and $\eta_2$ are compatible with $\lambda$, then so is $\eta_1\eta_2$ and we have $\iota_{\eta_1,\zeta_1}\circ \iota_{\eta_2,\zeta_2} = \iota_{\eta_1\eta_2,\zeta_1\zeta_2}$.
\end{prop}

\begin{lem}
\label{lem:MA_autos_inner_class}
    Let $\eta, \zeta$ be as in Definition \ref{defprop:autos}, and $g\in G_\mathrm{ad}\cong \Int(G)$.
    \begin{itemize}
        \item[(i)] $\Tilde{\iota}_{g,1} = \id$.
        \item[(ii)] $\Tilde{\iota}_{\eta,\zeta} = \Tilde{\iota}_{\eta g, \zeta}$. 
        \item[(iii)] If $\eta$ is compatible with $\lambda$, then so is $\eta g$, and we have $\iota_{\eta,\zeta} = \iota_{\eta g, \zeta}$.
    \end{itemize}
\end{lem}
\begin{proof}
    Part (i) follows at once from $G$-equivariance in the definition of $\Kir(\g)$, and (ii) follows from (i) via Proposition \ref{prop:MA_auto_composition}. Finally, (iii) amounts to checking that $g^*V^\lambda \cong V^\lambda$ for all $g\in G_\mathrm{ad}$. To see this, we lift $g$ to an element of $G$ (also denoted $g$) and use the identity $\Ad_g(x)\cdot v = g\cdot (x\cdot (g\inv\cdot v)))$, valid for all $x\in \g$ and $v\in V^\lambda$. This exhibits the action of $g\inv$ on $V^\lambda$ as the required isomorphism from $g^*V^\lambda$ to $V^\lambda$.
\end{proof}
\begin{cor}
\label{cor:MA_autos_CSA}
    For fixed $\zeta\in \C^\times$, the automorphism $\Tilde{\iota} \coloneqq \Tilde{\iota}_{\eta,\zeta}$ depends only on the image of $\eta$ in $\Out(G)$. In particular, we may assume without changing $\Tilde{\iota}$, that $\eta$ maps $\te$ to itself. Under this assumption, $\eta$ acts on the weight lattice, and compatibility with $\lambda$ is equivalent to $\eta(\lambda) \in W\cdot \lambda$. 
\end{cor}
\begin{proof}
    The first claim is a rephrasing of Lemma \ref{lem:MA_autos_inner_class}(ii). Next, $\eta\inv$ takes $\te$ to some Cartan subalgebra $\te'$, and since all Cartan subalgebras are conjugate we have $\te' = g(\te)$ for some $g\in \Int(G)$. Replacing $\eta$ with $\eta g$ arranges for $\te$ to be preserved. Assuming now that $\eta(\te) = \te$, we may compare the finite-dimensional irreducible representations $V^\lambda$ and $\eta^* V^\lambda$ using their weights. In fact, two such representations are isomorphic if and only if their sets of extremal weights (which form a single $W$-orbit) agree. But $\eta(\lambda) \in W\cdot \lambda$ is precisely the condition that $V^\lambda$ and $\eta^*V^\lambda$ have the same sets of extremal weights.
\end{proof}
In view of the corollary, the correct parameters for our automorphisms are really $\Out(G)\times \C^\times$:
\begin{defn}
\label{def:autos_via_out}
    Let $a\in \Out(G)$ and $\zeta\in \C^\times$.
    \begin{itemize}
        \item[(i)] We define $\Tilde{\iota}_{a,\zeta} \coloneqq \Tilde{\iota}_{\eta,\zeta}$ for any $\eta\in \Aut(G)$ in the inner class $a$.
        \item [(ii)] Let $\lambda\in \Lambda^+$. We say that $\lambda$ is \emph{compatible} with $a$ if it is compatible with any (hence all) $\eta \in \Aut(G)$ in the inner class $a$. In this case, we define $\iota_{a,\zeta}\coloneqq \iota_{\eta,\zeta}$. 
    \end{itemize}
\end{defn}

\subsection{Automorphisms: geometric picture}
\label{subsec:MA_autos_geometry}
We now connect these automorphisms to the equivariant cohomology setting of \S\ref{subsec:MA_eq_coh}. Due to admitting arbitrary connected reductive $G$, we have to address the fact that $\Gr_{\Check{G}}$ need not be connected. This will ultimately not make a difference for the $\M^\lambda(\g)$ since affine Schubert varieties are always connected, but it means -- as shown in Proposition \ref{prop:MA_eq_coh} -- that the universal medium algebra $\M(\g)$ sees only part of $H^*_{\Check{G}(\C)}(\Gr_{\Check{G}}(\C), \C)$. To better match the geometry, we therefore define a variant of the automorphisms $\Tilde{\iota}_{\eta,\zeta}$ of Definition \ref{defprop:autos}:
\begin{defn}
\label{def:MA_autos_ext}
    Let $\eta\in \Aut(G)$ and $\zeta\in \C^\times$. The action of $\eta$ on $\Lambda/\Q$ induces an automorphism of $\C^{\Lambda/Q}$ also denoted $\eta$. We define the automorphism $\Tilde{\iota}^{\,\mathrm{ext}}_{\eta,\zeta}$ of $\M(\g)\otimes_\C \C^{\Lambda/\Phi}$ by $\Tilde{\iota}^{\,\mathrm{ext}}_{\eta,\zeta} \coloneqq \Tilde{\iota}_{\eta,\zeta}\otimes \eta$.
\end{defn}
\begin{rmk}
    We shall freely identify $\M(\g)\otimes_\C \C^{\Lambda/Q}$ with $\M(\g)^{\Lambda/Q}$, and denote elements by $A=(A_{[\mu]})_{[\mu]\in \Lambda/Q}$.
\end{rmk}

As a first step towards the geometric picture, we now treat restriction to $\te$. In fact, the algebraic setting provides a clear candidate for a restriction map, namely the injective map
\begin{equation}
\label{eq:MA_res_t}
    \res_\te\colon \M(\g)^{\Lambda/Q}
    \to
    (\C[\te]^{\Lambda})^W,
    \quad
    A\mapsto \big(\lambda' \circ (A_{[\mu]})|_\te\big)_{\lambda\in \Lambda, \mu\in \lambda+Q},
\end{equation}
where the $W$-action on $\C[\te]^\Lambda$ is given by $w\cdot (f_\lambda)_\lambda = (w\cdot f_{w\cdot \lambda})_\lambda$, and $\lambda'\colon \te\to \C$ is the derivative of $\lambda$.
Here we have used two facts about $\M(\g)$. Firstly, for any $F\colon \g\to U(\g)$ in $\M(\g)$ we have $F(\te) \subseteq U(\te)$ as can be checked on generators, where it follows from equivariance. It therefore makes sense to compose $F_\te$ with $\lambda'\colon \te\to \C$. Secondly, $G$-equivariance implies that $F|_\te\colon \te\to U(\te)$ is $W$-equivariant, and it follows that $\res_\te$ indeed has values in $(\C[\te]^\Lambda)^W$. That $\res_\te$ is injective follows at once from Zariski-density of weights in $\te^*$, equivariance, and density of the semisimple locus in $\g$. The following lemma, which follows immediately from the definitions, shows how our automorphisms interact with $\res_\te$.
\begin{lem}
\label{lem:MA_autos_after_res}
    Let $\eta\in \Aut(G)$ such that $\eta(\te) = \te$, $\zeta\in \C^\times$, and $\Tilde{\iota}^{\,\mathrm{ext}}_{\eta,\zeta}$ as in Definition \ref{def:MA_autos_ext}. Then we have
    \begin{equation*}
        \res_\te \circ \Tilde{\iota}^{\,\mathrm{ext}}_{\eta,\zeta} = \overline{\iota}_{\eta,\zeta}\circ \res_\te,
    \end{equation*}
    where $\overline{\iota}\in \Aut((\C[\te]^\Lambda)^W)$ is defined to send $f=(f_\lambda)_{\lambda\in \Lambda}\in (\C[\te]^\Lambda)^W$ to $\overline{\iota}(f)$, given by
    \begin{equation*}
        \overline{\iota}(f)_\lambda(x)
        \coloneqq
        f_{\lambda\circ \eta}(\zeta\cdot \eta\inv(x)).
    \end{equation*}
\end{lem}
\begin{rmk}
    Note that, even though $\Tilde{\iota}^{\,\mathrm{ext}}_{\eta, \zeta}$ depends only on $[\eta]\in \Out(G)$, its extension $\overline{\iota}_{\eta, \zeta}$ can really depend on $\eta$ (away from the image of $\res_\te$). 
\end{rmk}

The restriction map \eqref{eq:MA_res_t} indeed corresponds to equivariant localisation:
\begin{lem}
\label{lem:MA_autos_diagram}
    We have the following commutative diagram:
    \begin{equation*}
    \begin{tikzcd}
        \M(\g)^{\Lambda/Q}
            \arrow[d, "\res_\s", "\cong"']
            \arrow[r, "\res_\te"]
        & (\C[\te]^\Lambda)^W
            \arrow[dd, "\cong"']
        \\
        U_\s(\jay_\mathrm{der})^{\Lambda/Q}
            \arrow[d, "\cong"']
        \\
        H^*_{\Check{G}(\C)}(\Gr_{\Check{G}}(\C),\C)
            \arrow[r, "\operatorname{Loc}_{\Check{T}}"]
        & H^*_{\Check{T}(\C)}(\Gr_{\Check{G}}(\C)^{\Check{T}}, \C)^W.
    \end{tikzcd}
    \end{equation*}
    Here, the upper left vertical map denotes restriction to $\s$, and the lower left vertical map is the isomorphism of Theorem \ref{thm:Gr_eq_coh}. The right vertical map arises from $\Gr_{\Check{G}}(\C)^{\Check{T}} = \Lambda$ and isomorphisms $H^*_{\Check{T}(C)}(\pt, \C) \cong \C[\check{\te}]\cong \C[\te]$.
\end{lem}
\begin{proof}
    The map $U_\s(\jay)\to (\C[\te]^\Lambda)^W$ along the lower half of the diagram is given in Proposition \ref{prop:Gr_eq_coh_loc_algebraic}, whose inspection readily shows that the map agrees with $\res_\te\circ \res_\s\inv$, so the diagram is indeed commutative.
\end{proof}

To now make the intended connection to equivariant cohomology, we are lead to look for self-homeomorphisms of $\Gr_{\Check{G}}(\C)$ which act on the coweights $\Lambda$ of $\Check{G}$ according to $\eta$. Moreover, there should be a compatible topological group automorphism of $\Check{T}(\C)$ realising the automorphism $\zeta\circ \eta$ of $\C[\te]\cong H^*_{\Check{T}(\C)}(\pt, \C)$. If $\zeta = 1$, the natural guess is to use automorphisms of $\Check{G}$ and the resulting automorphisms of $\Gr_{\Check{G}}$. The case $\zeta = -1$ turns out to similarly correspond to real structures. For this, we recall from \S\ref{subsubsec:real_red_setup} that real structures act on root data, so that one can assign to each inner class of real structures an element of the outer automorphism group (cf.\,Remark \ref{rmk:real_str_out}).
\begin{thm}
\label{thm:MA_autos_geometry}
    Let $\eta\in \Aut(G)$ such that $\eta(\te) = \te$, and $\lambda\in \Lambda_+$. 
    \begin{itemize}
        \item[(i)] 
        Let $\sigma\in \Aut(\Check{G})$ be an automorphism preserving $\Check{T}$ whose image in $\Out(\Check{G})\cong \Out(G)$ corresponds to that of $\eta\inv$.\footnote{The inverse is due to a mismatch of conventions: the action on medium algebras is constructed to be covariant in the automorphism slot, but the action on equivariant cohomology is not.} The induced automorphism $\sigma\in \Aut(\Gr(\C))$ preserves the affine Schubert variety $\Gr_{\Check{G}}^{\le \lambda}(\C)$ if and only if $\eta(\lambda) \in W\cdot \lambda$. In this case, the resulting automorphism $\sigma^*$ of $H^*_{\Check{G}(\C)}(\Gr_{\Check{G}}^{\le \lambda}(\C),\C)$ corresponds to $\iota_{\eta,1}$ under the isomorphism $H^*_{\Check{G}(\C)}(\Gr_{\Check{G}}^{\le \lambda}(\C),\C)\cong \M^\lambda(\g)$ of Proposition \ref{prop:MA_eq_coh}. 
        \item[(ii)] 
        Let $\sigma$ be a real structure of $\Check{G}$ preserving $\Check{T}$ whose associated class in $\Out(\Check{G})\cong \Out(G)$ corresponds to the image of $\eta$. The induced real structure $\sigma\colon \Gr(\C)\to \Gr(\C)$ preserves the affine Schubert variety $\Gr_{\Check{G}}^{\le \lambda}(\C)$ if and only if $\eta(\lambda) \in W\cdot \lambda$. In this case, the resulting automorphism $\sigma^*$ of $H^*_{\Check{G}(\C)}(\Gr_{\Check{G}}^{\le \lambda}(\C),\C)$ corresponds to $\iota_{\eta,-1}$ under the isomorphism $H^*_{\Check{G}(\C)}(\Gr_{\Check{G}}^{\le \lambda}(\C),\C)\cong \M^\lambda(\g)$ of Proposition \ref{prop:MA_eq_coh}. 
    \end{itemize}
\end{thm}
\begin{proof}
    We start with part (i). In view of Lemmas \ref{lem:MA_autos_after_res} and \ref{lem:MA_autos_diagram}, it suffices to identify the automorphism $\sigma^*$ of $H^*_{\Check{T}(\C)}(\Gr_{\Check{G}}(\C)^{\Check{T}},\C)$ with $\overline{\iota}_{\eta, 1}$. Indeed, $\sigma^*$ is compatible with $\operatorname{Loc}_T$ so it restricts to the automorphism of $H^*_{\Check{G}(\C)}(\Gr_{\Check{G}}(\C),\C)$ induced by $\sigma$. In turn, this automorphism restricts to the desired automorphism of $H^*_{\Check{G}(\C)}(\Gr^{\le \lambda}_{\Check{G}}(\C),\C)$ whenever $\sigma$ maps $\Gr^{\le \lambda}_{\Check{G}}$ to itself. It is easy to check that this condition is equivalent to $\sigma$ fixing the $W$-orbit of $\lambda$, and the action of $\eta$ on $W$-orbits of weights agrees with that of $\sigma$ by construction. 
    
    To identify $\sigma^* \in \Aut(H^*_{\Check{T}(\C)}(\Gr_{\Check{G}}(\C)^{\Check{T}},\C))$ with $\overline{\iota}_{\eta,1}$, we note that both automorphisms are combinations of an automorphism of the base ring $H^*_{\Check{T}(\C)}(\pt,\C) \cong \C[\te]$ and a permutation of $\Lambda$. At this point, we have to temporarily impose further restrictions on $\sigma$. Namely, by replacing $\sigma$ by another automorphism in its inner class (and still preserving $\Check{T}$), we may arrange that $\sigma$ acts on $\check{\te}\cong \te^*$ dually to the action of $\eta\inv$ on $\te$. This change does not affect the action of $\sigma$ on $\Check{G}(\C)$-equivariant cohomology because inner automorphisms act trivially there. With the imposed conditions, both the base ring automorphism and the permutations of $\Lambda$ match, so $\overline{\iota}_{\eta,1}$ indeed corresponds to $\sigma^*$. This finishes the proof of (i).

    For (ii), we first treat the case where $\eta$ is an inner automorphism, so that its image in $\Out(G)$ is the identity. In fact, by Lemma \ref{lem:MA_autos_inner_class}, we may assume $\eta = \id$. As before, the action of a real structure on $\Check{G}(\C)$-equivariant cohomology depends only on the inner class, and a valid choice of $\sigma$ is provided by a split real structure such that $(\Check{T})^\sigma$ is a split real maximal torus. Repeating the localisation argument above, we are now in the situation where the permutations of $\Lambda$ on both sides are trivial, and only have to compare the automorphisms of base rings. To that end, we may identify $\Check{T}(\C)$ with $(\C^\times)^r$ (for $r\in \N$) and $\sigma$ with complex conjugation. We have $H^*(B\C^\times,\C) = H^*(\C\Pe^\infty,\C)\cong \C[c]$, where $c$ can be thought of as the first Chern class of the universal line bundle $\mathcal{L}$ on $\C\Pe^\infty$. Pullback by complex conjugation takes $\mathcal{L}$ to its conjugate line bundle $\overline{\mathcal{L}}$, and we have $c_1(\overline{\mathcal{L}}) = - c_1(\mathcal{L}) = -c$. Hence, complex conjugation acts on $H^*(B\C^\times,\C)$ by $c\mapsto -c$, and therefore on $H^*_{\Check{T}(\C)}(\pt,\C) \cong \C[\check{\te}]$ by $f(x)\mapsto f(-x)$. This matches the base ring action of $\overline{\iota}_{\id,-1}$, which finishes the proof of (ii) in the special case $\eta\in \Int(G)$. 

    To conclude the proof of (ii), consider an arbitrary $\eta$. It is easy to see that $\overline{\iota}_{\eta,-1} = \overline{\iota}_{\eta,1}\circ \overline{\iota}_{\id, -1}$. We now choose an involution $\theta\in \Aut(\Check{G})$ commuting with $\sigma$ and preserving $\Check{T}$ such that $\sigma'\coloneqq\sigma\theta$ is in the inner class of a split real structure (in fact, we may as well assume $\sigma'$ to be split itself). The actions of $\sigma'$ and $\theta$ on equivariant cohomology clearly compose to the action of $\sigma$, and by what has already been proved we may match $(\sigma')^*$ to $\overline{\iota}_{\id,-1}$ and $\theta^*$ to $\overline{\iota}_{\eta, 1}$. Moreover, $\theta$ preserves an affine Schubert variety if and only if $\sigma$ does, and that is the case if and only if $\eta$ preserves the $W$-orbit of the labelling weight. This completes the proof.
\end{proof}

\subsection{Action on weights}
\label{subsec:MA_weight_action}
Fix parameters $a \in \Out(G)$ and $\zeta\in \C^\times$, as well as a dominant weight $\lambda\in \Lambda_+$ compatible with $a$, so that $\iota\coloneqq \iota_{a, \zeta}\in \Aut(\M^\lambda(\g))$ is well-defined. Then $\iota$ acts both on $\M^\lambda(\g)$ and on its base ring $\C[\g]^G$, and we may take spectra to obtain compatible automorphisms $\iota$ of $\Spec \M^\lambda(\g)$ and $\g/\!/G$. The fibres over closed points of the finite morphism $\Spec \M^\lambda(\g)\to \g/\!/G$ are given by Proposition \ref{prop:MA_ev=fibre}. For a closed fixed point $y\in (\g/\!/G)^\iota$, $\iota$ restricts to an automorphism of the fibre, and in view of Proposition \ref{prop:MA_sesi_fibre} this is particularly interesting for regular semisimple $y$, where $(\Spec \M^\lambda(\g))_y \cong \wt(\lambda)$. The resulting action on the set of weights is as follows:
\begin{prop}[{\cf\cite[Prop.\,7.3]{Elkner2026_minuscule}}]
\label{prop:MA_weight_action}
    Let $\iota=\iota_{a,\zeta}$ and $\lambda$ be as above.
    \begin{itemize}
        \item[(i)] If there exist $\eta\in a$ such that $\eta(\te) = \te$, as well as a regular element $h\in \te$ such that $\eta(h) = \zeta(h)$, then the image $\underline{h}$ of $h$ in $\te/\!/W \cong \g/\!/G$ is fixed by $\iota$. Moreover, the action of $\iota$ on $\Spec(\M^\lambda(\g))_{\underline{h}}$ agrees with the action of $\eta$ on $\wt(\lambda)$. 
        \item[(ii)] If $\zeta = 1$, then there exists a pinning-preserving automorphism $\eta_+\in a$ such that $\eta_+(\te) = \te$ and a regular element $h\in \te$ such that $\eta_+(h) = h$. Thus, the action of $\iota$ on the fibre of $\underline{h}$ is given by the action of $\eta_+$ on $\wt(\lambda)$.
        \item[(iii)] If $\zeta = -1$ and $\eta$ has order (at most) 2, then there exists an involutive automorphism $\eta_-\in a$ such that $\eta_-(\te) = \te$ and a regular element $h\in \te$ such that $\eta_-(h) = -h$. Thus, the action of $\iota$ on the fibre of $\underline{h}$ is given by the action of $\eta_-$ on $\wt(\lambda)$. In particular, if $a = 1$, then the action on weights is that of the longest element $w_0$ of $W$.
    \end{itemize}
\end{prop}
\begin{rmk}
\label{rmk:principal_invos}
    The involutive automorphisms $\eta_-$ with a regular semisimple eigenvector for $-1$ appearing in (iii) above are the Cartan involutions of quasi-split real structures, and sometimes called \emph{principal involutions}, see \cite[Def.\,6.3]{AV1992_projective}. They include the Cartan involutions of split real structures, also called \emph{Chevalley involutions}. A compatibly chosen Chevalley involution acts on $\te$ -- and hence on weights -- by $-1$ and is an inner involution if and only if the Weyl group contains $-1$.
\end{rmk}
\begin{proof}[Proof of Proposition \ref{prop:MA_weight_action}]
    By the Chevalley restriction theorem, we have $\C[\g]^G\cong \C[\te]^W$, and $\iota$ acts on $\C[\te]^W$ by sending $f(x)$ to $f(\zeta\cdot\eta\inv(x))$ for any $\eta\in a$ preserving $\te$. Closed points $\underline{x}$ of $\te/\!/W$ can be identified with $W$-orbits $W\cdot x$ of elements $x\in\te$, and $\iota$ fixes $\underline{x}$ if and only if $\eta(x) \in W\cdot \zeta x$. Now, if $\eta$ and $h$ are as in the assumption of (i), then $\underline{h}$ is clearly fixed by $\iota$. Moreover, Proposition \ref{prop:MA_ev=fibre} identifies the ring of functions on the fibre of $\underline{h}$ with the image of $U(\te)$ in $\End(V^\lambda)$, on which $\iota$ acts through $\eta\in \Aut(U(\te))$. It is then clear from the proof of Proposition \ref{prop:MA_sesi_fibre} that this corresponds to the action of $\eta$ on $\wt(\lambda)$.

    It remains to check the existence statements in (ii) and (iii). Temporarily ignoring the compatibility with $\te$, we first find $\eta_+$ or $\eta_-$ in the class $a$ with a regular semisimple eigenvector $h\in \g$ of eigenvalue $1$ resp.\,$-1$. This is accomplished in the former case by Theorem \ref{thm:pi_pr_auto_chara}(iv) below, and in the latter case by \cite[Cor.\,6.15]{AV1992_projective}. The centraliser of $h$ is some Cartan subalgebra $\te'$, so there exists $g\in G(\C)$ such that $g\cdot \te = \te'$ and it suffices to replace $\eta_\pm$ by $g\eta_\pm g\inv$ and $h$ by $g\cdot h$. 
\end{proof}

\subsection{Examples}
\label{subsec:examples}
We now discuss some illustrative explicit examples of medium algebras, their automorphisms, and their coinvariant algebras. The presentations in the first two examples are taken from \cite[Sec.\,4]{Hausel2024_avatars}, and a wealth of descriptions of more complicated medium algebras in type A can be found in \cite[Sec.\,8]{Ngo2026_BigAlgebra}. Here we focus on the connection to actions on weights via Proposition \ref{prop:MA_weight_action} and the appearance of torsion in coinvariant algebras, relevant to Theorem \ref{thm:coinv_main}. The latter shows (in most cases) that, after quotienting out torsion, the coinvariant algebras are free of rank equal to the number of fixed weights.

\begin{expl}[Medium algebras for $\SL_2$]
\label{expl:MA_A1}
    If $G=\SL_2$, then all dominant weights are of the form $k\varpi$ for $k\in \N$, where $\varpi$ is the fundamental weight. The universal medium algebra is $\M(\esl_2(\C)) = \C[c_2,M_1]$ where $c_2$ is a homogeneous generator of $\C[\esl_2(\C)]^{\SL_2}$ and $M_1$ the identity map $\esl_2(\C)\to \esl_2(\C)$. For $k\in \N$ and $h$ an element of the Cartan subalgebra $\te$, $M_1(h)$ acts on the weight spaces of $V^{k\varpi}$ by precisely its weight. Identifying $\te$ with $\C$ (compatibly with $\varpi$), we may choose $c_2$ so that $c_2(h) = h^2$ for all $h\in \te\cong\C$. Then the polynomials $f(c_2,M_1)$ vanishing on the weight space of weight $k\varpi - 2\ell$, $\ell\in 0,\ldots k$, are the multiples of $M_1^2-2(k-2\ell)^2c_2$ if $k-2\ell \ne 0$, and of $M_1$ if $k-2\ell = 0$. Altogether, one obtains
    \begin{equation*}
        \M^{k\varpi}(\esl_2(\C))
        \cong
        \begin{cases}
            \C[c_2,M_1]/((M_1^2-k^2c_2)(M_1^2-(k-2)^2c_2)\cdots(M_1^2 - 4c_2)M_1
            & \text{if }k\text{ is even,}\\
            \C[c_2,M_1]/((M_1^2-k^2c_2)(M_1^2-(k-2)^2c_2)\cdots(M_1^2 - c_2)
            & \text{if }k\text{ is odd.}
        \end{cases}
    \end{equation*}
    (Here we have not argued why the relation found by looking at weight spaces is the only relation. One can check this by comparing ranks, using the fact that medium algebras are free of rank equal to the number of weights.) 

    Since $\Out(\SL_2)$ is trivial, we only have the automorphism $\iota_{1,-1}=(-1)^{\deg}$ in this example. It acts trivially on the base ring $\C[c_2]$, and its coinvariant medium algebras are the $\C[c_2]$-algebras
    \begin{equation*}
        \M^{k\varpi}(\esl_2(\C))_{\iota_{1,-1}}
        \cong
        \begin{cases}
            \C[c_2] & \text{if }k \text{ is even,}\\
            \C[c_2]/(c_2^{(k+1)/2})& \text{if }k \text{ is odd.}
        \end{cases}
    \end{equation*}
    Comparing this to Proposition \ref{prop:MA_weight_action}, the even case matches the fact that the longest Weyl group element $-1$ fixes precisely the zero weight, which appears in $V^{k\varpi}$ when $k$ is even. For $k$ odd, no weights of $V^{k\varpi}$ are fixed, so a general fibre of $\Spec(\M^{k\varpi}(\esl_2(\C))_{\iota_{1,-1}})$ is empty. This also matches: $\C[c_2]/(c_2^{(k+1)/2})$ consists entirely of torsion elements with respect to $\C[c_2]$, so its torsion-free quotient is $0$. 
\end{expl}

\begin{expl}[The adjoint representation of $\SL_3$]
\label{expl:MA_A2_adj}
    The adjoint representation of $\SL_3$ has highest weight $\varpi_1+\varpi_2$. By \cite[Eq.\,4.9]{Hausel2024_avatars}, its medium algebra is
    \begin{equation*}
        \C[c_2,c_3]
        =
        \C[\esl_3(\C)]^{\SL_3}
        \rightarrow
        \M^{\varpi_1+\varpi_2}(\esl_3(\C))
        \cong
        \C[c_2,c_3,M_1,M_2]/(f_1,f_2,f_3)
    \end{equation*}
    where
    \begin{align*}
        &f_1 = M_1^2M_2 + c_2M_2 + 3c_3M_1,\\
        &f_2 = M_1^4 + 4c_2M_1^2 + 3M_2^2,\\
        &f_3 = 3M_1M_2^2 + 9c_3M_2 - c_2M_1^3 - 4c_2^2M_1,
    \end{align*}
    and $c_2,c_3,M_1,M_2$ are homogeneous of degrees $2,3,1,2$ (in this order). 
    
    The nontrivial class $a$ in $\Out(G) = \{1,a\}$ acts on isomorphism classes of representations by dualisation; since the adjoint representation is self-dual, it is compatible with $a$. More explicitly, $a$ is the class of a Chevalley involution, which acts on the Cartan subalgebra $\te$ by $-1$. This implies that $\iota_{a,1}$ acts on the generators by
    \begin{equation*}
        \iota_{a,-1}(c_2) = c_2,
        \quad
        \iota_{a,-1}(c_3) = -c_3,
        \quad
        \iota_{a,-1}(M_1) = M_1,
        \quad
        \iota_{a,-1}(M_2) = -M_2,
    \end{equation*}
    while $\iota_{1,-1}$ again acts by $(-1)^{\deg}$. The third nontrivial automorphism $\iota_{a,-1}$ is obtained as the product of the other two. Since the generators are all eigenvectors of these automorphisms, it is easy to compute the coinvariant algebras -- one only has to delete the generators with eigenvalue not equal to $1$ from the presentation. We find:
    \begin{enumerate}
        \item
        \label{it:MA_A2_adj_1}
        $\M^{\varpi_1+\varpi_2}(\esl_3(\C))_{\iota_{a,1}}
        \cong
        \C[c_2,M_1]/(M_1^2(M_1^2 + 4c_2), c_2M_1(M_1^2 +4c_2))$ 
        over $\C[c_2]$,
        \item 
        \label{it:MA_A2_adj_2}
        $\M^{\varpi_1+\varpi_2}(\esl_3(\C))_{\iota_{1,-1}}
        \cong
        \C[c_2,M_2]/(c_2M_2, M_2^2)$ 
        over $\C[c_2]$, and
        \item 
        \label{it:MA_A2_adj_3}
        $\M^{\varpi_1+\varpi_2}(\esl_3(\C))_{\iota_{a,1}}
        \cong
        \C[c_2,c_3]$ 
        over $\C[c_2, c_3]$.
    \end{enumerate}
    We may observe that the first two coinvariant rings have torsion with respect to their base rings $\C[c_2]$. Quotienting out the $\C[c_2]$-torsion ideals yields:
    \begin{itemize}
        \item $\C[c_2,M_1]/(M_1^3+4c_2M_1)$ in item (\ref{it:MA_A2_adj_1}) above; we may recognise this as $\M^{2\varpi}(\esl_2(\C))$ from Example \ref{expl:MA_A1}. This matches the fact that a pinning-preserving outer automorphism of $\SL_3$ fixes three weights of the adjoint representation, namely the weight string of the highest root through zero. 
        \item In item (\ref{it:MA_A2_adj_2}), quotienting out the torsion ideal $(M_2)$ yields $\C[c_2]$ as an algebra over itself. That this has rank one matches the fact that the longest Weyl group element fixes no roots, hence only the zero weight of the adjoint representation. 
    \end{itemize}
    Similarly, for $\iota_{a,-1}$, Proposition \ref{prop:MA_weight_action} says that the rank should be the number of weights fixed by a Chevalley involution, which again fixes precisely the zero weight.
\end{expl}

\begin{expl}[The second fundamental representation of $\SL_4$]
\label{expl:A3_2nd_fd}
    For a nontrivial example without torsion, consider the second fundamental weight $\varpi_2$ of $\SL_4$, with $V^{\varpi_2} = \Lambda^2\C^4$. We can find a presentation for its medium algebra following a similar strategy as in Example \ref{expl:MA_A1}. We use base ring generators $c_2,c_3,c_4\in \C[\esl_4(\C)]^{\SL_4}$ which restrict to the diagonal Cartan subalgebra as elementary symmetric polynomials. However, for generators of $(\C[\esl_4(\C)]\otimes_\C \esl_4(\C))^{\SL_4}$, it turns out to be more convenient to use derivatives of $p_k\colon A\mapsto \tr(A^k), A\in \esl_4(\C)$ (which restrict to the Cartan subalgebra as power sums). Explicitly, we take $M_i$ to be the derivative of $p_k/k$, viewed as a map to $\esl_4(\C)$ via the isomorphism $\esl_4(\C)^*\to \gl_4(\C)$ afforded by the trace form $(A,B)\mapsto \tr(AB)$ and projected to $\esl_4(\C)$. With these conventions, we obtain
    \begin{equation*}
        \M^{\varpi_2}(\esl_4(\C))
        \cong
        \C[c_2,c_3,c_4,M_1,M_2]/(M_1M_2-c_3,\, M_1^4 - M_2^2 + 2c_2M_1^2 - c_2^2 - 4c_4).
    \end{equation*}

    This representation is again compatible with the nontrivial element $a\in \Out(\SL_4)$, and we compute
    \begin{enumerate}
        \item $\M^{\varpi_2}(\esl_4(\C))_{\iota_{a,1}}
        \cong
        \C[c_2,c_4,M_1]/(M_1^4+2c_2M_1^2 - c_2^2 - 4c_4)$ over $\C[c_2,c_4]$. This is free of rank $4$, matching the fact that four weights of $\Lambda^2\C^4$ are fixed by a compatible pinning-preserving involution.\footnote{To check this, observe that such a pinning-preserving involution should act on the diagonal maximal torus by sending $\diag(x_1,x_2,x_3,x_4)$ to $\diag(x_4\inv,x_3\inv,x_2\inv,x_1\inv)$. The weights of $\Lambda^2\C^2$ are the restrictions of $x_ix_j$, $1\le i< j\le 4$, to the Cartan subalgebra, which has the determinant condition $x_1x_2x_3x_4 = 1$. From this, one easily sees that the fixed weights are the restrictions of $x_1x_2$, $x_1x_3$, $x_2,x_4$ and $x_3x_4$.}
        \item $\M^{\varpi_2}(\esl_4(\C))_{\iota_{1,-1}}
        \cong
        \C[c_2,c_4,M_2]/(M_2^2 + c_2^2 + 4c_4)$, which is free of rank $2$ over $\C[c_2,c_4]$, matching the fact that the longest Weyl group element fixes two weights of $\Lambda^2\C^4$.
        \item  $\M^{\varpi_2}(\esl_4(\C))_{\iota_{a,-1}}
        \cong
        \C[c_2,c_3,c_4]/(c_3, c_2^2 + 4c_4)$ over $\C[c_2,c_3,c_4]$. This consists entirely of torsion elements with respect to $\C[c_2,c_3,c_4]$, so the torsion-free quotient is the zero module. Again, this matches Proposition \ref{prop:MA_weight_action}(iii) since the Chevalley involution $A\mapsto A^{-t}$ of $\SL_4$ fixes no weights of $\Lambda^2\C^4$. 
    \end{enumerate}
\end{expl}

\section{Pinning-preserving automorphisms and quasi-anisotropic groups}
\label{sec:ppqa}
Section \ref{subsec:MA_autos_geometry} describes a class of automorphisms of medium algebras induced by automorphisms or real structures of a complex reductive group, where only the inner classes turn out to be relevant for the action on medium algebras. To use this freedom of choice, we now review distinguished representatives of inner classes for both automorphisms and real structures. In the former case, we will use pinning-preserving automorphisms, already introduced in \S\ref{subsubsec:autos_setup}. For the latter, the best choice turn out to be the \emph{quasi-anisotropic real structures}, whose definition and properties we review here. In this section, we shall use $\theta$ for complex automorphisms and reserve $\sigma$ for real structures.

\subsection{Characterisation and basic properties}
\label{subsec:ppqa_basics}
Recall that a real reductive group $G_0$ is called \emph{quasi-split} if it has a Borel subgroup, that is, a subgroup $B\le G$ such that $B_\C$ is a Borel subgroup of $(G_0)_\C$. In terms of a real structure $\sigma$ of a complex reductive group $G$, a real form is then quasi-split if and only if the corresponding real structure $\sigma$ preserves a Borel subgroup of $G$. In turn, this is equivalent to the existence of a pinning-preserving involution $\theta$ such that $\theta\sigma = \sigma \theta$ and such that $\theta\sigma$ is a split real structure; thus, ``quasi-split'' may be interpreted as ``split up to a compatible pinning-preserving involution''. It is then natural to make the following definition:
\begin{defn}
    A real structure $\sigma$ of a complex reductive group is called \emph{quasi-anisotropic} if a (equivalently any) Cartan involution of $\sigma$ is pinning-preserving, that is, if there exists a pinning-preserving involution $\theta$ such that $\theta \sigma = \sigma \theta$ and such that $\theta\sigma$ is anisotropic. A real reductive group is called \emph{quasi-anisotropic} if the real structure it induces on its complexification is quasi-anisotropic.
\end{defn}
\begin{rmk}
\label{rmk:quasicompact}
    In the literature (including the author's earlier work but also, for instance, \cite[Rmk.\,8.13]{AT_2018_Galois_and_Cartan}), quasi-anisotropic real structures (or forms etc.) are usually referred to as ``quasi-compact''. Since anisotropy seems more meaningful for algebraic groups, as well as due to the unfortunate conflicting usage of ``quasi-compact'' in algebraic geometry, we have chosen to deviate from that convention.
\end{rmk}
\begin{prop}
\label{prop:qa_inner_class}
    Let $G$ be a connected complex reductive group.
    \begin{itemize}
        \item[(i)] Every inner class of automorphisms in $\Aut(G)$ contains a unique inner-isomorphism class of pinning-preserving automorphisms.
        \item[(ii)] Every inner class of real structures of $G$ contains a unique inner-isomorphism class of quasi-anisotropic real structures. 
    \end{itemize}
\end{prop}
\begin{proof}
    Part (ii) follows from Proposition \ref{prop:aut_SES_splitting} and Corollary \ref{cor:pi_pr_auto_uniqueness}. Part (ii) follows from part (i) and Theorem \ref{thm:Cartan_bijection}.
\end{proof}

\begin{expl}[Anisotropic groups are quasi-anisotropic.]
\label{expl:a_is_qa}
    An anisotropic real structure admits the identity as its Cartan involution. Therefore, anisotropic real reductive groups are also quasi-anisotropic. 
\end{expl}

\begin{expl}[Complex reductive groups are quasi-anisotropic.]
    Let $G$ be complex reductive and let $\sigma_0$ be any real structure on $G$. Then $\sigma(g_1,g_2) \coloneqq (\sigma_0(g_2), \sigma_0(g_1))$ defines a real structure of $G\times G$. Choose a Cartan involution $\theta_0$ for $\sigma_0$, and let $\tau_0 \coloneqq \theta_0\sigma_0$ be the corresponding anisotropic real structure. Then the involution $\theta\in \Aut(G\times G)$ given by $\theta(g_1,g_2) = (\theta(g_2), \theta(g_1))$ commutes with $\sigma$ and $\theta\sigma = \tau_0\times \tau_0$ is anisotropic. Thus, $\theta$ is a Cartan involution for $\sigma$ and by taking pinnings of $G$ related via $\theta_0$ one sees that $\theta$ is pinning-preserving. The real form $(G\times G)^\sigma$ can be identified with $G$ viewed as a real reductive group via Weil restriction. In this sense, all complex reductive groups are quasi-anisotropic.
\end{expl}

\begin{expl}[Quasi-anisotropic groups with simple complexification]
\label{expl:simple_qa}
    If $G$ is a connected almost simple group with Lie algebra $\g$, then the possible quasi-anisotropic forms are restricted by the automorphisms of the Dynkin diagram of $\g$, see Remark \ref{rmk:diagram_autos}. When $G$ is of adjoint type, any diagram automorphism induces a pinning-preserving automorphism of $G$.\footnote{Indeed, one can always lift an automorphism from $\g$ to $G_\mathrm{sc}$, and since automorphisms preserve the centre, one can then also descend it to $G_\mathrm{ad}$.} In this case, the inner-isomorphism classes of quasi-anisotropic real structures are in bijection with involutive automorphisms of the Dynkin diagram. The trivial diagram automorphism corresponds to the anisotropic class (see Example \ref{expl:a_is_qa}). If $G$ is not simply laced, there are no other diagram automorphisms. For the simply laced adjont types the isomorphism classes of quasi-anisotropic but not anisotropic real structures are listed in Table \ref{tab:simple_qa}.\footnote{We list isomorphism classes rather than inner-isomorphism classes to avoid digressions about type $D_4$, where there are four quasi-anisotropic inner-isomorphism classes, three of which are in the same isomorphism class. In all other types, this makes no difference.} It may appear surprising that in type $A_{2n}$ the split real structure is also quasi-anisotropic. This coincidence is forced to happen (via Proposition \ref{prop:qa_inner_class}) by the fact that the corresponding inner class of real structures consists of a single inner-isomorphism class.
\end{expl}

\begin{table}[ht]
    \centering
    \caption{Strictly quasi-anisotropic real forms $G_0$ (up to isomorphism) for complex simple adjoint type groups $G$, and fixed subgroups $G^\theta$ of the corresponding pinning-preserving Cartan involutions $\theta$. We write $\PGL_{n,\Ha}$ for the quaternionic projective linear group; this is a real reductive group. The orbits of nodes under the diagram automorphism is indicated by gray lines.}
    \label{tab:simple_qa}
    \begin{tabular}{|c|c|c|c|c|}
        \hline
        Type & $G /\C$ & $G_0 /\R$ & $G^\theta /\C$ & Diagram involution
        \\ \hline
        $A_{2n}, n\ge 1$
        & $\PGL_{2n+1}$
        & $\PGL_{2n+1}$
        & $\PSO_{2n+1}$
        & \dynkin[fold]A{**.**.**}
        \\ \hline
        $A_{2n-1}, n\in \Z_{>0}$
        & $\PGL_{2n}$
        & $\PGL_{n,\Ha}$
        & $\PSp_{2n}$
        & \dynkin[fold]A{**.***.**}
        \\ \hline
        $D_n, n\ge 4$
        & $\PSO_{2n}$
        & $\PSO_{1, 2n-1}$
        & \makecell{$\PSO_{2n-1}$\\ $(\cong \SO_{2n-1})$}
        & \dynkin[fold]D{**.***}
        \\ \hline
        $E_6$
        & $E_{6,\ad}$
        & $E_{6,-26, \ad}$
        & $F_4$
        & \dynkin[fold]E{******}
        \\\hline
    \end{tabular}
\end{table}

We will require the following basic properties of the fixed points under a pinning-preserving automorphism.\footnote{A similar Lemma is given in \cite[\S4]{Elkner2026_minuscule}. In any case, its content is well-known.}
\begin{lem}
\label{lem:pi_pr_fp_pinning}
    Let $G$ be a connected complex reductive group and $\theta\in \Aut(G)$ a finite order automorphism.
    \begin{itemize}
        \item[(i)] $G^\theta$ is reductive. 
        \item[(ii)] If $\theta$ preserves a pinning $(\be, \te, \{e_\alpha\}_\alpha)$ of $\g$, then $\te^\theta$ (resp.\,$\be^\theta$) is a Cartan (resp.\,Borel) subalgebra of $\g^\theta$. The corresponding simple roots of $\g^\theta$ are the restrictions to $\te^\theta$ of the simple roots in $\Pi$, and this yields a bijection of the orbit set $\Pi/\theta$ with the simple roots of $\g^\theta$. Setting $e_A\coloneqq \sum_{\alpha\in A}e_\alpha$ for $A\in \Pi/\theta$, the triple $(\be^\theta,\te^\theta, \{\overline{e_\alpha}\}_{\overline{\alpha}\in \Pi/\theta})$ is a pinning of $\g^\theta$.  
    \end{itemize}
\end{lem}
\begin{proof}
    For the first part, since the unipotent radical is connected, it is enough to show that $(G^\theta)^\circ$ is reductive. In turn, it suffices to show that $\g^\theta$ is reductive. Splitting off the centre, we may further assume that $\g$ is semisimple, so that its Killing form $\kappa$ is nondegenerate. Since $\theta$ has finite order and preserves the $\kappa$ (as all automorphisms do), $\g$ has an orthogonal eigenspace decomposition with respect to $\kappa$, and it follows that the restriction of $\kappa$ to $\g^\theta$ is nondegenerate. Via the adjoint representation of $\g$, $\g^\theta$ is then a subalgebra of $\gl(\g)$ on which the inner product $(A,B)\mapsto \tr(AB)$ restricts nondegenerately; such subalgebras of general linear Lie algebras are well known to be reductive.

    For part (ii), let $\Phi$ be the root system of $\g$ with respect to $\te$, and $\g = \te \oplus \bigoplus_{\alpha\in \Phi} \g_\alpha$ the root space decomposition; note that $\theta$ acts naturally on $\Phi$. Then $\g^\theta = \te^\theta \oplus \bigoplus_{A\in \Phi/\theta} (\g^\theta)_A$ where $(\g^\theta)_A \coloneqq \big(\bigoplus_{\alpha\in A} \g_\alpha\big)^\theta$. Call an orbit $A\in \Phi/\theta$ \emph{admissible} if $\theta^{\#A}$ restricts to the identity on $\g_\alpha$ for any (equivalently all) $\alpha\in A$, and write $(\Phi/\theta)_\mathrm{adm}$ for the set of admissible root orbits; then
    \begin{equation}
    \label{eq:lem_pi_pr_fp_pinning1}
        \g^\theta = \te^\theta \oplus \bigoplus_{A\in (\Phi/\theta)_\mathrm{adm}} (\g^\theta)_A.
    \end{equation}
    Each $(\g^\theta)_A$ is spanned by $\overline{x_\alpha}\coloneqq \sum_{i=1}^{\#A}\theta^i(x_\alpha)$ for any choice of $\alpha\in A$ and $x_\alpha\in \g_\alpha\setminus\{0\}$, and $\overline{x_\alpha}$ is nonzero if and only if $A$ is admissible. Moreover, we have 
    \begin{equation}
    \label{eq:lem_pi_pr_fp_pinning2}
        [h, \overline{x_\alpha}] = \alpha(h) \overline{x_\alpha}
        \quad
        \text{for all } h\in \te^\theta
    \end{equation}
    
    To show that $\te^\theta$ is a Cartan subalgebra of the reductive Lie algebra $\g^\theta$, it suffices to show that it is maximal abelian and consists of semisimple elements. Since semisimplicity can be checked in $\g$, the latter is clear. From the observations of the previous paragraph, checking that $\te^\theta$ is maximal abelian reduces to showing that $[\te^\theta, (\g^\theta)_A]\ne 0$ for all $A\in (\Phi/\theta)_\mathrm{adm}$, but this follows from \eqref{eq:lem_pi_pr_fp_pinning2}. Thus, $\te^\theta$ is indeed a Cartan subalgebra, and now \eqref{eq:lem_pi_pr_fp_pinning2} shows that the $(\g^\theta)_A$, $A\in (\Phi/\theta)_\mathrm{adm}$, are the root spaces, with corresponding root given by restriction of any $\alpha\in A$ to $\te^\theta$.

    Since $\theta$ preserves $\Pi$, it takes positive roots to positive roots, and we infer that $(\Phi/\theta)_\mathrm{adm}^+ \coloneqq (\Phi^+/\theta)\cap (\Phi/\theta)_\mathrm{adm}$ defines a system of positive roots for $(\g^\theta, \te^\theta)$. The corresponding positive root spaces together with $\te^\theta$ clearly span $\be^\theta$, so that $\be^\theta$ is a Borel subalgebra. Since $\theta$ permutes the $e_\alpha$ for $\alpha\in \Pi$, all simple roots are admissible, from which it follows that their orbits form the simple roots in $(\Phi/\theta)_\mathrm{adm}^+$. (To see that there are no other simple roots, note that the restrictions of simple roots span $(\te^\theta)^*$.) This finishes the proof of (ii).
\end{proof}
\begin{rmk}
\label{rmk:fixed_subgroup_pi0}
    The fixed subgroup $G^\theta$ need not be connected. However, if $G$ is \emph{simply connected}, then $G^\theta$ is connected for any finite-order automorphism $\theta\in \Aut(G)$ \cite[Thm.\,8.1]{Steinberg1967_Chevalley_groups}.
\end{rmk}

We take now the opportunity to collect some equivalent characterisations of pinning-preserving automorphisms and quasi-anisotropic real structures. 

\begin{thm}
\label{thm:pi_pr_auto_chara}
    Let $G$ be a connected complex reductive group, $\g$ its Lie algebra, and $\theta \in \Aut(G)$ a finite order automorphism. The following are equivalent:
    \begin{itemize}
        \item[(i)] $\theta$ is pinning-preserving.
        \item[(ii)] The regular elements of $\g^{\theta}$ are also regular elements of $\g$, that is, $(\g_\mathrm{reg})^{\theta} = (\g^{\theta})_\mathrm{reg}$.
        \item[(iii)] $\g^{\theta}$ contains a regular nilpotent element of $\g$.
        \item[(iv)] $\g^{\theta}$ contains a principal three-dimensional subalgebra of $\g$.\footnote{Recall that a principal three-dimensional subalgebra is a subalgebra isomorphic to $\esl_2(\C)$ whose nonzero nilpotent elements are regular in $\g$ (cf.\,\cite{Kostant1959_TDS}).}
        \item[(v)] $G^\theta$ contains a principal three-dimensional subgroup of $G$.\footnote{By a principal three-dimensional subgroup, we mean a closed subgroup whose Lie algebra is principal three-dimensional.}    
    \end{itemize}
\end{thm}
\begin{proof}
    If $\theta$ is pinning-preserving, then by definition there exist a Cartan subalgebra $\te\le \g$ and a choice of simple roots $\Pi$ and nonzero simple root vectors $e_\alpha \in \g_\alpha$, $\alpha\in \Pi$, such that $\theta$ preserves $\te$ and permutes the $e_\alpha$. The element $e\coloneqq \sum_{\alpha \in \Pi}e_\alpha$ is a regular nilpotent of $\g$ \cite[p.\,996]{Kostant1959_TDS}. Moreover, $e\in \g^\theta$ by construction, and Lemma \ref{lem:pi_pr_fp_pinning} shows that $e$ is also a regular nilpotent of $\g^\theta$. We may then (using the Jacobson--Morozov theorem) complete $e$ to an $\esl_2$-triple $(e,h,f)$ in $\g^\theta$; this is a principal $\esl_2$-triple for both $\g$ and $\g^\theta$. Now, $\s \coloneqq e+ C_\g(f)$ and $\s^\theta = e + C_{\g^\theta}(f)$ are Kostant slices for $\g$ resp.\,$\g^\theta$. By \cite[p.\,382, Thm.\,8]{Kostant1963_LGRoPR}, every regular element of $\g^\theta$ is conjugate to an element of $\s^\theta$, and every element of $\s$ is regular in $\g$. Since regularity is preserved by conjugation, this shows that (i) implies (ii).
    
    Assuming (ii), let $x\in \g^\theta$ be any regular nilpotent element. Then $x$ is also a nilpotent element of $\g$ (for instance by compatibility of Jordan--Chevalley decompositions), and by assumption also regular in $\g$. Thus, (ii) implies (iii).

    If $\g^\theta$ contains a regular nilpotent $e$ of $\g$, then as above any extension to an $\esl_2$-triple in $\g^\theta$ yields a principal $\esl_2$-triple of $\g$ contained in $\g^\theta$. The Lie algebra it generates is a principal three-dimensional subalgebra of $\g$, so (iii) implies (iv). Moreover, (iv) and (v) are equivalent since any morphism from $\esl_2$ to $\g$ (or $\g^\theta$) lifts to a morphism from $\SL_2$ to $G$ (or $G^\theta$). 

    To finish the proof, we show that (iv) implies (i). Thus, let $(e,h,f)$ be a principal $\esl_2$-triple of $\g$ contained in $\g^\theta$. Then $h$ is regular semisimple in $\g$, so that $\te\coloneqq C_\g(h)$ is a Cartan subalgebra. Moreover, the regular nilpotent $e$ is contained in a unique Borel subalgebra $\be$ of $\g$ (\cf\cite[Cor.\,5.6 and Thm.\,3.6]{Kostant1959_TDS}). Let $\Pi$ be the system of simple roots corresponding to $\te$ and $\be$. Since $\theta$ fixes $h$ and $e$, we have $\theta(\te) = \te$ and $\theta(\be) = \be$, and it follows that the action of $\theta$ on $\te^*$ permutes $\Pi$. By construction, $e$ is contained in $\bigoplus_{\alpha \in \Pi} \g_\alpha$, and by \cite[Thm.\,5.3]{Kostant1959_TDS}, the component $e_\alpha$ of $e$ in $\g_\alpha$ is nonzero for every $\alpha \in \Pi$. Since $\theta(e) = e$, it follows that $(\be, \te, \{e_\alpha\}_{\alpha\in \Pi})$ is a pinning of $\g$ preserved by $\theta$, so (iv) implies (i).
\end{proof}

\begin{lem}
\label{lem:pi_pr_invo_chara}
    Let $G$ be a connected complex reductive group, $\g$ its Lie algebra, and $\theta \in \Aut(G)$ an involution. If $\theta$ is pinning-preserving, then the restriction map $\C[\g]^\g\to \C[\g^\theta]^{\g^\theta}$ is surjective. If $\theta$ is an involution, the converse also holds. 
\end{lem}
\begin{proof}
    If $\theta$ is pinning-preserving, we have seen in the proof of Theorem \ref{thm:pi_pr_auto_chara} that there is a Kostant slice $\s$ of $\g$ such that $\s^\theta$ is a Kostant slice of $\g^\theta$. By \cite[p.\,382, Thm.\,8]{Kostant1963_LGRoPR}, the natural maps $\C[\g]^\g\to \C[\s]$ and $\C[\g^\theta]^{\g^\theta}\to C[\s^\theta]$ are isomorphisms. Since $\s^\theta$ is a closed subscheme of the affine scheme $\s$, the restriction map $\C[\s]\to \C[\s^\theta]$ is surjective, which proves the first assertion. 

    In the case of involutions, it has been shown in \cite[Rmk.\,4.6]{PanyushevYakimova2021_Sginvolutions} (by classification) that surjectivity $\C[\g]^\g\to \C[\g^\theta]^{\g^\theta}$ is equivalent to item (iii) of Theorem \ref{thm:pi_pr_auto_chara}, hence to $\theta$ being pinning-preserving.
\end{proof}
\begin{rmk}
\label{rmk:res_surjectivity_higher_order}
    By \cite[Thm.\,4.3(1.4)]{HanWen2014_restrictions}, there exists an order 4 automorphism $\theta$ of $\g=\esl_{2n+1}$ (with fixed point subalgebra $\esp_{2n}$) such that $\C[\g]^\g\to \C[\g^\theta]^{\g^\theta}$ is surjective but $\theta$ is not pinning-preserving. This shows that the involutivity assumption in Proposition \ref{lem:pi_pr_invo_chara} is necessary, and explains why a classification-free proof may be difficult.
\end{rmk}

\begin{cor}
\label{cor:qa_chara}
    Let $G$ be a connected real reductive group and $\g$ its Lie algebra. Let $\theta$ be a Cartan involution on $G_\C$ corresponding to $G$. The following are equivalent:
    \begin{itemize}
        \item[(i)] $G$ is quasi-anisotropic.
        \item[(ii)] The regular elements of $(\g_\C)^{\theta}$ are also regular elements of $\g_\C$.
        \item[(iii)] $(\g_\C)^{\theta}$ contains a regular nilpotent element of $\g_\C$.
        \item[(iv)] $(\g_\C)^{\theta}$ contains a principal three-dimensional subalgebra of $\g_\C$.
        \item[(v)] $G$ contains an anisotropic principal three-dimensional subgroup.\footnote{Here, a subgroup $H\le G$ is called a principal three-dimensional subgroup if $H_\C\le G_\C$ is a principal three-dimensional subgroup.}
        \item[(vi)] The restriction map $\C[\g_\C]^{\g_\C} \to \C[(\g_\C)^{\theta}]^{(\g_\C)^\theta}$ is surjective.
        
    \end{itemize}
\end{cor}

\subsection{Maximal tori}
\label{subsec:qa_tori}
For a quasi-anisotropic group $G_0$, the definition in terms of the Cartan involution yields a preferred maximal torus, and it is not clear that it should be maximally split (as required by the standard setup in \S\ref{subsubsec:real_red_setup}). This is one of several instances where the split quasi-anisotropic groups in types $A_{2n}$, $n\in \Z_{>0}$, are more complicated -- for other simple types, all maximal tori of quasi-anisotropic groups turn out to be conjugate. The reason for this can be traced back to properties of root systems:
\begin{lem}
\label{lem:qa_real_and_nci_roots}
    Let $G$ be a connected complex complex reductive group with simple Lie algebra $\g$, and $G_0$ a quasi-anisotropic real form of $G$ with corresponding real structure $\sigma$ and Cartan involution $\theta$. Define roots of $G$ with respect to the maximal torus $T$ of a pinning preserved by $\theta$. Then, unless $\g=\esl_{2n+1}(\C)$ for some $n\in \Z_{>0}$ with nontrivial $\theta$, there are neither real nor noncompact imaginary roots with respect to $\sigma$. 
\end{lem}
\begin{proof}
    The proof follows \cite[Prop.\,6.104]{KnappLG}. We note first that we may assume $\g$ to be simply laced, since otherwise $\theta$ must be trivial, in which case the claim is obvious.
    
    Let $\Pi\subseteq \Phi$ be a system of simple roots involved in a pinning preserved by $\theta$. Recall that a root $\alpha\in \Phi$ is real if and only if $\theta(\alpha) = - \alpha$. Since $\theta$ takes positive roots to positive roots, there are no real roots, as claimed. 
    
    Similarly, $\alpha\in \Phi$ is noncompact imaginary if and only if $\theta(\alpha) = \alpha$ and $\theta$ acts on the root space $\g_\alpha$ by $-1$. We now show by induction on $n$ that no root $\alpha$ expressible as the sum of up to $n$ simple roots is noncompact imaginary. Since $\theta$ permutes a set of simple root vectors, there are no noncompact imaginary simple roots, which establishes the case $n=1$. If $\alpha \in \Phi_+$ is an imaginary root given by the sum of $n\ge 2$ simple roots, we can find a simple root $\beta$ and a positive root $\gamma$ such that $\alpha = \beta + \gamma$ and $\g_\alpha = [\g_\beta,\g_\gamma]$.  Assume first that $\beta$ is imaginary. Then $\gamma$ must be imaginary, and the scalar $\theta|_{\g_\alpha}$ is the product $\theta|_{\g_\beta}\cdot \theta|_{\g_\gamma} = 1\cdot 1 =1$, by induction, so $\alpha$ is compact imaginary. 
    
    If, instead, $\beta$ is complex, we claim that $\alpha'\coloneqq \gamma - \theta(\beta)$ is also a root. To see this, it suffices to show that $(\gamma, \theta(\beta)) > 0$ where $(\cdot, \cdot)$ denotes the pairing induced by the Killing form. Here we use the assumption excluding $\esl_{2n+1}(\C)$ with nontrivial $\theta$, which ensures (by classification of diagram automorphisms) that $(\beta,\theta(\beta)) = 0$. Thus,
    \begin{equation*}
        (\gamma,\theta(\beta))
        =
        (\gamma, \theta(\beta))
        +
        (\beta, \theta(\beta))
        =
        (\alpha, \theta(\beta))
        =
        (\alpha, \beta)
        =
        (\gamma, \beta) + (\beta, \beta)
        \ge 
        -1 + 2
        > 
        0
    \end{equation*}
    where we have used that $\g$ is simply laced. Clearly, the root $\alpha'$ is imaginary, and the induction hypothesis ensures that $\alpha'$ is compact imaginary. Another simple calculation shows that $\alpha$ then must be compact imaginary as well. This finishes the induction, so there are no noncompact imaginary roots, as claimed.
\end{proof}

The relevance of Lemma \ref{lem:qa_real_and_nci_roots} to maximal tori comes from the theory of Cayley transforms (see e.\,g.\,\cite[pp.\,83--85]{AV1992_projective} or \cite[Sec.\,VI.7]{KnappLG}). Cayley transforms change a given maximal torus $T$ of a real reductive group $G_0$ and come in two types: one uses a real root for $T_\C$ to produce a new maximal torus of lower-dimensional split part; the other uses a noncompact imaginary root to produce a new maximal torus of higher-dimensional split part. It is known that all $G_0(\R)$-conjugacy classes of maximal tori can be reached from (any) $T$ by successive Cayley transforms. In particular, if $T_\C$ has neither real nor noncompact imaginary roots, it supports no Cayley transforms, so there cannot be any other conjugacy classes. Thus, the lemma implies the following:
\begin{prop}
\label{prop:qa_max_tori}
    Let $G_0$ be a quasi-anisotropic connected real reductive group. If no simple factors of its Lie algebra are isomorphic to $\esl_{2n+1}(\R)$, $n\in \Z_{>0}$, then all maximal tori of $G_0$ are conjugate through $G_0(\R)$. In particular, they are all maximally split.
\end{prop}

\subsection{Further topological properties}
\label{subsec:ppqa_topological}
Here we describe topological properties of quasi-anisotropic real structures and their Cartan involutions to be used later. The first set of properties is related to the corresponding compact symmetric spaces (see \S\ref{subsubsec:corresponding_sym_sp}).
\begin{prop}
\label{prop:qa_nhz}
    Let $G$ be a connected complex reductive group, $\sigma$ a quasi-anisotropic real structure, and $\theta$ a Cartan involution for $\sigma$. Define $K\coloneqq G(\C)^\tau$ and $H\coloneqq K^\theta$, where $\tau = \theta\sigma$. Then:
    \begin{itemize}
        \item[(i)] The identity component $H^\circ$ is nonhomologous to zero in $K$ (see Definition \ref{def:nhz}). 
        \item[(ii)] $K/H$ is rationally homotopy-simple, that is, if $k$ is a field of characteristic zero, then the action of $\pi_1(K/H)$ on $\pi_n(K/H)\otimes_\Z k$ is trivial for all $n\ge 0$. 
    \end{itemize}
\end{prop}
\begin{proof}
    For part (i), it suffices by Theorem \ref{thm:Samelson} to show that the restriction map $\R[\kay]^K \to \R[\h]^{H^\circ}$ is surjective, where $\kay, \h$ are the Lie algebras of $K$ and $H$, respectively. After (faithfully flat) base change to $\C$, this is equivalent to $\C[\g]^\g = \C[\g]^G\to \C[\g^\theta]^{G^{\theta\circ}} = \C[\g^\theta]^{\g^\theta}$ being surjective, which is the case by Corollary \ref{cor:qa_chara}(vi).

    For (ii), let $p\colon \Tilde{K}\to K$ be the universal covering of $K$, $\Tilde{\theta}$ the lift of $\theta$ to $\Aut(\Tilde{K})$ and $\Tilde{H}\coloneqq \Tilde{K}^{\Tilde{\theta}}$. Then $\Tilde{H}$ is connected and $\Tilde{X}\coloneqq \Tilde{K}/\Tilde{H}$ is the universal covering of $X$. The Deck group of $\Tilde{X}\to X$ is then identified with $p\inv(H)/\Tilde{H}$, acting on $\Tilde{X}$ via right multiplication. For $n\ge 2$ we have $\pi_n(X)\cong \pi_n(\Tilde{X})$ and the action of $\pi_1(X)$ can be identified with the Deck group action on $\pi_n(\Tilde{X})$.

    Thus, we now aim to show that the Deck group $D\coloneqq p\inv(H)/\Tilde{H}$ acts trivially on $\pi_n(\Tilde{X})$ for $n\ge 2$. To this end, recall the Hurewicz homomorphism $h_n\colon \pi_n(\Tilde{X})\otimes_\Z k\to H_n(\Tilde{X},k)$, which is easily seen to be $D$-equivariant. Moreover, it follows from part (i) (applied after lifting to universal coverings), Corollary \ref{cor:nhz_formal} and Theorem \ref{thm:RHT_main}(i) that each $h_n$ is injective. It therefore suffices to show that $D$ acts trivially on $H_*(\Tilde{X},k)$. By the universal coefficient theorems, this is equivalent to $D$, or equivalently $p\inv(H)$, acting trivially on $H^*(\Tilde{X}, \R)$. But, again by Theorem \ref{thm:Samelson}, the canonical map $H^*(\Tilde{X},\R)\to H^*(\Tilde{K}, \R)$ is injective. Moreover, it is $p\inv(H)$-equivariant, which finishes the proof upon observing that the connected group $\Tilde{K} \ge p\inv(H)$ acts trivially on its own cohomology,   
\end{proof}

For another, unfortunately less uniform property, recall the decomposition of $\Gr_G(\R)$ into Iwahori cells (see \S\ref{subsubsec:Gr_points} above). 
\begin{prop}
\label{prop:even_Iwahori_dimensions}
    If $G$ is quasi-anisotropic with simple Lie algebra $\g$ and $\g$ is not isomorphic to $\esl_{2n+1}(\R)$ for $n\in \Z_{>0}$, then all restricted root multiplicities $n_\alpha$ are even.  Hence, the Iwahori cells of $\Gr_G$ all have even dimension, their real points define a basis of the (integral) cohomology of $\Gr_G(\R)$, and the cohomology is concentrated in even degrees. Moreover, the restriction map $H^*(\Gr_G(\R),\Z)\to H^*(\Gr_G^{\le \lambda}(\R), \Z)$ is surjective for every $\lambda\in \Lambda^\vee_{\R,+}$.
\end{prop}
\begin{proof}
    In light of the dimension formula \eqref{eq:Iwahori_cell_dimension} and standard arguments for affine pavings, only the statement about root multiplicities requires proof. In the setting of \S\ref{subsubsec:real_red_setup}, the restriction of $\alpha\in \Phi$ to the split part $S\le T$ vanishes precisely when $\alpha$ is imaginary. Moreover, if $\alpha$ is complex, then $\alpha$ and $\sigma(\alpha)$ have the same restriction, and the pair $\{\alpha, \sigma(\alpha)\}$ contributes a two-dimensional subspace $(\g_\alpha \oplus \g_{\sigma(\alpha)})^\sigma$ to the restricted root space for $\alpha|_{S}$. Now, by Lemma \ref{lem:qa_real_and_nci_roots} and Proposition \ref{prop:qa_max_tori} (which apply due to our assumptions on $G$), there are only imaginary and complex roots, and it follows immediately that all restricted root spaces have even dimensions as claimed.
\end{proof}

When $\g = \esl_{2n+1}(\R)$, say for $G = \PGL_{2n+1,\R}$, there can be odd-dimensional Iwahori cells, and the integral cohomology is not concentrated in even degrees. Nevertheless, we shall see in \S\ref{subsec:coinv_global} that the cohomology of $\Gr_{\PGL_{2n+1}}(\R)$ with complex (equivalently rational) coefficients \emph{is} concentrated in even degrees. The situation becomes more complicated for the affine Schubert varieties, though experimentation suggests the following conjecture.
\begin{conj}
\label{conj:typeA2_Gr_cohomology}
    The canonical maps $H^*(\Gr_{\PGL_{2n+1}}(\R),\Q)\to H^*(\Gr_{\PGL_{2n+1}}^{\le \lambda}(\R),\Q)$ are surjective for all $\lambda\in \Lambda^\vee_{\R, +}$. Moreover, the dimension of $H^*(\Gr_{\PGL_{2n+1}}^{\le \lambda}(\R),\Q)$ equals the number of Iwahori cells in $\Gr_{\PGL_{2n+1}}^{\le \lambda}(\R)$ whose labelling coweights are fixed by the longest element of the Weyl group. 
\end{conj}
\begin{expl}
    To illustrate the second half of the conjecture, we consider two simple examples for $n=2$. Using the usual diagonal maximal torus, the first fundamental coweight of $\PGL_5$ is given by $\check{\varpi}_1\colon \G_m\ni a\mapsto \diag(a,1,1,1,1)$ and the corresponding affine Schubert variety is $\R\Pe^4$. The Iwahori cells (which agree here with usual Schubert cells) in this affine Schubert variety correspond to the coweights in the Weyl group orbit of $\check{\varpi}_1$. Among these, only $a\mapsto\diag(1,1,a,1,1)$ is fixed by the longest Weyl group element (which reverses the order of diagonal entries), so this matches the conjecture. For the second fundamental coweight $a\mapsto\diag(a,a,1,1,1)$ there are two fixed Iwahori cells, which matches the fact that the Grassmannian $\Gr_2(\R^5)$ has Poincaré polynomial $1+q^4$. These examples show also that the analogous conjecture with $\Z$ in place of $\Q$ cannot hold, but it may well hold over $\Z[1/2]$.
\end{expl}

\subsection{Fundamental groups and affine Grassmannian components}
\label{subsec:qa_pi1}
In the remainder of this section, we study how automorphisms interact with the component groups of affine Grassmannians. We first treat the case of pinning-preserving automorphisms, where we work in the equivalent setting of fundamental groups.

\begin{lem}
\label{lem:pi1_image}
    Let $G$ be a connected complex reductive group with a finite order automorphism $\theta$. Then $\theta$ acts on $\pi_1(G(\C))$ and the canonical map $\pi_1(G^\theta(\C))\to \pi_1(G(\C))$ has image $\pi_1(G(\C))^\theta$.
\end{lem}
\begin{proof}
    By conjugating $\theta$ by an inner automorphism if necessary, we may assume that $\theta$ preserves $B$ and $T$. We then have $\pi_1(G(\C)) = \Lambda^\vee/Q^\vee$ and $\pi_1(G^\theta(\C)) = (\Lambda^\vee)^\theta/Q'$, where $Q'$ is the coroot lattice of $G^\theta(\C)$. Now let $k\in \N$ be the order of $\theta$. The short exact sequence $0\to Q^\vee \to \Lambda^\vee \to \pi_1(G(\C)) \to 0$ yields the following exact sequence in group cohomology, where the cyclic group $C_k\coloneqq \Z/k\Z$ acts through $\theta$:
    \begin{equation*}
        (\Lambda^\vee)^\theta 
        \rightarrow
        \pi_1(G(\C))^\theta
        \rightarrow
        H^1(C_k, Q^\vee).
    \end{equation*}
    Since $\theta$ preserves the Borel subgroup $B$ and maximal torus $T$, it permutes the simple coroots, so that we may identify $Q^\vee$ with a direct sum of $C_k$ modules $\Z[C_k/C_\ell]$ where $C_\ell\le C_k$ is a cyclic subgroup. 
    
    We further observe that $\Z[C_k/C_\ell]$ is the induced $C_k$-module from the trivial $C_\ell$-module $\Z$. By Shapiro's lemma (\cf\cite[Lem.\,6.3.2 and 6.3.4]{Weibel_HA}), we then have
    \begin{equation*}
        H^1(C_k, \Z[C_k/C_\ell])
        \cong
        H^1(C_\ell, \Z)
        =
        0.
    \end{equation*}
    It follows that $(\Lambda^\vee)^\theta\to \pi_1(G(\C))^\theta$, hence also $\pi_1(G^\theta(\C))\to \pi_1(G(\C))^\theta$, is surjective.
\end{proof}

\begin{lem}
\label{lem:pi1_kernel}
    Let $G$ be a connected complex reductive group with a pinning-preserving automorphism $\theta$. Let $d$ denote the number of $\theta$-orbits of simple summands $\h$ of $\g = \Lie(G)$ for which there exists $k\in \Z_{>0}$ with $\theta^k(\h) = \h $ and $\theta^k|_{\h}$ exchanging two adjacent nodes of the Dynkin diagram. Equivalently, $d$ is the number of $\theta$-orbits of summands isomorphic to $\esl_{2n+1}(\C)$ (for $n\in \Z_{>0}$) on which some power of $\theta$ acts nontrivially. Then the kernel of $\pi_1(G^\theta(\C))\to \pi_1(G(\C))$ has order $2^d$.
\end{lem}
\begin{proof}
    We have $\pi_1(G^\theta(\C)) = (\Lambda^\vee)^\theta/Q'$, where $Q'$ is the coroot lattice of $G^\theta(\C)$, and the kernel in question is $(Q^\vee)^\theta/Q'$. Since this involves only coroots, we may work in the Lie algebra $\g$ and assume it to be semisimple, say $\g = \bigoplus_{i=1}^\ell \h_i$ for simple Lie algebras $\h_1,\ldots, \h_\ell$. By treating each orbit separately, we may also assume that $\theta$ permutes the $\h_i$ transitively. Then $\theta^\ell$ restricts to an automorphism of each $\h_i$ and we have $\g^\theta \cong \h_1^{\theta_\ell}$. Moreover, $\theta^\ell$ exchanges two adjacent nodes of the Dynkin diagram of $\h_1$ if and only if does so for any $\h_i$. With this in mind, we can replace $\g$ by $\h_1$ and $\theta$ by $\theta^\ell$, and assume $\g$ to be simple. The claim then becomes that the index of $Q'$ in $(Q^\vee)^\theta$ is $2$ if $\theta$ exchanges two adjacent nodes, and $1$ otherwise. 
    
    Now, by Lemma \ref{lem:pi_pr_fp_pinning}(ii), the simple coroots of $\g^\theta$ are in bijection with $\Pi^\vee/\theta$, and inspecting the proof we see that the simple coroot corresponding to $A\in \Pi^\vee/\theta$ has to be some multiple of $\sum_{\check{\alpha}\in A}\check{\alpha}$. From $\g$ being simple it follows that $\g^\theta$ is also simple and the restriction of the Killing form of $\g$ is a multiple of the Killing form of $\g^\theta$. Some calculations, which we skip here, then show that the simple coroot corresponding to $A$ is $2/(2-c_A)\cdot \sum_{\check{\alpha}\in A}\check{\alpha}$, where $c_A$ is the number of $\check{\beta}\in A\setminus \{\check{\alpha}\}$ adjacent to $\check{\alpha}$ in the Dynkin diagram, for some fixed $\check{\alpha}\in A$.

    Since the elements $\sum_{\check{\alpha}\in A}\check{\alpha}$, $A\in \Pi^\vee/\theta$ form a basis of $(Q^\vee)^\theta$, the index of $Q'$ in $(Q^\vee)^\theta$ is given by the product of $2/(2-c_A)$ over $A\in \Pi^\vee/\theta$. Inspecting the possible automorphisms of Dynkin diagram shows that $c_A$ is either $0$ or $1$, with $1$ occurring exactly for the two middle nodes in type $A_{2n}$, swapped by its nontrivial diagram automorphism. Thus, the factor $2/(2-1) = 2$ occurs at most once, and only in that case. This finishes the proof.
\end{proof}

Lemmas \ref{lem:pi1_image} and \ref{lem:pi1_kernel} characterise the map $\pi_0(\Gr_G(\C)^\theta)\to \pi_0(\Gr_G(\C))$ for a pinning-preserving automorphism $\theta$ of a connected complex reductive group $G$. For a real structure $\sigma$, the situation is different: as recalled in \S\ref{subsubsec:Gr_points}, $\pi_0(\Gr_G(\C)^\sigma) = \pi_0(\Gr_{G^\sigma}(\R))$ injects into $\pi_0(\Gr_G(\C))^\sigma$ by \cite[Lem.\,3.8.1]{Nadler2005_PSoRLG}. However, unlike the complex case, this map need not be surjective. In other words, there can exist connected components of $\Gr_G(\C)$ preserved by $\sigma$ which do not intersect the real locus.
\begin{defn}
\label{def:adm_comps}
    Let $G$ be a connected complex reductive group with a real structure $\sigma$. A connected component of $\Gr(\C)$ is called \emph{admissible (with respect to $\sigma$)} if it intersects $\Gr(\C)^\sigma$. We denote the union of admissible components by $\Gr(\C)_\mathrm{adm}$.
\end{defn}
The definition of admissibility is made for the following to hold:
\begin{lem}
\label{lem:pi0_real_Gr}
    Let $G$ be a connected complex reductive group with a real structure $\sigma$. Then $\pi_0(\Gr_G(\C)^\sigma)\to \pi_0(\Gr_G(\C))^\sigma$ is injective with image $\pi_0(\Gr_G(\C)_\mathrm{adm})$.  
\end{lem}

Identifying $\pi_0(\Gr(\C))$ with $\Lambda^\vee/Q^\vee$, we can be more explicit about admissibility. We assume the setup of \S\ref{subsubsec:real_red_setup}, so in particular the implicit maximal torus $T$ and Borel subgroup $B$ are chosen compatibly with $\sigma$. Now let $\gamma\in \Lambda^\vee/Q^\vee$ and suppose its component $\Gr_G(\C)_\gamma$ has a real point. By Lemma \ref{lem:pi0_real_Gr}, $(\Gr_G(\C)_\gamma)^\sigma$ is then a single connected component of $\Gr_{G^\sigma}(\R)$, and by the decomposition \eqref{eq:Schubert_cell_partitions} of $\Gr_{G^\sigma}(\R)$ into affine Schubert cells, it must contain a closed affine Schubert cell $\Gr_{G^\sigma}^\lambda$ for some $\lambda\in \Lambda^\vee_{\R,+}$. That this cell is closed is equivalent to $\lambda$ being a minimal dominant coweight under the partial ordering $\le$ (\cf\eqref{eq:aSc_closures}), in other words, $\lambda$ is minuscule. Since every $Q^\vee$-coset contains a unique minuscule coweight, we conclude:
\begin{prop}
\label{prop:adm_minu}
    In the above setting, the connected component of $\Gr_G(\C)$ labelled by $\gamma\in\Lambda^\vee/Q^\vee$ is admissible if and only if the unique minuscule coweight $\lambda\in \gamma$ is in $\Lambda^\vee_\R$.
\end{prop}

Thus, admissible components correspond to minuscule coweights fixed by $\sigma$. We end this section by compiling these coweights for quasi-anisotropic real structures of simple groups of adjoint type in Table \ref{tab:admissible_minuscules} below (see also Example \ref{expl:simple_qa}).\footnote{The relevant computations can be found in \ref{sec:appx_setup}.}
\begin{table}[ht]
\label{tab:admissible_minuscules}
    \centering
    \caption{Minuscule real coweights for quasi-anisotropic real forms $G_0$ of simple groups $G$ of adjoint type. We use $\check{\varpi}_i$ to denote the $i$-th fundamental coweight. The zero coweight is always real; we label the remaining minuscule real coweights in the Dynkin diagram of $\Lie(G)$ by putting a cross at the corresponding nodes.}
    \begin{tabular}{|c|c|c|}
        \hline
        $G_0$ & real minuscule coweights & Dynkin diagram \\
        \hline
        $\PGL_{2n+1}$, $n\in \Z_{>0}$ & all minuscule coweights & \dynkin A{tt.t}\\
        \hline
        $\PGL_{n,\Ha}$, $n\in \Z_{>0}$ & $0$ and $\check{\varpi}_{2i}$, $i=1,\ldots, n-1$ & \dynkin A{otot.oto}\\
        \hline
        $\PSO_{1,2n-1}$, $n\in \Z_{>0}$ & $0$ and $\check{\varpi}_1$ & \dynkin D{too.ooo}\\
        \hline
        $E_{6,-26,\mathrm{ad}}$ & all minuscule coweights & \dynkin E{toooot}\\
        \hline
        anisotropic & $0$ & unmarked\\
        \hline
    \end{tabular}
\end{table}

\section{Coinvariant algebras and fixed loci}
\label{sec:coinv}
\subsection{Overview}
We now come to our main results concerning medium algebras, relating their coinvariant algebras to equivariant cohomology. Recall that the automorphisms $\iota\colon \M^\lambda(\g)\to \M^\lambda(\g)$ introduced in Section \ref{subsec:MA_autos} preserve $\C[\g]^G\subseteq \M^\lambda(\g)$, so that they can be viewed as automorphisms of the algebra $\C[\g]^G\to \M^\lambda(\g)$. As in the introduction, the \emph{coinvariant algebra} of $\iota$ is then the algebra
\begin{equation*}
    \C[\g]^G/(f-\iota(f)\mid f\in \C[\g]^G) = \C[\g]^\G_\iota \longrightarrow \M^\lambda(\g)_\iota = \M^\lambda(\g)/(A-\iota(A)\mid A\in \M^\lambda(\g)),
\end{equation*}
the algebraic counterpart of the map of fixed points $(\Spec \M^\lambda(\g))^\iota\to (\g/\!/G)^\iota$. 

Now, $\M^\lambda(\g)$ may be viewed as the equivariant cohomology ring of the $\Check{G}$-affine Schubert variety $\Gr^{\le \lambda}_{\Check{G}}(\C)$. As explained in \S\ref{subsec:MA_autos_geometry}, our automorphisms are then given by the equivariant cohomology action $\sigma^*$ for an automorphism or real structure $\sigma$ of $\Check{G}$, and the goal is to describe the coinvariant algebra in terms of the $\Check{G}(\C)^\sigma$-equivariant cohomology of $\Gr_{\Check{G}}^{\le\lambda}(\C)^\sigma$. 

Before realising this description for any particular $\lambda$, we study the case of full affine Grassmannians. In \S\ref{subsec:coinv_global} we show that essentially the canonical map $H^*_{\Check{G}(\C)}(\Gr_{\Check{G}}(\C),\C)\to H^*_{\Check{G}(\C)^\sigma}(\Gr_{\Check{G}}(\C)^\sigma,\C)$ is surjective and factors through an isomorphism $H^*_{\Check{G}(\C)}(\Gr_{\Check{G}}(\C),\C)_{\sigma^*} \cong H^*_{\Check{G}(\C)^\sigma}(\Gr_{\Check{G}}(\C)^\sigma,\C)$. Here we have ignored issues related to connected components (see \S\ref{subsubsec:conn_components}), and the actual statement has to be corrected accordingly. Using equivariant formality, this global step can be reduced to nonequivariant cohomology, where we may use the rational homotopy theory methods reviewed in \S\S\ref{subsec:prelim_RHT}--\ref{subsec:prelim_blp}. 

In \S\ref{subsec:asV_fixed_loci}, we give a precise description of the fixed loci in affine Schubert varieties, before turning finally to coinvariants of medium algebras and the equivariant cohomology of these fixed loci in \S\ref{subsec:coinv_aSv}. Unlike the global case, these coinvariant algebras can have torsion, whereas their equivariant cohomology models are free (as modules). It is therefore necessary to take the torsion-free quotient of the coinvariant algebra. Let us also note that the proofs of our main theorem \ref{thm:coinv_main} about medium algebra coinvariants and affine Schubert varieties do not fully use the main global result of \S\ref{subsec:coinv_global}. Nevertheless, these results are clearly related, and we discuss this further in \S\ref{subsec:coinv_global_to_aSv}. This also connects naturally to the aforementioned issue of torsion.

\subsubsection{Connected components}
\label{subsubsec:conn_components}
We now address three complications regarding connected components. Firstly, $\Check{G}(\C)^\sigma$ need not be connected. In the case of a complex (pinning-preserving) automorphism, this really means that the algebraic group $\Check{G}^\sigma$ need not be connected, as noted in Remark \ref{rmk:fixed_subgroup_pi0}. In the case of a (quasi-anisotropic) real structure, the real form $\Check{G}^\sigma$ of $\Check{G}$ has to be connected as an algebraic group, but $\Check{G}(\C)^\sigma = \Check{G}^\sigma(\R)$ need not be connected in the classical topology. In either case, we will have to take equivariant cohomology with respect to the identity component $\Check{G}(\C)^{\sigma \circ}$.

Secondly, as discussed in \S\ref{subsec:qa_pi1}, there can be connected components of $\Gr_{\Check{G}}(\C)$ preserved by a real structure $\sigma$ which nevertheless have no real points. These components contribute to the coinvariant cohomology ring (at least in degree zero), and thereby obstruct injectivity of this ring to the cohomology of the real locus. To avoid this, we restrict to admissible components (see Definition \ref{def:adm_comps}).

Lastly, if $\sigma$ is a complex (pinning-preserving) automorphism, several connected components of $\Gr_{\Check{G}}(\C)^\sigma$ can be contained in a single connected component of $\Gr_{\Check{G}}(\C)$. For the restriction map in (equivariant) cohomology, it means that a direct factor of $H^*(\Gr_{\Check{G}}(\C))$ corresponding to a connected component can map into several factors of $H^*(\Gr_{\Check{G}}(\C)^\sigma)$.  However, all connected components of $\Gr_{\Check{G}}(\C)^\sigma$ are  isomorphic (as can be seen by viewing it as $\Gr_{(\Check{G})^\sigma}(\C)$), and their cohomology rings are canonically isomorphic. The equivariant cohomology rings are also canonically isomorphic (cf.\,\S\ref{subsubsec:Gr_cohomology_review}). 

It will therefore be sufficient to keep track of the set of components of $\Gr_{\Check{G}}(\C)^\sigma$ contained in a given component of $\Gr_{\Check{G}}(\C)$. Now, the component sets in question are isomorphic to the fundamental groups $\pi_1(\Check{G}^\sigma)$ and $\pi_1(\Check{G})$, respectively, so this is about cosets of the kernel of $\pi_1(\Check{G}^\sigma)\to \pi_1(\Check{G})$. This component issue will be addressed using the following symmetric group action:
\begin{defn}
\label{def:components_action}
    If $\sigma$, as above, is a pinning-preserving automorphism, we put $\kappa_\sigma \coloneqq \ker(\pi_1(\Check{G}^\sigma(\C))\to \pi_1(\Check{G}(\C)))$. We let the symmetric group $S(\kappa_\sigma)$ act on $H^*(\Gr_{\Check{G}}(\C)^\sigma,\C)$ and $H^*_{\Check{G}(\C)^{\sigma\circ}}\Gr_{\Check{G}}(\C)^\sigma,\C)$ by permuting component factors within $\kappa_\sigma$ cosets using the canonical isomorphisms between cohomology rings of connected components.
\end{defn}
\begin{rmk}
    By Lemma \ref{lem:pi1_kernel}, $\kappa_\sigma$ is always finite, and is trivial if $\Check{G}$ has no type $A_{2n}$ (almost) factors fixed by some power of $\sigma$ which acts on it nontrivially.
\end{rmk}
\begin{rmk}
    The purpose of the group $S(\kappa_\sigma)$ is that the maps $H^*(\Gr_{\Check{G}}(\C),\C)\to H^*(\Gr_{\Check{G}}(\C)^\sigma,\C)$ and $H^*_{\check{G}(\C)}(\Gr_{\Check{G}}(\C),\C)\to H^*_{\check{G}(\C)^{\sigma\circ}}(\Gr_{\Check{G}}(\C)^\sigma,\C)$ have values in the $S(\kappa_\sigma)$-fixed parts. Let us briefly sketch why this is the case. 
    
    Namely, suppose that $\lambda_1,\lambda_2\in \Lambda^\sigma$ define the same component $Y$ of $\Gr_{\Check{G}}(\C)$ but different components $X_{\lambda_1}, X_{\lambda_2}$ of $\Gr_{\Check{G}^\sigma}(\C)$. Left multiplication with $t^{\mu}$, where $\mu = \lambda_1-\lambda_2$ defines a homeomorphism $\phi_\mu\colon X_{\lambda_2}\to X_{\lambda_1}$. Since the component of $LG^\sigma(\C)$ containing $t^{\mu}$ does not depend on the choices of $\lambda_1$, $\lambda_2$, all of these homeomorphisms are homotopic and define the canonical isomorphism $\phi^*\colon H^*(X_{\lambda_1},\C)\to H^*(X_{\lambda_2},\C)$ used above. In $LG(\C)$, $t^\mu$ is contained in the identity component, so that $\phi^*$ indeed acts trivially on the image of $H^*(Y,\C)$.

    The equivariant version of this (as in the arXiv version 4 of \cite{BezrukavnikovFinkelberg2008_EquivariantSatake}, pp.\,8--9) again starts from $\phi_\mu^*$ as above, but has to introduce a correction depending on $\mu$. This is constructed in such a way that the composition of $\phi_\mu^*$ with that correction does not depend on $\mu$, and for the same reason it is the identity when $\mu$ is in the coroot lattice. By assumption, $\mu$ is indeed in the coroot lattice of $\Check{G}$, so as before this component identification acts trivially on the image of $H^*_{\Check{G}(\C)}(Y,\C)$.   
\end{rmk}

\subsection{Coinvariants for the full affine Grassmannian}
\label{subsec:coinv_global}
In this subsection, we prove the following:

\begin{thm}
\label{thm:coinv_global}
    Let $G$ be a connected complex reductive group and $\sigma$ either a pinning-preserving automorphism or a quasi-anisotropic real structure of $G$. 
    \begin{itemize}
        \item[(i)] The canonical map $H^*_{G(\C)}(\pt, \C)\to H^*_{G(\C)^{\sigma\circ}}(\pt, \C)$ factors through an isomorphism 
        \begin{equation*}
            H^*_{G(\C)}(\pt, \C)_{\sigma^*}
            \overset{\cong}{\rightarrow} 
            H^*_{G(\C)^{\sigma\circ}}(\pt, \C).
        \end{equation*}
        \item[(ii)] If $\sigma$ is a pinning-preserving automorphism, then the canonical map $H^*_{G(\C)}(\Gr_G(\C), \C)\to H^*_{G(\C)^{\sigma\circ}}(\Gr_G(\C)^\sigma, \C)$ factors through an isomorphism
        \begin{equation*}
            H^*_{G(\C)}(\Gr_G(\C), \C)_{\sigma^*}
            \overset{\cong}{\rightarrow}
            H^*_{G(\C)^{\sigma\circ}}(\Gr_G(\C)^\sigma, \C)^{S(\kappa_\sigma)}
            =
            H^*_{G^{\sigma\circ}(\C)}(\Gr_{G^\sigma}(\C), \C)^{S(\kappa_\sigma)}.
        \end{equation*}
        \item[(iii)] If $\sigma$ is a quasi-anisotropic real structure, then the canonical map $H^*_{G(\C)}(\Gr_G(\C)_\mathrm{adm}, \C)\to H^*_{G(\C)^{\sigma\circ}}(\Gr_G(\C)^\sigma, \C)$ factors through an isomorphism
        \begin{equation*}
            H^*_{G(\C)}(\Gr_G(\C)_\mathrm{adm}, \C)_{\sigma^*}
            \overset{\cong}{\rightarrow}
            H^*_{G(\C)^{\sigma\circ}}(\Gr_G(\C)^\sigma, \C)
            =
            H^*_{G^{\sigma}(\R)^\circ}(\Gr_{G^\sigma}(\R), \C).
    \end{equation*}
    \end{itemize}
\end{thm}

We begin with the base rings. We may identify $H^*_{G(\C)}(\pt, \C) = H^*(BG(\C),\C)$ with $\C[\g]^G$. If $\sigma$ is a complex automorphism, then we similarly have $H^*_{G(\C)^{\sigma\circ}}(\pt, \C) \cong \C[\g^\sigma]^{G^{\sigma\circ}}$, and the canonical map $H^*_{G(\C)}(\pt, \C)\to H^*_{G(\C)^{\sigma\circ}}(\pt, \C)$ is identified with restriction of invariant polynomials. In the case of a quasi-anisotropic real structure, we note that $G(\C)^{\sigma\circ}$ is homotopy-equivalent to its maximal compact subgroup $(G(\C)^{\theta\circ})^\sigma$ where $\theta$ is the Cartan involution of $\sigma$. But this is also a maximal compact subgroup of (hence homotopy-equivalent to) $G(\C)^{\theta\circ}$, so we have $H^*_{G(\C)^\sigma}(\pt, \C)\cong \C[\g^\theta]^{G^{\theta\circ}}$ in that case. Similarly, one sees that the isomorphism $H^*_{G(\C)}(\pt, \C)\cong \C[\g]^G$ takes $\sigma^*$ to the action of $\sigma$ (that is, $f(x)\mapsto f(\sigma(x))$) in the former case, and to the action of $\theta$ in the latter. Since Cartan involutions of quasi-anisotropic real structures are by definition pinning-preserving, we have reduced to the action of a pinning-preserving automorphism on the ring of invariant polynomials. The following is proved in \cite[Lem.\,4.7 and Cor.\,4.8]{Elkner2026_minuscule}:
\begin{lem}
\label{lem:inv_poly_coinvariants}
    Let $\theta$ be a pinning-preserving automorphism of $G$, and let $\theta^*\colon \C[\g]^G\to \C[\g]^G$ denote precomposition by $\theta$. Then the restriction map $\C[\g]^G\to \C[\g^\theta]^{G^\theta}$ factors through an isomorphism $\C[\g]^G_{\theta^*}\to \C[\g^\theta]^{G^\theta} = \C[\g^\theta]^{G^{\theta\circ}}$.  
\end{lem}
\begin{cor}[{\cf\cite[Cor.\,4.10]{Elkner2026_minuscule}}]
\label{cor:base_ring_coinvariants}
    For $\sigma$ as above, the canonical map $H^*_{G(\C)}(\pt, \C)\to H^*_{G(\C)^{\sigma\circ}}(\pt, \C)$ factors through an isomorphism $H^*_{G(\C)}(\pt, \C)_{\sigma^*}\to H^*_{G(\C)^{\sigma\circ}}(\pt, \C)$.
\end{cor}
\begin{rmk}
    Lemma \ref{lem:inv_poly_coinvariants} shows in particular that the map of base rings $H^*_{G(\C)}(\pt, \C)\to H^*_{G(\C)^{\sigma\circ}}(\pt,\C)$ is surjective, which is certainly necessary for it to induce an isomorphism. By Proposition \ref{lem:pi_pr_invo_chara}, this surjectivity characterises pinning-preserving involutions. Thus, we are essentially forced to use quasi-anisotropic real structures. (For higher-order automorphisms this is not quite the case, see Remark \ref{rmk:res_surjectivity_higher_order}). 
\end{rmk}

To continue uniformly, let us agree to write $X=\Gr_G(\C)$ if $\sigma$ is a complex automorphism and $X=\Gr_G(\C)_\mathrm{adm}$ if $\sigma$ is a real structure. We now come to the map $i^*\colon H^*_{G(\C)}(X,\C)\to H^*_{G(\C)^{\sigma\circ}}(X^\sigma,\C)$ induced by inclusion. Since $i^*\circ \sigma^* = i^*$ and the coinvariant ring projection is the coequaliser of $\sigma^*$ and the identity, $i^*$ automatically factors through a map $\underline{i^*}\colon H^*_{G(\C)}(X,\C)_{\sigma^*}\to H^*_{G(\C)^{\sigma\circ}}(X^\sigma,\C)$. Moreover, by Corollary \ref{cor:base_ring_coinvariants}, we may view $\underline{i^*}$ as a map of graded modules over the graded ring $H^*_{G(\C)}(\pt, \C)_{\sigma^*}$. We now aim to apply the following algebraic lemma (see \S\ref{subsubsec:generalities} for graded algebra notions):
\begin{lem}
\label{lem:graded_Nakayama_iso}
    Let $A=\bigoplus_{n\ge 0}A^n$ be a connected graded $\C$-algebra with augmentation $\epsilon$, and $f\colon M\to N$ be a degree-preserving morphism of graded $A$-modules. Suppose that 
    \begin{itemize}
        \item[(a)] $N$ is free, and
        \item[(b)] $f\otimes_A 1 \colon M\otimes_A A/\ker(\epsilon)\to N\otimes_A A/\ker(\epsilon)$ is an isomorphism.
    \end{itemize}
    Then $f$ is an isomorphism.
\end{lem}
\begin{proof}
    Choose a $\C$-basis $\{m_i'\}_{i\in I}$ of the $\C$-vector space $M\otimes_A A/\ker(\epsilon)\cong M/\ker(\epsilon)\cdot M$ and, for each $i$, a homogeneous preimage $m_i\in M$ of $m_i'$. By the graded Nakayama lemma, the set $\{m_i\}_{i\in I}$ generates $M$, and we claim that it is in fact a basis. To see this, we note that $f\otimes_A A/\ker(\epsilon)$ sends the $m_i'$ to a basis of $N\otimes_A A/\ker(\epsilon)$ by assumption (b). Another application of the graded Nakayama lemma and assumption (a) show that $\{f(m_i)\}_{i\in I}$ is a basis of $N$. Any $A$-linear relation of the generating set $\{m_i\}_{i\in I}$ then must become trivial upon applying $f$, hence has to be trivial in the first place. Therefore, $\{m_i\}_{i\in I}$ is indeed a basis of $M$ which $f$ sends bijectively to a basis of $N$, so $f$ is indeed an isomorphism. 
\end{proof}

Intending to apply the lemma with $A = H^*_{G(\C)}(\pt,\C)_{\sigma^*}\cong H^*_{G(\C)^{\sigma\circ}}(\pt, \C)$ and $f = \underline{i^*}$, we have to check the assumptions (a) and (b). To this end, we first record a simple algebraic observation, regarding quotients of coinvariants:
\begin{prop}
\label{prop:coinv_and_quotients}
    Let $\phi$ be an automorphism of an algebra $A\to B$. (That is, $A$, $B$ are rings, and $\phi$ is to consist of compatible ring automorphisms of $A$ and $B$, both denoted $\phi$.) Suppose that $I\subseteq A$ is an ideal such that $\phi(I) = I$. Then $\phi$ descends to an automorphism $\phi_I$ of the algebra $A/I \to B/(B\cdot I)$, $I$ descends to an ideal $I_\phi$ of the coinvariant ring $A_\phi$, and we have compatible isomorphisms
    \begin{equation*}
        A_\phi/I_\phi \cong (A/I)_{\phi_I},
        \quad
        B_\phi/(I_\phi B_\phi) \cong (B/(I\cdot B))_{\phi_I}.
    \end{equation*}
\end{prop}
\begin{proof}
    The algebra automorphism $\phi_I$ is defined uniquely by the commutative diagrams
    \begin{equation*}
    \begin{tikzcd}
        A 
            \arrow[r, "\phi"]
            \arrow[d, twoheadrightarrow]
        &A
            \arrow[d, twoheadrightarrow]
        \\
        A/I 
            \arrow[r, "\phi_I"]
        &A/I
    \end{tikzcd}
    \quad
    \text{and}
    \quad
    \begin{tikzcd}
        B 
            \arrow[r, "\phi"]
            \arrow[d, twoheadrightarrow]
        &B
            \arrow[d, twoheadrightarrow]
        \\
        B/(IB) 
            \arrow[r, "\phi_I"]
        &B/(IB).
    \end{tikzcd}
    \end{equation*}
    The ideal $I_\phi$ is the image of $I$ under the surjection $A\twoheadrightarrow A_\phi$. Both $A_\phi/I_\phi$ and $(A/I)_{\phi_I}$ are quotients of $A$ by the ideal $I+(a-\phi(a)\mid a\in A)$, so they are canonically isomorphic. Similarly, both $B_\phi/(I_\phi B_\phi)$ and $(B/IB)_{\phi_I}$ are quotients of $B$ by the ideal $IB + (b-\phi(b)\mid b\in B)$, so they are also canonically isomorphic, compatibly with the structure map $A\to B$.
\end{proof}

Now, since $H^*(X,\C)$ is concentrated in even degrees, $X$ is an equivariantly formal $G(\C)$-space, so we have $H^*_{G(\C)}(X,\C)/(H^+_{G(\C)}(\pt, \C))\cong H^*(X, \C)$ (see Proposition \ref{prop:eq_formal_sp_fibre}). Proposition \ref{prop:coinv_and_quotients} then allows us to commute coinvariants with quotienting by $H^+_{G(\C)}(\pt, \C)$, yielding
\begin{equation}
\label{eq:coinv_special_fibre_step}
    H^*_{G(\C)}(X,\C)_{\sigma^*}/
    ((H^*_{G(\C)}(\pt,\C)_{\sigma^*})^+)
    \cong
    H^*(X,\C)_{\sigma^*},
\end{equation}
where the right-hand subscript $\sigma^*$ denotes coinvariants for the action of $\sigma$ on the usual cohomology. In terms of Lemma \ref{lem:graded_Nakayama_iso}, the left-hand side of \eqref{eq:coinv_special_fibre_step} corresponds to $M\otimes_A A/\ker(\epsilon)$. If we also knew $X^\sigma$ to be $G(\C)^{\sigma\circ}$-equivariantly formal, we would firstly obtain the freeness condition (a) of Lemma \ref{lem:graded_Nakayama_iso}. Secondly, we would get a variant of \eqref{eq:coinv_special_fibre_step} for $\sigma$-fixed points, which would reduce condition (b) in the lemma to an isomorphism of $H^*(X,\C)_{\sigma^*}\cong $ with $H^*(X^\sigma, \C)^{S(\kappa_\sigma)}$ or $H^*(X^\sigma, \C)$. Happily, this isomorphism is provided by Theorem \ref{thm:noneq_coinv} below and will easily imply the tentatively assumed equivariant formality.
\begin{thm}
\label{thm:noneq_coinv}
    Let $G$ be a connected complex reductive group with affine Grassmannian $\Gr = \Gr_G$.
    \begin{itemize}
        \item[(i)] If $\sigma$ is a pinning-preserving automorphism of $G$, the canonical homomorphism $H^*(\Gr(\C), \C)\to H^*(\Gr(\C)^\sigma, \C)$ factors through an isomorphism 
        \begin{equation*}
            H^*(\Gr(\C), \C)_{\sigma^*}\to H^*(\Gr(\C)^\sigma, \C)^{S(\kappa_\sigma)}.
        \end{equation*}
        \item[(ii)] If $\sigma$ is a quasi-anisotropic real structure of $G$, the canonical homomorphism $H^*(\Gr(\C)_\mathrm{adm}, \C)\to H^*(\Gr(\C)^\sigma, \C)$ factors through an isomorphism
        \begin{equation*}
            H^*(\Gr(\C)_\mathrm{adm}, \C)_{\sigma^*}\to H^*(\Gr(\C)^\sigma, \C).  
    \end{equation*}
    \end{itemize}
\end{thm}
Before giving the proof, we abstract from it two arguments related to direct products. In terms of affine schemes, they correspond to essentially evident facts about fixed points and connected components.
\begin{lem}
\label{lem:coinv_components}
    Let $A, B$ be rings equipped with an automorphism $\varphi\colon A\to A$ and a homomorphism $f\colon A\to B$. Moreover, let $I, J$ be finite sets equipped with a permutation $s\colon I\to I$ and a function $t\colon J\to I$. Consider the product rings $A^I, B^J$ equipped with the automorphism $\Phi\colon A^I\to A^I$, and homomorphism $F\colon A^I\to B^J$ given as follows on $\mathbf{a}=(a_i)_{i\in I}\in A^I$:
    \begin{equation*}
        \Phi(\mathbf{a})_i \coloneqq \varphi(a_{s(i)}),
        \quad
        F(\mathbf{a})_j \coloneqq f(a_{t(j)}),
        \quad i\in I.
    \end{equation*}
    Suppose that
    \begin{enumerate}
        \item $f\circ \varphi = f$,
        \item $t$ is injective with image $I^s$, and
        \item $f$ factors through an isomorphism $A_\varphi\to B$.
    \end{enumerate}
    Then $F$ factors through an isomorphism $(A^I)_\Phi\to B^J$.
\end{lem}
\begin{proof}
    Assumptions (1) and (2) guarantee $F\circ \Phi = F$, so $F$ certainly factors through a homomorphism $(A^I)_\Phi\to B^J$. We have to check that $F$ is surjective and that its kernel is contained in the ideal generated by the elements $\mathbf{a}-\Phi(\mathbf{a})$ for $\mathbf{a}\in A^I$. 
    
    For surjectivity, let $\mathbf{b} = (b_j)_{j\in J}\in B^J$ be arbitrary. We construct a preimage $\mathbf{a} = (a_i)_{i\in I}\in A^I$ as follows: for $i\in I^s$, there is a unique $j\in J$ such that $i=t(j)$, by assumption (2). By assumption (3), we may find $a\in A$ such that $f(a_i) = b_i$, and we define $a_i\coloneqq a$. For $i\in I\setminus I^s$, we may define $a_i$ arbitrarily. Then $F(\mathbf{a}) = \mathbf{b}$ by construction, so $F$ is surjective.

    Now suppose $\mathbf{a}= (a_i)_{i\in I} \in \ker(F)$. If $i\in I^s$, then $f(a_i) = 0$, so assumption (2) yields elements $x_{i,1},\ldots, x_{i,n_i}$ and $y_{i,1},\ldots, y_{i,n_i}$ of $A$ such that $a_i = \sum_{k=1}^{n_i} (x_{i,k} - \varphi(x_{i,k}))y_{i,k}$. Since $I$ is finite, the numbers $n_i$ for $i\in I^s$ are uniformly bounded by some $n\in \N$, and we can define $\mathbf{x}_k \coloneqq (x_{i,k})_{i\in I}, \mathbf{y}_k \coloneqq (y_{i,k})_{i\in I} \in A^I$, $k= 1,\ldots, n$ by setting $x_{i,k} = y_{i,k} = 0$ if $i\not \in I^s$ or $k>n_i$. For every $i\in I\setminus I^s$, we define $\mathbf{z}_i\in A^I$ to be $\mathbf{z}_i\coloneqq a_i\delta_i$, where $\delta_i\in A^I$ has $i$-coordinate 1 and all other coordinates 0. We then have
    \begin{equation*}
        \mathbf{a} = \sum_{k=1}^n (\mathbf{x}_k - \Phi(\mathbf{x}_k)) + \sum_{i\in I\setminus I^s}(\mathbf{z}_i - \Phi(\mathbf{z}_i))\cdot \delta_i,
    \end{equation*}
    which exhibits $\mathbf{a}$ as an element of $(\mathbf{c} - \Phi(\mathbf{c}) \mid \mathbf{c}\in A^I)$, as desired.
\end{proof}

\begin{lem}
\label{lem:coinv_dir_prod}
    Let $A, B$ be Noetherian rings equipped with an automorphism $\varphi\colon A\to A$ and a homomorphism $f\colon A\to B$. Let $I$ be a set, and consider the automorphism $\varphi^I\colon A^I\to A^I$ and homomorphism $f^I\colon A^I\to B^I$. If $f$ factors through an isomorphism $A_\varphi\to B$, then $f^I$ factors through an isomorphism $(A^I)_{\varphi^I}\to B^I$. 
\end{lem}
\begin{proof}
    Both $(A^I)_{\varphi^I}$ and $B^I$ are quotients of $A^I$, and the claim reduces to an equality of the ideals of $A^I$ defining them. Setting $K\coloneqq (a-\varphi(a)\mid a\in A)$, we have $\ker(f^I) = K^I$, and this clearly contains $(\mathbf{a}-\varphi^I(\mathbf{a})\mid \mathbf{a}\in A^I)$. For the converse inclusion, we note that $K$ admits a finite generating set $\{k_1,\ldots, k_m\}$ because $A$ is Noetherian, and we may assume $k_\ell = x_\ell - \varphi(x_\ell)$ for some $x_\ell\in A$, for every $\ell\in \{1,\ldots, m\}$. Thus, defining $\mathbf{x}_\ell\in A^I$ to have all entries equal to $x_\ell$, we see that every element of $K^I$ is an $A^I$-linear combination of the elements $\mathbf{x}_\ell - \varphi^i(\mathbf{x}_\ell)$, $\ell = 1,\ldots, m$. This establishes the claimed identity of ideals, finishing the proof.
\end{proof}

\begin{proof}[Proof of Theorem \ref{thm:noneq_coinv}]
    To avoid cohomology of infinite type, we may reduce to the case of semisimple $G$ as follows: it suffices to consider one $\sigma$-orbit of connected components of $\Gr_G(\C)$ at a time, and since $\sigma$ has finite order, this is a finite union of connected components. As recalled in \S\ref{subsubsec:Gr_to_sesi}, there is a map $\Gr_G\to \Gr_{G_\mathrm{ad}}\times \Gr_{G_\mathrm{ab}}$ which identifies connected components of $\Gr_G$ with connected components of $\Gr_{G_\mathrm{ad}}$. In fact, a given $\sigma$-orbit of components will land in a single copy of $\Gr_{G_\mathrm{ad}}$. This map is compatible with the action of $\sigma$, and allows us to replace $G$ with the semisimple group $G_\mathrm{ad}$, using Lemma \ref{lem:coinv_dir_prod}. Thus, we now assume $G$ to be semisimple; in particular, $\pi_0(\Gr(\C))$ is finite.

    Let us now begin with part (i), where $\sigma$ is a complex automorphism. Theorem \ref{thm:mitchell--quillen} yields a commutative diagram
    \begin{equation}
    \label{eq:noneq_coinv_diag}
    \begin{tikzcd}
        \Gr(\C)^\sigma
            \arrow[r]
            \arrow[d]
        & \Gr(\C)
            \arrow[d]
        \\
        \Omega H
            \arrow[r, "\Omega i"]
        & \Omega K,
    \end{tikzcd}
    \end{equation}
    where $K$ is a maximal compact subgroup of $G(\C)$ stable under $\sigma$, so that $H\coloneqq K^\sigma$ is also a maximal compact subgroup of $G(\C)^\theta$. The horizontal arrows are inclusions (with $i$ the inclusion $H\to K$), and the vertical arrows are homotopy equivalences. Moreover, it is easy to see that the automorphism $\sigma$ of $\Gr(\C)$ extends to the automorphism of $\Omega K$ given by postcomposition of loops with $\sigma$. 

    Next, the compact Lie groups $K$, $H$ are homotopy-simple and of finite type, so Proposition \ref{prop:coh_of_disconnected_loop_space} yields isomorphisms 
    \begin{align}
    \label{eq:noneq_coinv_cplx_step}
        &H^*(\Omega K,\C) \cong U(V_{\Tilde{K}}[1])\otimes_\C \C[\pi_1(K)]\cong U(V_{\Tilde{K}}[1])^{\pi_1(K)},
        \\
        &H^*(\Omega H,\C) \cong U(V_{\Tilde{H}}[1])\otimes_\C \C[\pi_1(H)]\cong U(V_{\Tilde{H}}[1])^{\pi_1(H)},\notag
    \end{align}
    where $\Tilde{K}$, $\Tilde{H}$ are the universal covering groups, and $V_{\Tilde{K}}$, $V_{\Tilde{H}}$ the graded vector spaces underlying their $\C$-minimal models. We have lifted homomorphisms $\Tilde{i}\colon \Tilde{H}\to \Tilde{K}$ and $\Tilde{\sigma}\colon \Tilde{K}\to \Tilde{K}$, which induce morphisms $Q(\Tilde{i})\colon V_{\Tilde{K}}\to V_{\Tilde{H}}$ and $Q(\Tilde{\sigma})\colon V_{\Tilde{K}}\to V_{\Tilde{K}}$. The isomorphisms in \eqref{eq:noneq_coinv_cplx_step} identify the canonical map $H^*(\Omega K,\C)\to H^*(\Omega H, \C)$ with
    \begin{equation*}
        F\coloneqq 
        U(Q(\Tilde{i}))\otimes (i_*)^*\colon
        U(V_{\Tilde{K}}[1])\otimes_\C \C[\pi_1(K)] \to
        U(V_{\Tilde{H}}[1])\otimes_\C \C[\pi_1(H)],
    \end{equation*}
    where $(i_*)^*$ is precomposition with $i_*\colon \pi_1(H)\to \pi_1(K)$, and the automorphism $\sigma^*$ with
    \begin{equation*}
        \Phi\coloneqq
        U(Q(\Tilde{\sigma}))\otimes (\sigma_*)^*\colon
        U(V_{\Tilde{K}}[1])\otimes_\C \C[\pi_1(K)] \to
        U(V_{\Tilde{K}}[1])\otimes_\C \C[\pi_1(K)].
    \end{equation*}
    Moreover, the $S(\kappa_\sigma)$-invariant part of $H^*(\Gr(\C)^\sigma,\C)$ is identified with the subalgebra of $U(V_{\Tilde{H}}[1])^{\pi_1(H)}$ whose elements have the same entries on any given $\kappa_\sigma$-coset. We may also identify this algebra with $U(V_{\Tilde{H}}[1])^{\pi_1(H)/\kappa_\sigma}$.

    We are now in the setting of Lemma \ref{lem:coinv_components}, with $\phi = U(Q(\Tilde{\sigma}))$, $f = U(Q(\Tilde{i}))$, $s=\sigma_*$ and $t$ the map $\pi_1(H)/\kappa_\sigma\to \pi_1(G)$ induced by $i_*$, so we turn to checking the Lemma's conditions (1)--(3). That condition (1) holds follows from $\Tilde{\sigma}\circ \Tilde{i} = \Tilde{i}$ and functoriality. Condition (2) follows from $\kappa_\sigma = \ker i_*$ and Lemma \ref{lem:pi1_image}. As for condition (3), we claim that $Q(\Tilde{i})$ identifies with the projection of $V_{\Tilde{K}}$ to $V_{\Tilde{K}}^{Q(\Tilde{\sigma})}$. Since $U(V_{\Tilde{K}}[1])$ and $U(V_{\Tilde{H}}[1])$ are polynomial rings, this will immediately yield the desired isomorphism of the condition.  
    
    To now prove this claim, note that $\Tilde{K}$ and $\Tilde{H}$ are connected, simply connected formal spaces with free cohomology algebras. Lemma \ref{lem:formal_min_model} therefore yields $V_{\Tilde{K}}\cong Q(H^*(\Tilde{K}, \C)) \cong P_{\Tilde{K}}$ and $V_{\Tilde{H}}\cong Q(H^*(\Tilde{H}, \C))\cong P_{\Tilde{H}}$, with notation as in \S\ref{subsec:prelim_Hsp}. Transgression then identifies $P_{\Tilde{K}}[1]\cong Q_{B\Tilde{K}}$ and $P_{\Tilde{H}}[1]\cong Q_{B\Tilde{H}}$, and further identifies them with the generators of the graded polynomial rings $\C[\kay]^{\Tilde{K}}$ and $\C[\h]^{\Tilde{H}}$ (\cf diagram \eqref{eq:transgression_naturality_diagram}). Following the maps $Q(\Tilde{\sigma})$ and $Q(\Tilde{i})$ through these steps, we see that they correspond to the usual $\sigma^*\in \Aut(\C[\kay]^{\Tilde{K}})$ and restriction $\C[\kay]^{\Tilde{K}}\to \C[\h]^{\Tilde{H}}$, respectively. This reduces the claim to Lemma \ref{lem:inv_poly_coinvariants}, and finishes the verification of condition (3). We may now apply Lemma \ref{lem:coinv_components} to infer part (i) of the theorem.

    Now consider part (ii), where $\sigma$ is a quasi-anisotropic real structure. Here the overall argument runs parallelly to the previous case, but some details differ. Firstly, the role of diagram \eqref{eq:noneq_coinv_diag} is played by
    \begin{equation}
    \begin{tikzcd}
        \Gr(\C)^\sigma
            \arrow[r]
            \arrow[d]
        & \Gr(\C)_\mathrm{adm}
            \arrow[d]
        \\
        \Omega X
            \arrow[r, "\Omega C"]
        & (\Omega K)_\mathrm{adm},
    \end{tikzcd}
    \end{equation}
    where $K$ and $H$ are as above, $C$ denotes the Cartan embedding $X\coloneqq K/H\to K$, and $(\Omega K)_\mathrm{adm}$ is the union of connected components of $\Omega X$ corresponding to admissible components of $\Gr(\C)$. By Theorem \ref{thm:mitchell--quillen} and Proposition \ref{prop:inclusion_to_loop_Cartan}, this is again a commutative diagram whose vertical arrows are homotopy equivalences. We then argue analogously to the previous case to identify $(\Omega C)^*$ with 
    \begin{equation*}
        F'\coloneqq 
        U(Q(\Tilde{C}))\otimes (C_*)^*\colon
        U(V_{\Tilde{K}}[1])\otimes_\C \C[\pi_0((\Omega K)_\mathrm{adm})] \to
        U(V_{\Tilde{X}}[1])\otimes_\C \C[\pi_0(\Omega X)],
    \end{equation*}
    where $\Tilde{X} = \Tilde{K}/\Tilde{K}^{\Tilde{\sigma}}$ is the universal covering of $X$, $\Tilde{C}\colon \Tilde{X}\to \Tilde{K}$ its Cartan embedding (lifting $C$), and $\sigma^*$ with
    \begin{equation*}
        \Phi'\coloneqq
        U(Q(\Tilde{\sigma}))\otimes (\sigma_*)^*\colon
        U(V_{\Tilde{K}}[1])\otimes_\C \C[\pi_0((\Omega K)_\mathrm{adm})] \to
        U(V_{\Tilde{K}}[1])\otimes_\C \C[\pi_0((\Omega K)_\mathrm{adm})].
    \end{equation*}

    We now again aim to apply Lemma \ref{lem:coinv_components}, this time with $\phi = U(Q(\Tilde{\sigma}))$, $f=U(Q(\Tilde{C}_*))$, $s=\sigma_*$ and $t=C_*$. For condition (1), note that $\Omega \Tilde{C}$ is homotopic to the map sending a loop $\gamma\in \Omega \Tilde{X}$ to $\gamma\ast (\Tilde{\sigma}\circ\gamma^-)$, where $\ast$ denotes concatenation and the superscript minus denotes reversal. Similarly, $\Tilde{\sigma}$ acts on $\Omega K$ by a map homotopic to $\gamma\mapsto \sigma\circ \gamma^-$. It follows that $\Tilde{\sigma} \circ \Omega \Tilde{C}$ is homotopic to $\Tilde{C}$, which implies $Q(\Tilde{C}) \circ Q(\Tilde{\sigma}) = Q(\Tilde{C})$. Condition (2) follows from Lemma \ref{lem:pi0_real_Gr}. Similarly as in the complex case, we may reduce condition (3) to the following claim, proved below: $Q(\Tilde{C})$ identifies with the projection of $V_{\Tilde{K}}[1]$ to $V_{\Tilde{K}}[1]^{\Tilde{\sigma}}$.

    Before proving this claim, let us compare the natural actions of $\Tilde{\sigma}$ on $\pi_*(\Omega \Tilde{K})$ and $\pi_*(\Tilde{K})$. Namely, $\Tilde{\sigma}$ acts on $\Omega \Tilde{K}$ by postcomposition with $\sigma$ and precomposition with complex conjugation. Due to this complex conjugation, the isomorphism 
    \begin{equation*}
        V_{\Omega\Tilde{X}}
        \cong
        (\pi_*(\Omega\Tilde{X})\otimes_\Z \C)^*
        \cong
        (\pi_*(\Tilde{X})\otimes_\Z \C)^*[1]
        \cong
        V_{\Tilde{X}}[1]
    \end{equation*}
    implicit in the previous step identifies the action of $\Tilde{\sigma}$ on $V_{\Tilde{X}}[1] = V_{\Omega \Tilde{X}}$ with \emph{minus} the natural action of $\Tilde{\sigma}$ on $V_{\Tilde{X}}$.
    
    Now, since $\sigma$ is quasi-anisotropic, so is $\Tilde{\sigma}$, whence $\Tilde{K}^{\Tilde{\sigma}}$ is nonhomologous to zero in $\Tilde{K}$ by Proposition \ref{prop:qa_nhz}(i). It then follows from Corollary \ref{cor:nhz_formal} that $\Tilde{X}$ is $\C$-formal with free cohomology algebra, so we get $V_{\Tilde{X}}\cong Q(H^*(\Tilde{X}, \C))$ in addition to the already established isomorphism $V_{\Tilde{K}}\cong P_{\Tilde{K}}$. In Proposition \ref{prop:nhz_Cartan}, the restriction of $\Tilde{C}^*$ to $P_{\Tilde{K}}$ is identified (up to a nonzero scalar) with the projection to the $-1$-eigenspace for the natural action of $\Tilde{\sigma}$ on $P_{\Tilde{K}}$. As we have just noted, this is the $+1$-eigenspace for the $\Tilde{\sigma}$-action on $P_{\Tilde{K}}[1]$ arising from the based loop space setting. This establishes the claim, so an application of Lemma \ref{lem:coinv_components} finishes the proof.
\end{proof}

We can now collect the results of this subsection to prove its main theorem:
\begin{proof}[Proof of Theorem \ref{thm:coinv_global}]
    Part (i) is Corollary \ref{cor:base_ring_coinvariants}. For (ii) and (iii), Theorem \ref{thm:noneq_coinv} implies that the cohomology of $\Gr(\C)^\sigma$ is concentrated in even degrees, since this is the case for $\Gr(\C)$. In particular, it is equivariantly formal as a $G(\C)^{\sigma\circ}$-space, so its equivariant cohomology is free (see Proposition \ref{prop:eq_formality}). This means that condition (a) of Lemma \ref{lem:graded_Nakayama_iso} applies to our situation, and its condition (b) is equivalent to Theorem \ref{thm:noneq_coinv} via Proposition \ref{prop:coinv_and_quotients}. The remaining claims of Theorem \ref{thm:coinv_global} therefore follow from Lemma \ref{lem:graded_Nakayama_iso}. 
\end{proof}

\subsection{Fixed loci of affine Schubert varieties}
\label{subsec:asV_fixed_loci}
Let $\sigma$ be as in the previous subsection, that is, either a pinning-preserving automorphism or a quasi-anisotropic real structure of a connected complex reductive group $G$. We use the standard setup of \S\ref{subsubsec:cplx_red_setup}, and in the real structure case also the refined setup of \S\ref{subsubsec:real_red_setup}. In the complex automorphism case, we instead assume that the maximal torus $T$ and Borel subgroup $B$ are part of a pinning preserved by $\sigma$. 

Now, consider a dominant coweight $\lambda\in \Lambda^\vee$. The corresponding affine Schubert variety $\Gr_G^{\le \lambda}(\C)$ is preserved by $\sigma\in \Aut(\Gr_G(\C))$ if and only if $\sigma(\lambda)$ lies in the Weyl group orbit of $\lambda$. In the case of a pinning-preserving automorphism, this is equivalent to $\sigma(\lambda) = \lambda$ (because $\sigma$ then preserves $\Lambda^\vee_+$), but this is not so for real structures. In any case, it is clear that $\Gr_G(\C)^\sigma$ is a union of points of affine Schubert varieties for $G^\sigma$ -- of complex points in the complex automorphism case, and of real points in the real structure case. The purpose of this subsection is to characterise when this union consists of only a single affine Schubert variety.  

If $\sigma$ is a pinning-preserving automorphism, then compatibility with $\lambda\in \Lambda^\vee_+$ entails $\sigma(\lambda) = \lambda$, so $(\Gr_G(\C)^{\le \lambda})^\sigma$ contains at least the affine Schubert variety $\Gr_{G^\sigma}^{\le \lambda}(\C)$. It is then clear from the closure relations of affine Schubert cells that any $G^\sigma$-affine Schubert varieties contained in $(\Gr_G(\C)^{\le \lambda})^\sigma$ but not $\Gr_{G^\sigma}^{\le \lambda}(\C)$ must lie in a different connected component. Therefore, Lemma \ref{lem:pi1_kernel} implies:
\begin{prop}
\label{prop:cplx_fixed_locus}
    Let $\sigma$ be a pinning-preserving automorphism such that $\Lie(G)$ contains no copies of $\esl_{2n+1}(\C)$, $n\in\Z_{>0}$, on which $\sigma$ acts nontrivially. Then, we have $(\Gr_G(\C)^{\le \lambda})^\sigma = \Gr_{G^\sigma}^{\le \lambda}(\C)$ for all $\lambda\in (\Lambda^\vee_+)^\sigma$.
\end{prop}

For quasi-anisotropic real structures, the result turns out to be cleaner, even though the proof -- given in Appendix \ref{appx:case_proofs} -- is more cumbersome:
\begin{prop}
\label{prop:unique_max_fixed_cwt}
    Let $\sigma$ be a quasi-anisotropic real structure of $G$, and $\lambda\in \Lambda^\vee_+$ with $\sigma(\lambda) \in W\cdot \lambda$. Write $(\Lambda^\vee_+)_{\le \lambda}$ for the set of dominant coweights less than or equal to $\lambda$. Then the partially ordered set $(\Lambda^\vee_+)_{\le \lambda}\cap \Lambda^\vee_\R$ is either empty or contains a unique maximal element $\mu_\mathrm{max}$. It follows that $\Gr_G^{\le\lambda}(\C)^\sigma$ is either empty (when its component is nonadmissible) or equal to $\Gr_{G^\sigma}^{\le \mu_\mathrm{max}}(\R)$. 
\end{prop}

\subsection{Coinvariants for affine Schubert varieties and medium algebras}
\label{subsec:coinv_aSv}
Our main result, stated below, uses the isomorphism 
\begin{equation}
\label{eq:MA_geo_repeated}
    \M^\lambda(\g) \cong H^*_{\Check{G}(\C)}(\Gr_{\Check{G}}^{\le \lambda}(\C),\C)
\end{equation}
of Proposition \ref{prop:MA_eq_coh}(ii) to describe coinvariants of the medium algebra automorphisms $\iota_{a, \zeta}$ of Definition \ref{def:autos_via_out}. We shall use the following notation:
\begin{notn}
    If $M$ is a module (or algebra) over an integral domain $A$, we write $M_\mathrm{tf}$ for its torsion-free quotient module (resp.\,algebra), that is, the quotient of $M$ by its submodule (resp.\,ideal) of torsion elements.
\end{notn}
\begin{thm}
\label{thm:coinv_main}
    Let $G$ be a connected complex reductive group and $\lambda$ a dominant weight of $G$. Consider an automorphism $\iota\in \Aut(\M^\lambda(\g))$ defined in one of the following two ways:
    \begin{itemize}
        \item[(A)] either $\iota \coloneqq \iota_{a,1}$ for $a\in \Out(G)$ compatible with $\lambda$,
        \item[(B)] or $\iota \coloneqq \iota_{a, -1}$ for $a\in \Aut(G)$ compatible with $\lambda$ such that $a$ has order at most 2.
    \end{itemize}
    Define $\sigma\colon \Check{G}\to \Check{G}$ to be
    \begin{itemize}
        \item a pinning-preserving automorphism of $\Check{G}$ whose image in $\Out(\Check{G})\cong \Out(G)$ matches $a\inv$ in case (A), or
        \item a quasi-anisotropic real structure of $\Check{G}$ whose action on the root datum matches that of $a$ ($=a\inv$) in case (B).
    \end{itemize}
    Then $\sigma$ preserves the $\Check{G}$-affine Schubert variety labelled by $\lambda$, so it induces an automorphism $\sigma^*$ on the equivariant cohomology algebra of $\Gr_{\Check{G}}^{\le \lambda}(\C)$, and the isomorphism \eqref{eq:MA_geo_repeated} identifies $\iota$ with $\sigma^*$. 
    
    Assume further that $\Lie(\Check{G})$ contains no copy of $\esl_{2n+1}(\C)$, $n \in \Z_{>0}$, on which
    \begin{itemize}
        \item some power of $\sigma$ restricts to a nontrivial automorphism, in case (A), or
        \item $\sigma$ restricts to a split real structure, in case (B).
    \end{itemize}
    Then \eqref{eq:MA_geo_repeated} and the equivariant cohomology maps induced by inclusions yield isomorphisms
    \begin{equation}
    \label{eq:coinv_main_base}
        \C[\g]^G_\iota 
        \cong 
        H^*_{\Check{G}(\C)}(\pt,\C)_{\sigma^*}
        \cong
        H^*_{\Check{G}(\C)^{\sigma\circ}}(\pt,\C)
    \end{equation}
    and
    \begin{equation}
    \label{eq:coinv_main_total}
        \big(\M^\lambda(\g)_\iota\big)_\mathrm{tf}
        \cong 
        \big(H^*_{\Check{G}(\C)}(\Gr_{\Check{G}}^{\le \lambda}(\C),\C)_{\sigma^*}\big)_\mathrm{tf}
        \cong
        H^*_{\Check{G}(\C)^{\sigma\circ}}(\Gr_{\Check{G}}^{\le \lambda}(\C)^\sigma,\C)
    \end{equation}
    where torsion-free quotients are taken with respect to the coinvariant base ring $\C[\g]^G_\iota \cong H^*_{\Check{G}(\C)}(\pt, \C)_{\sigma^*}$. 
\end{thm}
The proof will use the following lemma.\footnote{The idea of using this lemma in the present context was inspired by \cite[Sec.\,8]{Zveryk2024_dynkin}.}
\begin{lem}
\label{lem:coinv_torsion_lemma}
    Let $A$ be a Noetherian integral domain and $f\colon M\to N$ a morphism of finitely generated $A$-modules. Suppose that
    \begin{itemize}
        \item[(a)] $N$ is free,
        \item[(b)] $f$ is surjective, and
        \item[(c)] there exists a maximal ideal $\m$ of $A$ such that $f\otimes_A 1\colon M\otimes_A A/\m \to N\otimes_A A/\m$ is an isomorphism.
    \end{itemize}
    Then $f$ factors through an isomorphism $M_\mathrm{tf}\to N$, where $M_\mathrm{tf}$ is the quotient of $M$ by its submodule of torsion elements.
\end{lem}
\begin{proof}
    Assumptions (a) and (b) imply that $f$ factors through a surjection $M_\mathrm{tf}\to N$, so it suffices to show that $M'\coloneqq \ker f$ consists of torsion elements. 
    
    Now, if $d$ is the minimum of the dimensions $\dim_{A/\mathfrak{n}}(M\otimes_A A/\mathfrak{n}M)$ over all maximal ideal $\mathfrak{n}$ of $A$, then by semicontinuity the fibre dimension $d$ is attained on a nonempty open subset $U$ of $\Spec A$. By assumption (b), we have $d\ge \operatorname{rk}_A N$, and so assumption (c) implies $d=\operatorname{rk}_A N$. It follows that $f$ induces an isomorphism $M/\mathfrak{n} M\to N/\mathfrak{n} N$ for all maximal ideals $\mathfrak{n}$ in $U$. 
    
    Since $A$ is Noetherian and $M$ finitely generated, $M'$ is also finitely generated. Therefore, we may use Nakayama's lemma and assumption (c) to show that $U$ is disjoint from the support of $M'$. In particular, there exists $f\in A$ such that the basic open subset $D(f)\subseteq \Spec A$ is disjoint from $\operatorname{Supp}M'$, whence $M'_f = 0$. This shows that all elements of $M'$ are annihilated by some power of $f$, as required.
\end{proof}
\begin{proof}[Proof of Theorem \ref{thm:coinv_main}]
    That \eqref{eq:MA_geo_repeated} takes $\iota$ to $\sigma^*$ was shown in Theorem \ref{thm:MA_autos_geometry}. It yields an isomorphism of coinvariant algebras 
    \begin{equation*}
        (\C[\g]^G_\iota \to \M^\lambda(\g)_\iota)
        \cong
        (H^*_{\Check{G}(\C)}(\pt,\C)_{\sigma^*} \to H^*_{\Check{G}(\C)}(\Gr_{\Check{G}}^{\le \lambda}(\C),\C)_{\sigma^*}\big),
    \end{equation*}
    which of course implies an isomorphism of their torsion-free quotients. We can now restrict our attention to the equivariant cohomology setting. There, the base ring isomorphism $H^*_{\Check{G}(\C)}(\pt,\C)_{\sigma^*} \cong H^*_{\Check{G}(\C)^{\sigma\circ}}(\pt,\C)$ has been established in Corollary \ref{cor:base_ring_coinvariants}. It remains to establish the isomorphism 
    \begin{equation*}
        \big(H^*_{\Check{G}(\C)}(\Gr_{\Check{G}}^{\le \lambda}(\C),\C)_{\sigma^*}\big)_\mathrm{tf} \overset{!}{\cong} H^*_{\Check{G}(\C)^{\sigma\circ}}(\Gr_{\Check{G}}^{\le \lambda}(\C)^\sigma,\C),
    \end{equation*}
    for which we intend to use Lemma \ref{lem:coinv_torsion_lemma}. 

    Write $X\coloneqq \Gr_{\Check{G}}(\C)$ and $X^\lambda \coloneqq \Gr_{\Check{G}}^{\le \lambda}(\C)$. Then we have the following commutative diagram of maps induced by inclusions:
    \begin{equation*}
    \begin{tikzcd}
        H^*_{\Check{G}(\C)}(X,\C) \arrow[r] \arrow[d]
        & H^*_{\Check{G}(\C)^{\sigma\circ}}(X^\sigma,\C) \arrow[d]\\
        H^*_{\Check{G}(\C)}(X^\lambda,\C)\arrow[r, "r"]
        & H^*_{\Check{G}(\C)^{\sigma\circ}}((X^\lambda)^\sigma,\C).
    \end{tikzcd}
    \end{equation*}
    In this diagram, the top horizontal map is surjective by Theorem \ref{thm:coinv_global} -- note that in the complex case (A), this uses the additional assumption and Lemma \ref{lem:pi1_kernel}, which imply that $\kappa_\sigma$ is trivial. The left vertical map is surjective by Lemma \ref{lem:eq_coh_restriction_surjective} and the surjectivity of $H^*(X,\C)\to H^*(X^\lambda,\C)$, obtained from compatible complex affine pavings. Thus, the map labelled $r$ is surjective if and only if the right vertical map is. By another application of Lemma \ref{lem:eq_coh_restriction_surjective}, this is the case if 
    \begin{equation*}
        r'\colon H^*(X^\sigma,\C)\to H^*((X^\lambda)^\sigma,\C) 
    \end{equation*}
    is surjective. In the complex case (A), surjectivity of $r'$ again comes from compatible complex affine pavings; in the real case (B), it follows from Proposition \ref{prop:even_Iwahori_dimensions}. Altogether, we conclude that $r$ is surjective.

    That $r'$ is surjective implies also (together with Theorem \ref{thm:noneq_coinv}) that $H^*((X^\lambda)^\sigma,\C)$ is concentrated in even degrees, so that $H^*_{\Check{G}(\C)^{\sigma\circ}}((X^\lambda)^\sigma,\C)$ is free by Proposition \ref{prop:eq_formality}. The map $r$ therefore satisfies conditions (a) and (b) of Lemma \ref{lem:coinv_torsion_lemma}, and it remains only to check condition (c).

    To that end, it suffices to find a maximal ideal $\m$ of $\C[\g]^G$, such that the dimension of $(\M^\lambda(\g)_\iota)/(\m\cdot (\M^\lambda(\g)_\iota))$ is equal to the rank of the free module $H^*_{(\Check{G})^{\sigma\circ}}((X^\lambda)^\sigma, \C)$, that is, to the dimension of $H^*((X^\lambda)^\sigma,\C)$. A basis of $H^*(X^\lambda(\C),\C)$ is given by classes of Iwahori cells (more precisely, the cell closures), see Proposition \ref{prop:even_Iwahori_dimensions}, so the dimension in question equals the cardinality of $\wt(\lambda)^\sigma$. 
    
    In the case (A) of a pinning-preserving automorphism, the sought-after maximal ideal (or closed point) is provided by Proposition \ref{prop:MA_weight_action}(ii). In case (B), some care is now needed regarding the choice of maximal torus. Namely, by the standard setup, the maximal torus $\Check{T}$ of $\Check{G}$ should be maximally split with respect to $\sigma$. This is relevant here for the identification of $\wt(\lambda)^\sigma$ with a basis of the singular cohomology of $(X^\lambda)^\sigma$, so we cannot easily change $\Check{T}$. However, by Proposition \ref{prop:qa_max_tori}, for quasi-anisotropic real reductive groups without any split type $A_{2n}$ factors, all maximal tori are conjugate. Thus, our additional assumption in case (B) guarantees that $\Check{T}$ can be chosen freely among $\sigma$-stable maximal tori. 
    
    In particular, we may assume that $\Check{T}$ is part of a pinning preserved by the Cartan involution $\theta$ of $\sigma$. Recall that the actions of $\theta$ and $\sigma$ on coweights differ by a sign. In view of Proposition \ref{prop:MA_weight_action}, we want to now identify the action of $\sigma$ on coweights (of $\Check{G}$) with the action of $\eta_-$ on weights (of $G$), where $\eta_-$ is an involutive automorphism of $G$ preserving $\te$ and taking a regular element $h$ to $-h$. To do so, let $\theta'$ be a pinning-preserving automorphism of $G$ dual to $\theta$, and $\omega\in \Aut(G)$ a Chevalley involution (\cf\cite[Sec.\,2]{AV2016_contragradient}) commuting with $\theta'$, both assumed to preserve $\te$. Then, on the one hand, the actions on weights of $\theta'\omega$ and $\theta'$ differ by a sign, and on the other hand $\theta'\omega$ takes a regular element of $\te$ to its negative, so it can be chosen for $\eta_-$. Combining these observations gives the desired identification of (co)weight actions, so Proposition \ref{prop:MA_weight_action}(iii) yields a suitable closed point. This finishes the proof.
\end{proof}

\begin{rmk}
    The additional assumption in the second half of Theorem \ref{thm:coinv_main} looks similar in cases (A) and (B), but is used for different purposes. Indeed, in case (A) it prevents a disconnected fixed locus (see Proposition \ref{prop:cplx_fixed_locus}). This is indeed necessary, since having several connected components in the fixed locus would make its degree zero cohomology too large. 
    
    In case (B), the assumption excludes occurrences of the split quasi-anisotropic group $\PGL_{2n+1,\R}$, or more generally those groups whose adjoint group contains such a factor. This exclusion is used essentially in two places in the proof. Firstly, it ensures that maps in cohomology from real affine Grassmannians to their affine Schubert varieties are surjective. This may in fact also work in the excluded case (\cf Conjecture \ref{conj:typeA2_Gr_cohomology}). Secondly, excluding $\PGL_{2n+1,\R}$ ensures that all real maximal tori are conjugate (through real points), which is certainly not the case for $\PGL_{2n+1,\R}$ with $n\ge 1$. To see that this distinction matters, we note that its real structure acts trivially on the coweights with respect to a split maximal torus, so the last step of the proof would have to identify this trivial action with the action of $\iota$ on a medium algebra fibre. But the corresponding automorphism $\iota = \iota_{1,-1}$ does not (generally) act trivially on any fibre. Despite these two issues, we expect that Theorem \ref{thm:coinv_main} holds in case (B) without the additional assumption. In fact, inspecting the proof one sees that this would follow from Conjecture \ref{conj:typeA2_Gr_cohomology}.
\end{rmk}

\subsection{Specialisation of coinvariants and torsion}
\label{subsec:coinv_global_to_aSv}
The main content of Theorem \ref{thm:coinv_main} is a description of coinvariants for the equivariant cohomology of an affine Schubert variety $\Gr_{\Check{G}}^{\le \lambda}$. Since the analogous description for the full affine Grassmannian is known from Theorem \ref{thm:coinv_global}, the straightforward approach would be the following: we have 
\begin{equation}
\label{eq:coinv_sp_I_lambda}
    H^*_{\Check{G}(\C)}(\Gr_{\Check{G}}^{\le \lambda}(\C),\C)
    =
    H^*_{\Check{G}(\C)}(\Gr_{\Check{G}}(\C)_\lambda,\C) / I_\lambda
\end{equation}
where the subscript $\Gr_{\Check{G}}(\C)_\lambda$ denotes the connected component labelled by $\lambda$, and $I_\lambda\unlhd H^*_{\Check{G}(\C)}(\Gr_{\Check{G}}(\C),\C)$ is some ideal. Clearly, $I_\lambda$ is stable under $\sigma^*$ if $\sigma$ preserves $\Gr_{\Check{G}}^{\le \lambda}$. Since coinvariants commute with quotients (\cf Proposition \ref{prop:coinv_and_quotients}), Theorem \ref{thm:coinv_global} implies
\begin{equation*}
    H^*_{\Check{G}(\C)}(\Gr_{\Check{G}}^{\le \lambda}(\C),\C)_{\sigma^*}
    \cong 
    H^*_{\Check{G}(\C)^{\sigma \circ}}(\Gr_{\Check{G}}(\C)_\lambda^\sigma,\C)_{\sigma^*} / (I_\lambda)_{\sigma^*}
\end{equation*}
where $(I_\lambda)_{\sigma^*}$ is the image of $I_\lambda$ under the map $H^*_{\Check{G}(\C)}(\Gr_{\Check{G}}(\C),\C)\to H^*_{\Check{G}(\C)^{\sigma \circ}}(\Gr_{\Check{G}}(\C)^\sigma,\C)$. On the other hand, we also have
\begin{equation}
\label{eq:coinv_sp_2}
    H^*_{\Check{G}(\C)^{\sigma\circ}}(\Gr_{\Check{G}}^{\le \lambda}(\C)^\sigma,\C)
    =
    H^*_{\Check{G}(\C)^{\sigma \circ}}(\Gr_{\Check{G}}(\C)_\lambda^\sigma,\C) / J_\lambda
\end{equation}
for some ideal $J_\lambda$. Thus, if $(I_\lambda)_{\sigma^*}$ was known to equal $J_\lambda$, we could infer a coinvariant algebra description as in Theorem \ref{thm:coinv_main}, but even \emph{without taking the torsion-free quotient}.

This cannot hold in general, since the coinvariant algebra really can have torsion, as we have seen in Examples \ref{expl:MA_A1}--\ref{expl:A3_2nd_fd}. Nevertheless, calculations in the real case suggest the following:
\begin{conj}
\label{conj:fixed_lambda_1}
    With notation as above, suppose $\sigma$ is a quasi-anisotropic real structure and $\Lie(\Check{G}^\sigma)$ does not contain a copy of $\esl_{2n+1}(\R)$ for any $n\in \Z_{>0}$. If $\sigma(\lambda) = \lambda$, then $(I_\lambda)_{\sigma^*}\cong J_\lambda$, and so we have the following strengthened version of \eqref{eq:coinv_main_total}:
    \begin{equation}
    \label{eq:conj_fixed_lambda}
        \M^\lambda(\g)_\iota
        \cong 
        H^*_{\Check{G}(\C)}(\Gr_{\Check{G}}^{\le \lambda}(\C),\C)_{\sigma^*}
        \cong
        H^*_{\Check{G}(\C)^{\sigma\circ}}(\Gr_{\Check{G}}^{\le \lambda}(\C)^\sigma,\C)
        =
        H^*_{\Check{G}^\sigma(\R)^\circ}(\Gr_{\Check{G}^\sigma}^{\le \lambda}(\R),\C).
    \end{equation}
\end{conj}
In \eqref{eq:conj_fixed_lambda}, the first isomorphism and the last identity hold always (the latter by Proposition \ref{prop:unique_max_fixed_cwt}). As for the second conjectured isomorphism, we always have at least a surjection to the free module on the right-hand-side. Since both sides of that surjection are graded modules, it would suffice to check the second isomorphism on the special fibres. One arrives this way at the following equivalent form of the conjecture:
\begin{conj}
\label{conj:fixed_lambda_2}
    With notation as in Conjecture \ref{conj:fixed_lambda_1}, let $\Gr\coloneqq \Gr_{\Check{G}^\sigma}$. If $\sigma(\lambda) = \lambda$, then
    \begin{equation*}
        H^*(\Gr^{\le \lambda}(\R),\C)
        \cong
        H^*(\Gr^{\le \lambda}(\C),\C)
        \otimes_{H^*(\Gr(\C),\C)}
        H^*(\Gr(\R),\C).
    \end{equation*}
\end{conj}

\section{Automorphisms of moduli spaces of Higgs bundles}
\label{sec:Higgs}
The original motivation for Hausel's study of medium algebras \cite{Hausel2024_avatars} was a conjectural relation to Higgs bundles (\cf \cite[Sec.\,3]{Hausel2023_ICM_mirror}). An established special case of this relation is explained in \cite[Sec.\,4]{GH2024_even_very_stable}, where it was used to describe certain fixed loci in Higgs bundle moduli spaces. In this section, we recall some details of this special case, communicated to the author by Miguel González, with the intention of connecting our medium algebra automorphisms to automorphisms of moduli spaces of Higgs bundles, extending \cite[Sec.\,4]{GH2024_even_very_stable}.

We first recall the setup, following \cite{González2025_vs_regular} throughout. Thus, let $G$ be a connected semisimple complex algebraic group, $C$ a smooth projective complex curve of genus at least 2 with canonical line bundle $K_C$, and $\mathbb{M}(G)$ the moduli space of polystable $G$-Higgs bundles on $C$. It has a $\C^\times$-action scaling the Higgs fields, and for a smooth $\C^\times$-fixed point $\E\in \mathbb{M}(G)$, one can consider its \emph{upward flow} $W_\E^+$, consisting of those $\E'\in \mathbb{M}(G)$ with $\lim_{z\to 0}z\cdot \E' = \E$. One further has the $\C^\times$-equivariant \emph{Hitchin map} $h\colon \mathbb{M}(G)\to \mathbb{A}_G$ to an affine space $\mathbb{A}_G$, which admits a section $s\colon \mathbb{A}_G\to \mathbb{M}(G)$. To be precise, we use here a construction of $s$ based on a principal $\esl_2$-triple (see \cite[Ex.\,3.27]{González2025_vs_regular}). The image $\E_0$ of the origin under $s$ is a smooth $\C^\times$-fixed point called the \emph{canonical uniformising Higgs bundle}. 

One way of constructing new elements of $\mathbb{M}(G)$ from old ones is via \emph{Hecke transformations} (\cf\cite[Sec.\,4.2]{HauselHitchin2022_very-stability} and \cite[Sec.\,4]{González2025_vs_regular}), which change the underlying principal $G$-bundle of a Higgs bundle locally at a point $c\in C$. For fixed $c$, the Hecke transformations of a Higgs bundle $\E$ at $c$ can be parametrised by an affine Springer fibre $\Gr_{G,X}$, where $X\in \g\otimes_\C\C((t))$ comes from the Higgs field locally around $c$. Here $\Gr_{G,X}$ is given on complex points by the set of $[g]\in LG(\C)/L^+G(\C)$ such that $g\inv \cdot X \in \g\otimes_\C \C[[t]]$. By allowing only parameters which lie inside an affine Schubert cell $\Gr_G^\lambda$ for some coweight $\lambda$, one obtains a notion of Hecke transformation of type $\lambda$. 

Now fix a point $c\in C$ and let $\lambda$ be a \emph{minuscule coweight} of $G$. Consider the following commutative square (\cf \cite[Eq.\,(3.8)]{Hausel2023_ICM_mirror}):
\begin{equation}
\label{eq:Hecke_diagram}
\begin{tikzcd}
    \mathcal{H}^\lambda
        \arrow[r, "\operatorname{Hecke}"]
        \arrow[d, "p_1"]
    & W^+_{\E_\lambda}
        \arrow[d, "h"]
    \\
    W^+_{\E_0}
        \arrow[r, "h"]
    & \mathbb{A}_G.
\end{tikzcd}
\end{equation}
Here, $\mathcal{H}^\lambda$ is the subspace of $W^+_{\E_0}\times \Gr_G^\lambda$ consisting of pairs $(\E, A)$ such that $A$ defines a Hecke transformation of $\E$ at $c$ (in other words, such that $A$ lies in the appropriate affine Springer fibre), $p_1$ is the first projection, and $\E_\lambda$ denotes a certain $\C^\times$-fixed point obtained as a Hecke transformation of $\E_0$. The map labelled ``$\operatorname{Hecke}$'' takes a pair $(\E,A)$ to the Hecke transformation of $\E$ at $c$ defined by $A$, and one can show (at least under the minusculity assumption) that it is an isomorphism. Since the restriction of $h$ to $W^+_{\E_0} = s(\mathbb{A}_G)$ is also an isomorphism, the diagram shows that the Hitchin map on $W^+_{\E_\lambda}$ is isomorphic to $p_1\colon \mathcal{H}^\lambda\to W^+_{\E_0}$.

Evaluating at $c$ in (the left of) diagram \eqref{eq:Hecke_diagram} yields
\begin{equation*}
\label{eq:Hecke_diagram_at_c}
\begin{tikzcd}
    \mathcal{H}^\lambda
        \arrow[r, "\ev_c"]
        \arrow[d, "p_1"]
    & \Tilde{\s}^\lambda
        \arrow[d]
    \\
    W^+_{\E_0}
        \arrow[r, "h"]
    & \s,
\end{tikzcd}
\end{equation*}
where $\s$ is a Kostant section (used in constructing $s$) and $\Tilde{\s}^\lambda$ is the subspace of $\s\times \Gr^\lambda$ consisting of pairs $(x,A)$ such that $A$ is in the affine Springer fibre of $x$. Since $\lambda$ is minuscule, $\Tilde{\s}^\lambda$ is isomorphic to a Grothendieck--Springer partial resolution, restricted to $\s$. Now, as explained in \cite{HauselRychlewicz2025_fps} one can interpret $\Tilde{\s}^\lambda$ as a fixed point scheme, and conclude that its ring of global regular functions is isomorphic to $H^*_G(\Gr^\lambda) \coloneqq H^*_{G(\C)}(\Gr^\lambda(\C),\C)$, hence also to the medium algebra $\M^\lambda(\g)$. In fact, by \cite[Thm.\,1.2]{HauselRychlewicz2025_fps} (or directly), $\Tilde{\s}^\lambda$ is affine, so this says $\Tilde{\s}^\lambda \cong \Spec(\M^\lambda(\check{\g}))$. We may summarise the situation in the following commutative diagram:
\begin{equation}
\label{eq:big_HB_MA_diagram}
\begin{tikzcd}
    W^+_{\E_\lambda}
        \arrow[d, "h"]
    & \mathcal{H}^\lambda
        \arrow[l, "\operatorname{Hecke}"', "\cong"]
        \arrow[d]
        \arrow[r, "\ev_c"]
    & \Tilde{\s}^\lambda
        \arrow[r, "\rho", "\cong"']
        \arrow[d]
    & \Spec H^*_G(\Gr_G^\lambda)
        \arrow[r, "\cong"]
        \arrow[d]
    & \Spec \M^\lambda(\check{\g})
        \arrow[d]
    \\
    \mathbb{A}_G
    & W^+_{\E^+_0}
        \arrow[l, "h"', "\cong"]
        \arrow[r, "\ev_c"]
    & \s
        \arrow[r, "\cong"]
    & \Spec H^*_G(\pt)
        \arrow[r, "\cong"]
    & \check{\g}/\!/\check{G}.
\end{tikzcd}
\end{equation}

The above relationship between upward flows in Higgs bundle moduli spaces and spectra of medium algebras is expected to generalise to arbitrary dominant coweights, but how to implement this correctly appears to be quite subtle. Its importance lies in possible connections to mirror symmetry \cite{Hausel2023_ICM_mirror}, but we shall not go into this here. Instead, we use \eqref{eq:big_HB_MA_diagram} to relate our medium algebra automorphisms to automorphisms of $\mathbb{M}(G)$ studied by García-Prada--Ramanan \cite{GPRa2019_involutions}. The latter depend on several parameters, some of which involve the curve $C$ and therefore cannot be expected to be seen by medium algebras. The remaining parameter set is $\Out(G)\times \C^\times$, which exactly matches the parameters appearing for the medium algebra automorphisms in Definition \ref{defprop:autos}. To see if this is more than a coincidence, we follow an $\mathbb{M}(G)$-automorphism through the isomorphisms in \eqref{eq:big_HB_MA_diagram}.

Now, for $a\in \Out(G)$ and $\zeta\in \{1,-1\}$, García-Prada and Ramanan define an automorphism $\iota(a,\zeta)$ of $\mathbb{M}(G)$ by sending a Higgs bundle $(E,\phi)$ to $(\theta(E), \zeta\cdot \theta(\phi))$ for any lift $\theta\in \Aut(G)$ of $a$. Here $E$ is a principal $G$-bundle, $\phi$ a global section of $E(\g)\otimes K_C$, $\theta(E)$ denotes the principal $G$-bundle $E\times^G G$ where $G$ acts on $G$ through $\theta$, and similarly $\theta(\phi)$ is the section of $\theta(E)(\g)\otimes K_C$ obtained from $\phi$ using the derivative of $\theta$. 

We claim that any such $\iota(a,\zeta)$ fixes $\E_0$ as constructed in \cite[Ex.\,3.27]{González2025_vs_regular}. Indeed, $\E_0$ is a $\C^\times$-fixed point, so we may assume $\zeta = 1$. Moreover, the underlying principal $G$-bundle $E_0$ of $\E_0$ is obtained by extension of structure groups from a principal $H$-bundle, where $H$ is a one-dimensional torus integrating the line through a principal semisimple element $h$. Similarly, the Higgs field $\phi_0$ of $\E_0$ can be viewed as constant and equal to the nilpotent $e$ of a principal $\esl_2$-triple $(e,h,f)$. By choosing a (pinning-preserving) automorphism $\theta_a$ lifting $a$ such that $\theta_a$ fixes $e$ and $h$, one obtains $\theta_a(E_0) = E_0$ and $\theta_a(\phi_0) = \phi_0$. It is also clear that $\iota(a,\zeta)$ commutes with the $\C^\times$-action, so it preserves the upward flow $W^+_{\E_0}$.

Let us now see how $\iota(a,\zeta)$ interacts with Hecke transformations at $c$. In principle, the parametrisation of Hecke transformations of a Higgs bundle $(E,\phi)$ by affine Springer fibres depends on the choice of trivialisations of $E$. However, all Higgs bundles in the upward flow $W^+_{\E_0} = s(\mathbb{A}_G)$, have a common underlying principal $G$-bundle $E_0$, so we may choose trivialisations for all of them at once. As argued above, there is a (unique) preferred lift $\theta_a\in \Aut(G)$ of $a$ which has the property of fixing $E_0$, hence also its transition functions. With this setup, one easily sees that 
\begin{equation}
\label{eq:HB_and_MA_hecke_step}
    \iota(a, \zeta) \circ \operatorname{Hecke}
    =
    \operatorname{Hecke} \circ (\iota(a,\zeta) \times \theta_a),
\end{equation}
where $\theta_a$ on the right-hand side denotes the automorphism induced on $\Gr_G$. If $\theta_a(\Gr_G^\lambda) = \Gr^\lambda_G$ (which depends only on $a$), then $\iota(a,\zeta)\times \theta_a$ defines an automorphism of $\mathcal{H}^\lambda$ corresponding to the restriction of $\iota(a,\zeta)$ to $W^+_{\E_\lambda}$. 

Next, evaluation at $c\in C$ for $W^+_{\E_0} = s(\mathbb{A}_G)$ is to be understood in terms of the Higgs field. By construction, the Higgs fields on $s(\mathbb{A}_G)$ can be viewed as having values in $\s$, and in this sense there is an evaluation sending $(E_0,\phi)\in W^+_{\E_0}$ to $\phi(c)\in \s$. It also follows from the construction that this evaluation is $\C^\times$-equivariant for a standard $\C^\times$-action on $\s$ -- explicitly we have $\s = e + \C\cdot\{v_1,\ldots,v_r\}$ where $v_1,\ldots, v_r\in \g$ are lowest weight vectors for the principal semisimple element $h$, of weights $-m_1,\ldots, -m_r \in \Z_{<0}$, and $z\in\C^\times$ acts by $e+\sum_{i=1}^r a_iv_i \mapsto e + \sum_{i=1}^r z^{m_i+1}a_iv_i$. For the map $\ev_c$ in \eqref{eq:big_HB_MA_diagram} we then find
\begin{equation}
\label{eq:HB_and_MA_ev_step}
    \ev_c \circ (\iota(a,\zeta)\times \theta_a)
    =
    ((\zeta\circ\theta_a)\times \theta_a) \circ \ev_c
\end{equation}
where $\zeta$ on the right-hand side denotes the aforementioned $\C^\times$-action on $\s$.

The isomorphism $\rho\colon \Tilde{\s}^\lambda \to \Spec H^*_G(\Gr_G)$ in \eqref{eq:big_HB_MA_diagram} proceeds by identifying $\Tilde{\s}^\lambda$ with $\mathcal{Z}_G(\Gr_G^\lambda)$, the zero scheme of the $G$-space $\Gr^\lambda_G$ along $\s$. The resulting isomorphism $\mathcal{Z}_G(\Gr_G^\lambda) \cong \Spec H^*_G(\Gr^\lambda_G)$ is functorial for certain admissible morphisms of the group $G$ and the space $\Gr_G^\lambda$, and matches the grading of equivariant cohomology with twice that coming from the $\C^\times$-action on $\s$ above \cite[Thm.\,4.1]{HauselRychlewicz2025_fps}.\footnote{More precisely, the cohomological grading matches the one on $\C[\s]$ which puts the coordinate function for $v_i$ in degree $2(m_i+1)$.} For $\zeta = 1$, the functoriality readily shows that $\theta_a\times \theta_a \in \Aut(\Tilde{\s}^\lambda)$ corresponds to the automorphism $(\theta_a\inv)^*$ induced by $\theta_a\inv$ on equivariant cohomology. For $\zeta=-1$, the matching of gradings shows, via results of \S\ref{sec:MA}, that $((-1)\times 1)$ corresponds to the action $\sigma^*$ ($=(\sigma\inv)^*$) of a real structure inner to a split real structure of $G$. Combining these observations, we see 
\begin{equation}
\label{eq:HB_and_MA_coh_step}
    \rho \circ ((\zeta\circ \theta_a)\times \theta_a)
    =
    \begin{cases}
        (\theta_a\inv)^* \circ \rho & \text{if } \zeta = 1\\
        \sigma_a^* \circ \rho & \text{if } \zeta = -1 \text{ and }a^2 =1,
    \end{cases}
\end{equation}
where $\sigma_a$ denotes any real structure of $G$ whose action on the root datum matches that of $a$.

Finally, the isomorphism $\M^\lambda(\check{\g})\cong H^*_G(\Gr_G^\lambda)$ has been treated in Theorem \ref{thm:MA_autos_geometry}. Altogether, we conclude the following:
\begin{thm}
\label{thm:HB_MA_autos}
    Let $G$ be a connected semisimple complex algebraic group, $C$ a smooth projective curve of genus at least $2$, and $\mathbb{M}(G)$ the moduli space of polystable $G$-Higgs bundles on $C$. Let $(a,\zeta)\in \Out(G)\times \{1,-1\}$ with $a^2=1$ if $\zeta=-1$, and let $\lambda$ be a minuscule coweight of $G$ such that $\theta(\lambda) \in G\cdot \lambda$ for any (hence all) lifts $\theta\in \Aut(G)$ of $a$. Lastly, let $\Check{G}$ be the Langlands dual of $G$, with $\check{\g} = \Lie(\Check{G})$, and let $\check{a}\in \Out(\Check{G})\cong \Out(G)$ correspond to $a$.
    Then $\iota(a,\zeta)\in \Aut(\mathbb{M}(G))$ restricts to an automorphism of $W^+_{\E_\lambda}$, which induces $\iota_{\check{a}, \zeta} \in \Aut(\M^\lambda(\check{\g}))$ via the morphisms in \eqref{eq:big_HB_MA_diagram}.    
\end{thm}
\begin{rmk}\leavevmode
\begin{enumerate}
    \item The case $(a,\zeta) = (1,-1)$ was observed in \cite[Sec.\,4]{GH2024_even_very_stable}.
    \item For minuscule $\lambda$ not preserved by $a$, $\iota(a,\zeta)$ is still an isomorphism of $W^+_{\E_\lambda}$ with its image, corresponding in the same way to an isomorphism of medium algebras.
    \item The assumption that $\lambda$ be minuscule is only used to justify that all horizontal arrows in \eqref{eq:big_HB_MA_diagram} except $\ev_c$ are isomorphisms, though generalising this to arbitrary $\lambda$ does not seem entirely straightforward. For the correspondence of $\M(G)$-automorphisms to medium algebra automorphisms, the minuscule case covered by Theorem \ref{thm:HB_MA_autos} can be viewed as a proof of concept.
\end{enumerate}    
\end{rmk}

\section{Cohomology and Nadler groups}
\label{sec:coh_context}
Let $G$ be a quasi-anisotropic real reductive group. By the results of \S\ref{sec:coinv}, the (equivariant) cohomology of $\Gr_G(\R)$ is given by the coinvariants of the (equivariant) cohomology of $\Gr_G(\C)_\mathrm{adm}$. In proving this, we have already used a description of $H^*(\Gr_G(\R),\C)$, so this observation does not circumvent an independent calculation. Nevertheless, it is an interesting consequence that \emph{any} description of $H^*(\Gr_G(\C),\C)$ yields one of $H^*(\Gr_G(\R),\C)$ via taking coinvariants. In this section, we sketch how to interpret results of O'Brien \cite{obrien2024_splitting_rank} on the homology from this perspective. Curiously, this incorporates Nadler groups of without appealing to the real geometric Satake equivalence of \cite{Nadler2005_PSoRLG}, their main raison d'être. For simplicity, we discuss only the nonequivariant setting.

We begin with complex coefficients, in \S\ref{subsec:NG_coh_cplx}, matching our setting so far. There, we also address the case of $\SL_{2n+1,\R}$ not covered by \cite{obrien2024_splitting_rank}. Since O'Brien's results hold for integral coefficients up to inverting a finite set of primes, we sketch in \S\ref{subsec:coefficients} the extent to which our approach should generalise to such coefficients. 

\subsection{Nadler groups and (co)invariants}
\label{subsec:NG_coh_cplx}
\subsubsection{Construction of Nadler groups}
\label{subsubsec:NG}
We begin by recalling an explicit definition of the Nadler group of a connected real reductive group $G$, adapted from \cite[\S 10]{Nadler2005_PSoRLG}. We use the setup from \S\ref{subsubsec:real_red_setup}. In particular, we have a maximally split maximal torus $T$ with split part $S$, a minimal parabolic subgroup $P\le G$ with Levi factor $M = C_G(S)$ and positivity defined in terms of a Borel subgroup $B\le G_\C$ such that $T_\C\le B\le P_\C$. We write $2\rho_M \in X$ for the sum of all positive roots of $M_\C$; equivalently $2\rho_M$ is the sum of all positive imaginary roots in $\Phi$. We write $\Check{G}_\C$ for the Langlands dual of $\G_C$. Via Langlands duality, we may view $2\rho_M$ as a coweight of $\Check{G}_\C$, and define $\Check{L}_0 \coloneqq C_{\Check{G}_\C}(2\rho_M) \le \Check{G}_\C$ to be the corresponding Levi subgroup. This Levi subgroup inherits a choice of positive roots from $\Check{G}_\C$, namely $\Phi^\vee_{+,0} = \{\alpha \in \Phi^\vee_+\mid \langle 2\rho_M, \alpha\rangle = 0\}$, and we write $\Pi^\vee_0\subseteq \Phi^\vee_{0,+}$ for the corresponding set of simple roots for $\Check{L}_0$. 

The real structure $\sigma$ of $G_\C$ corresponding to $G$ acts on the based root datum of $G_\C$ as in \S\ref{subsubsec:real_red_setup}, and this action preserves $\Pi^\vee_0$. Now call $\alpha\in \Pi^\vee_0$ \emph{good} if $\alpha$ is nonadjacent to $\sigma(\alpha)$ in the Dynkin diagram of $\Check{L}_0$. We then define $\Check{L}_1\le \Check{L}_0$ as the subgroup generated by $\Check{T}$ and the root groups of $\Check{L}_0$ for the sets of roots $\{\alpha\mid \alpha \in \Pi^\vee_0 \text{ good}\}$ and $\{\alpha + \sigma(\alpha) \mid \alpha \in \Pi_0^\vee \text{ not good}\}$. The construction equips $\Check{L}_1$ with a choice of based root datum preserved by $\sigma$. It therefore gives rise to a unique inner-isomorphism class of pinning-preserving involutions of $\Check{L}_1$. Let $\eta \in \Aut(\Check{L}_1)$ be any representative of that class.
\begin{defn}
\label{def:NG}
    The \emph{Nadler group} of $G$ is $\Check{G}_\mathrm{Nad}\coloneqq (\Check{L}_1)^\eta$. Its Lie algebra is denoted $\check{\g}_\mathrm{Nad}$.
\end{defn}
\begin{rmk}
    One can check directly that $\Check{G}_\mathrm{Nad}$ does not depend on the choices of $T$, $B$ and $P$ (up to isomorphism), but this also follows from the main result of \cite{Nadler2005_PSoRLG}, where $\Check{G}_\mathrm{Nad}$ plays the role of the Langlands dual in an extension of the geometric Satake equivalence (hence is unique up to isomorphism via Tannakian formalism).
\end{rmk}
\begin{lem}
\label{lem:qa_NG}
    If $G$ is a quasi-anisotropic real form of a connected simply connected semisimple complex group, then:
    \begin{itemize}
        \item[(i)] $\Check{L}_1 = \Check{L}_0$, and
        \item[(ii)] $2\rho_M\in 2\Lambda$, so $\rho_M$ is a well-defined element of $\Lambda$.
    \end{itemize}
\end{lem}
\begin{proof}
    See Lemma \ref{lem:appx_A_qa_NG_lemma}.
\end{proof}

\subsubsection{Universal centralisers and Nadler groups}
\label{subsubsec:uc_and_ng}
The following is related to \cite[Sec.\,3.4]{CMNO2025_quaternionic}. Keeping the notation of \S\ref{subsubsec:NG}, we further assume that $G$ is quasi-anisotropic and (for simplicity) simply connected. The half sum of all positive roots $\rho$ is a principal cocharacter of $\Check{G}_\C$, so the adjoint action $\Ad_{\rho(-1)}$ is a principal involution as in Remark \ref{rmk:principal_invos} and its derivative $h\in \check{\g} \coloneqq \Lie((G_\C)^\vee)$ can be extended to a principal $\esl_2$-triple $(e,h,f)$, with $e$ a sum of simple root vectors of $\check{\g}$. 

Since $G$ is simply connected, Theorem \ref{thm:Gr_eq_coh} (or already \cite{Ginzburg2000_GeomSat}) implies that $H^*(\Gr_G(\C),\C)$ is isomorphic to $U(\check{\g}_e)$ as a graded algebra. Moreover, Theorem \ref{thm:MA_autos_geometry} matches the action of $\sigma$ on equivariant cohomology with the universal medium algebra automorphism $\Tilde{\iota}_{a,-1}$ for a class $a\in \Out(\Check{G})$ depending on $\sigma$, and Proposition \ref{prop:MA_slice_restriction} identifies the universal medium algebra with $\Sym(\jay^\vee)$. As in \S\ref{sec:Higgs}, we write $\theta_a$ for the unique pinning-preserving involution in the class $a$ which fixes $e$, $h$ and $f$. Specialising to zero fibres, we find (as in Proposition \ref{prop:MA_weight_action}) that the action of $\sigma$ on $H^*(\Gr_G(\C),\C)$ corresponds to the action of $\theta_a\circ \Ad_{\rho(-1)}$ on $U(\check{\g}_e)$. 

We now wish to connected this to the Nadler group. This turns out to work for those quasi-anisotropic $G$ whose Lie algebras do not contain a copy of $\esl_{2n+1}(\R)$, $n\in \Z_{>0}$; these groups are called \emph{splitting-rank} groups in \cite{obrien2024_splitting_rank}.\footnote{See \cite[Sec.\,VI.4.2]{LoosSySpII} for the symmetric space property behind the name.}

\begin{lem}
\label{lem:NG_lemma}
    Let $G$ be a quasi-anisotropic real form of a simply connected almost simple complex group, and assume that $\Lie(G)$ does not contain a copy of $\esl_{2n+1}(\R)$, $n\in \Z_{>0}$. Then:
    \begin{itemize}
        \item[(i)] $\theta_a\circ \Ad_{\rho(-1)}$ and $\Ad_{\rho_M(-1)}$ fix the same elements of $\check{\g}_e$.
        \item[(ii)] $\check{\g}_e$ contains a regular nilpotent element $e'$ of $\check{\g}_\mathrm{Nad}$ such that $\check{\g}_e\cap \check{\g}_\mathrm{Nad} = (\check{\g}_\mathrm{Nad})_{e'}$.
    \end{itemize}
\end{lem}
\begin{rmk}
    To be precise, this Lemma holds for a particular realisation of $\check{\g}_\mathrm{Nad}$, that is, a particular choice of the automorphism $\eta$ in Definition \ref{def:NG}.
\end{rmk}
\begin{proof}
    See Lemma \ref{lem:appx_A_sr_NG_lemma}.
\end{proof}

We recover the following consequence of \cite{obrien2024_splitting_rank}:
\begin{thm}
\label{thm:splitting_rank_coh}
    Let $G$ be a simply connected quasi-anisotropic semisimple real group whose Lie algebra does not contain a copy of $\esl_{2n+1}(\R)$, $n\in \Z_{>0}$. Then $H^*(\Gr_G(\R),\C) \cong U((\check{\g}_\mathrm{Nad})_{e'})$ where $e'$ is a regular nilpotent element of $\check{\g}_\mathrm{Nad}$ and the grading of $(\check{\g}_\mathrm{Nad})_{e'}$ is determined by the principal grading of $\check{\g}_e$ via the inclusion in Lemma \ref{lem:NG_lemma}(ii).
\end{thm}
\begin{proof}
    Indeed, the above discussion and Theorem \ref{thm:noneq_coinv} yield $H^*(\Gr_G(\R),\C) \cong U((\check{\g}_e)^{\theta_a\circ \Ad_\rho(-1)})$. If $G$ is the underlying real group of a complex reductive group $H$, then $\check{\g}_\mathrm{Nad} = \Lie(\Check{H})$ and the claim reduces to the complex case. Otherwise, Lemma \ref{lem:NG_lemma}(i) yields $(\check{\g}_e)^{\theta_a\circ \Ad_\rho(-1)} = (\check{\g}_e)^{\Ad_{\rho_M(-1)}} = \check{\g}_e\cap \check{\g}_\mathrm{Nad}$, which by part (ii) of the lemma agrees with $(\check{\g}_\mathrm{Nad})_{e'}$.
\end{proof}
\begin{expl}
\label{expl:A2n-1_coh}
    Consider $G = \SL_{n,\Ha}$, the quaternionic real form of $\SL_{2n,\C}$, for $n\ge 1$. The Langlands dual Lie algebra is $\esl_{2n}(\C)$, and a regular nilpotent element is given by the matrix $e = \sum_{i=1}^{2n-1}E_{i,i+1}$, where $E_{i,j}$ denotes the elementary matrix whose only nonzero entry is $1$ at position $(i,j)$. Its centraliser is 
    \begin{equation*}
        \check{g}_e = \C\cdot \{e^k \mid k = 1,\ldots 2n-1\},
    \end{equation*}
    and the principal grading gives $e^k$ degree $k$. One can then realise $\check{\g}_\mathrm{Nad}$ as the subalgebra consisting of trace-free $(2\times 2)$-block-matrices made of blocks of the form $\begin{psmallmatrix}
        z & 0 \\
        0 & z
    \end{psmallmatrix}$, $z\in \C$. This subalgebra is isomorphic to $\esl_n(\C)$, and it intersects $\check{\g}_e$ in $\C\cdot \{e^{2\ell}\mid \ell = 1,\ldots, n-1\}$, which is also the centraliser of the regular nilpotent element $e'\coloneqq e^2$ of $\check{\g}_\mathrm{Nad}$. Altogether, we have isomorphisms
    \begin{equation}
    \label{eq:A2n-1_expl}
        H^*(\Gr_{\SL_n}(\C), \C) \cong
        U(\esl_n(\C)_{e'}) \cong 
        H^*(\Gr_{\SL_{n,\Ha}},\C),
    \end{equation}
    the second of which doubles degrees (or the first if the middle term is graded as in Theorem \ref{thm:splitting_rank_coh}). This matches related degree-doubling phenomena observed in \cite{CMNO2025_quaternionic}.
\end{expl}

\subsubsection{The case of \texorpdfstring{$\SL_{2n+1,\R}$}{SL2n+1R}}
\label{subsubsec:A2n-1_case}
For the quasi-anisotropic group $G=\SL_{2n+1,\R}$, the Nadler group is $\PGL_{2n+1,\C}$, which is too large to make Theorem \ref{thm:splitting_rank_coh} work in this case. It is nevertheless interesting to point out that the conclusion of Example \ref{expl:A2n-1_coh} adapts to this case. Indeed, we have $H^*(\Gr_G(\C),\C) \cong \C[x_1,\ldots, x_{2n}]$ with $\deg(x_i) = 2i$, and the real structure acts by $(-1)^{\mathrm{deg}/2}$. Thus, Theorem \ref{thm:noneq_coinv} identifies the cohomology of $\Gr_G(\R)$ with $\C[x_2,x_4,\ldots, x_{2n}]$. We can in turn identify this with the cohomology ring of $\Gr_{\SL_{n+1}}(\C)$ with doubled degrees. Thus, we altogether have a degree-doubling isomorphism
\begin{equation}
\label{eq:A2n_coh_iso}
    H^*(\Gr_{\SL_{n+1}}(\C), \C) \cong H^*(\Gr_{\SL_{2n+1}}(\R),\C)
\end{equation}
analogous to \eqref{eq:A2n-1_expl}. 

Let us stress that \eqref{eq:A2n_coh_iso} is not a deep result on its own, and we currently do not have an entirely satisfying explanation for the appearance of $\SL_{n+1}$. The point is rather that it may be possible to fit $\SL_{2n+1,\R}$ into the setting of \cite{obrien2024_splitting_rank}. Perhaps this could be of interest in the study of a real derived Satake equivalence, analogous to the splitting-rank cases discussed in \cite{CMNO2025_quaternionic,CO2024_lorentzian_octonionic}.

\subsection{Coefficient rings}
\label{subsec:coefficients}
We have so far used cohomology with complex coefficients, whereas O'Brien \cite{obrien2024_splitting_rank} uses coefficients in $\Z[1/n]$ for some natural number $n$ (depending on the group). This suggests that our approach should also work more generally. Thus, let us briefly review the key steps underlying the coinvariant descriptions in \S\ref{sec:coinv} to see what happens for such coefficient rings. For this, $G$ will denote a connected quasi-anisotropic real reductive group.
\begin{enumerate}
    \item Additively, the integral cohomology of complex affine Grassmannians is torsion-free. The same is true for quasi-anisotropic real affine Grassmannians, except when the Lie algebra contains a copy of $\esl_{2n+1}(\R)$, $n\in \Z_{>0}$ (see Proposition \ref{prop:even_Iwahori_dimensions}). Even in those cases, it seems likely that there is only $2$-torsion. 
    \item The cohomology of $BG(\C)$ is a polynomial ring (in even degrees) whenever the cohomology of $G(\C)$ is an exterior algebra, and this happens whenever the cohomology of either of them is torsion-free. (These are classical results of Borel \cite{Borel1961_sg_et_torsion}, see also \cite[Thm.\,7.40]{MNN2017_nilpotence}). We should therefore invert the finitely many torsion primes of $G(\C)$ and $G(\R)$. Then, if the affine Grassmannians involved have only even cohomology and the classifying spaces have no $\Z[1/n]$-torsion, then we can reduce to nonequivariant cohomology as in \S\ref{sec:coinv}.
    \item Unsurprisingly, rational homotopy theory does not generalise to $\Z[1/n]$-coefficients. It was used to understand the map $H^*(\Gr_G(\C),\C)\to H^*(\Gr_G(\R),\C)$ induced by inclusion, with both domain and codomain being polynomial rings. In fact, it is well known that these are not polynomial rings integrally. Thus, write $R\coloneqq \Z[1/n]$ for (some) $n\in \Z_{>0}$, and consider
    \begin{equation*}
        i^*\colon H^*(\Gr_G(\C),R)\to H^*(\Gr_G(\R),R)
    \end{equation*}
    and suppose for simplicity that $G$ is simply connected. Proposition \ref{prop:inclusion_to_loop_Cartan} identifies $i^*$ with the map $(\Omega C)^*$ induced by a looped Cartan embedding $C\colon X\to K$; $K$ being maximal compact in $G(\C)$ and $X = K/(K\cap G(\R))$. 
    
    Suppose that $n$ is chosen so that $H^*(K, R)$ and $H^*(X, R)$ are exterior algebras and $C^*$ is a surjection on generators -- inspecting the proof of Proposition \ref{prop:nhz_Cartan} shows that this is the case when $n$ is divisible by the torsion primes of $K$ and $X$ as well as $2$. Although $H^*(\Omega K, R)$ and $H^*(\Omega X,R)$ are no longer polynomial rings, they are then still computable from the exterior algebras $H^*(K, R)$ and $H^*(X,R)$ via Eilenberg--Moore spectral sequences. More precisely, the spectral sequences collapse on the second page and yield (infinite) generating sets for the resulting associated graded algebras, corresponding to the exterior algebra generators of $H^*(K, R)$ and $H^*(X, R)$. Moreover, Proposition \ref{prop:nhz_Cartan} should adapt to showing a coinvariant description for these associated graded algebras, which should then also yield one for the actual cohomology algebras over $R$.
\end{enumerate} 
Working out the various details in the above list is left for future work.

\section{Medium algebras for symmetric pairs}
\label{sec:sym_pairs}
The Kirillov algebras discussed in Section \ref{subsec:MA_def} generalise to symmetric pairs as follows. Let $\g$ be a complex reductive Lie algebra with an involutive automorphism $\eta\in \Aut(\g)$. Write $\kay\coloneqq \g^\eta$ and $\p\coloneqq \kay^{-\eta}$, so that $\g = \kay \oplus \p$. Then $\kay$ is again a complex reductive Lie algebra, so we may consider a finite dimensional irreducible $\kay$-representation $V^\lambda$ of highest weight $\lambda$. We then define its \emph{symmetric pair Kirillov algebra} to be
\begin{equation*}
    \Kir^\lambda(\g,\eta)
    \coloneqq
    (\C[\p]\otimes \End(V^\lambda))^\kay,
\end{equation*}
where the $\kay$-action on $\p$ is by restriction of the adjoint action. This is a graded algebra over $\C[\p]^\kay$ with similar structural properties to Kirillov algebras. Moreover, if $\g = \g_0 \oplus \g_0$ for some complex reductive $\g_0$, with $\eta(x,y) = (y,x)$, then we get $\kay \cong \g_0 \cong \p$, and $\Kir^\lambda(\g,\eta) \cong \Kir^\lambda(\g_0)$, so this generalises the usual Kirillov algebras. 

The definition can be motivated by its relation to Harish-Chandra modules. Indeed, the symmetric pair $(\g,\kay)$ is an algebraic incarnation of the pair $(\g,K_0)$ where $K_0$ is the maximal compact subgroup of a real reductive group $G_0$ whose Lie algebra is a real form of $\g$. Such pairs and their bimodules arise in the representation theory of $G_0$. As observed in \cite{Higson2011_Analogy}, replacing $G_0$ by its \emph{Cartan motion group} $K_0\ltimes \p_0$ (where $\p_0 = \p\cap \Lie(G)$) takes certain Hecke algebras for isotypic components relevant to the representation theory of $G_0$ to symmetric pair Kirillov algebras. In this sense, the algebras $\Kir^\lambda(\g,\eta)$ are degenerations of algebras relevant to representation theory of real reductive groups. See also \cite{Hausel_Center_of_Kostant} for a related discussion.

We now define the \emph{symmetric pair medium algebra} $\M^\lambda(\g,\eta)\subseteq \Kir^\lambda(\g,\eta)$ to be the unital $\C[\p]^\kay$-subalgebra generated by $\kay$-equivariant morphisms $\p\to \kay$ (composed with the representation $\kay\to \End(V^\lambda)$). These algebras should be particularly well-behaved when $\eta$ is a principal involution (see Remark \ref{rmk:principal_invos}), in other words, the Cartan involution of a quasi-split real form.\footnote{Indeed, it is well known that we then have a surjection $\C[\g]^\g\to \C[\p]^\kay$, \eg using a Kostant--Rallis section \cite{KostantRallis_1971}.} In this case, we\footnote{The ideas sketched here were informed by discussions with Oscar García-Prada, Miguel González and Tamás Hausel.} expect a connection to the automorphisms of \ref{subsec:MA_autos} of the following form.  
\begin{conj}
\label{conj:principal_coinv}
    Let $\g$, $\eta$, $\kay$ and $\p$ be as above, and assume that $\eta$ is principal. Let $\lambda$ be a dominant integral weight of $\kay$. Then there exist a dominant integral weight $\lambda'$ of $\g$ such that the automorphism $\iota = \iota_{\eta,-1}$ of $\M^{\lambda'}(\g)$ is defined, and we have
    \begin{equation*}
        \M^\lambda(\g,\eta) \cong \M^\lambda(\g)_\iota,
    \end{equation*}
    possibly after taking the torsion-free quotient on the right-hand side, compatibly with an isomorphism $\C[\p]^\kay \cong \C[\g]^\g_\iota$.
\end{conj}

More details on and progress towards Conjecture \ref{conj:principal_coinv} will be given elsewhere. It offers a purely algebraic description of the coinvariants, which seems more closely related to the fixed points of Higgs bundle automorphisms given in \cite[Thm.\,6.3]{GPRa2019_involutions}. Moreover, comparing the conjecture with the geometric description in Theorem \ref{thm:coinv_main} suggests an interesting correspondence between quasi-split real structures and quasi-anisotropic real structures of the Langlands dual.

\clearpage
\appendix

\section{Case analyses for real structure actions on root systems}
\label{appx:case_proofs}
In this appendix, we prove several results about quasi-anisotropic real structures cited in the main text by case distinctions. All of them reduce easily to the four families listed in Table \ref{tab:simple_qa}, so we first establish some setup in each of them. When working with matrices, $E_{i,j}$ denotes the matrix with $(i,j)$-entry $1$ and all other entries $0$, and $I_n$ the $n\times n$ identity matrix. 

\subsection{Setup}
\label{sec:appx_setup}
In each case, we describe the action of a real structure $\sigma$ of a connected complex reductive group $G$ on the root system (with respect to a maximally split $\sigma$-stable maximal torus), as well as a system of positive roots corresponding to a Borel subgroup $B$ chosen as in \S\ref{subsubsec:real_red_setup}. This requirement on $B$ translates to $\sigma$ taking nonimaginary positive roots to positive roots, so it can be checked directly on the positive roots. Since (for geometric reasons) we work with the adjoint forms, some statements about roots are checked in the equivalent setting of coroots. On the way, we compute the real minuscule coweights, which correspond to the admissible components of $\Gr_G(\C)$ as listed in Table \ref{tab:admissible_minuscules}, the sum $2\check{\rho}_M$ of positive imaginary coroots, and the class $a\in\Out(G_\C)$ corresponding to $\sigma$.

\subsubsection{\texorpdfstring{$\PGL_{2n+1}$}{PGL2n+1}}
\label{subsubsec:case1}
Here $G=\PGL_{2n+1}$, $n\in \Z_{>0}$, with real structure $\sigma(A) = \overline{A}$. This is a split real structure; the standard diagonal maximal torus and upper triangular Borel subgroup fulfil the requirements, and $\sigma$ acts trivially on the root datum. This also shows that all connected components of $\Gr_G(\C)$ are admissible. There are no imaginary coroots, so $2\check{\rho}_M = 0$, and the class $a=1$.

\subsubsection{\texorpdfstring{$\PGL_{2n}$}{PGL2n}}
\label{subsubsec:case2}
Here $G=\PGL_{2n}$, $n\in \Z_{>0}$, with $\sigma(A) = J\overline{A}J\inv$, where $J$ is a $2\times 2$-block-diagonal matrix all of whose diagonal blocks are equal to $\begin{psmallmatrix}
    0 & -1\\
    1 & 0
\end{psmallmatrix}$. We use the standard diagonal maximal torus $T$ and standard upper triangular Borel subgroup $B$. We have $\Lambda^\vee \cong \Z^{2n}/\Z\cdot \mathbf{1}$ where $\mathbf{1} = (1,1,\ldots, 1)$, and a coweight $\lambda = [(\lambda_1, \ldots, \lambda_{2n})]$ is dominant if and only if $\lambda_1\ge \lambda_2 \ge \cdots \ge \lambda_{2n}$. The coroots are of the form $[e_i - e_j]$, $i\ne j$, where $e_i$ denotes the tuple with $1$ in position $i$ and all other entries zero. Positivity corresponds to $i<j$ and the simple coroots are $\check{\alpha}_i = [e_i - e_{i+1}]$, $i=1,\ldots, 2n-1$.  The action of $\sigma$ on $\Lambda^\vee$ is
\begin{equation*}
    \lambda =
    [(\lambda_1,\ldots, \lambda_{2n})]
    \mapsto
    \sigma(\lambda) =
    [(\lambda_2, \lambda_1, \lambda_4, \lambda_3, \ldots, \lambda_{2n}, \lambda_{2n-1}],
\end{equation*}
which implies that
\begin{itemize}
    \item $\lambda\in \Lambda^\vee_\R$ if and only if $\lambda_{2i-1} = \lambda_{2i}$ for all $i=1,\ldots, n$;
    \item a coroot $[e_i-e_j]$ is imaginary if $j=i+1$ and $i$ is odd, and complex otherwise;
    \item nonimaginary positive coroots are sent to positive coroots.
\end{itemize}
It follows from Proposition \ref{prop:qa_max_tori} that the chosen maximal torus is maximally split, and the last item in the above list shows that $B$ fulfils the requirements. The nonzero minuscule coweights are $\check{\varpi}_1,\ldots\check{\varpi}_{2n-1}$, where $\check{\varpi}_i = [e_1+\cdots + e_i]$. Thus, the nonzero real minuscule coweights are $\check{\varpi}_{2i}$, $i=1,\ldots, n-1$. Moreover, the second item shows that $2\check{\rho}_M$ is represented by $\sum_{i=1}^{n}e_{2i-1}-e_{2i}$. Lastly, $\sigma$ acts on coweights by a permutation, hence matches the action of a Weyl group element, so $a=1$.
    
\subsubsection{\texorpdfstring{$\PSO_{2n}$}{PSO2n}}
\label{subsubsec:case3}
Here $G=\PSO_{2n}$, $n\ge 4$, which we realise as the complex projective orthogonal group preserving the bilinear form with matrix $\begin{psmallmatrix}
        0 & I_n\\
        I_n & 0 
\end{psmallmatrix}$. As the real structure, we may take $\sigma(A) = R\overline{A}R\inv$ for the block matrix $R = \begin{psmallmatrix}
        E_{1,1} & I_n - E_{1,1}\\
        I_n - E_{1,1} & E_{1,1}
\end{psmallmatrix}$.
In other words, $R$ is the permutation matrix swapping $i$ and $n+i$ for $i=2,\ldots, n$. We use the standard diagonal maximal torus and the standard upper triangular Borel subgroup $B$. The coweight lattice is then isomorphic to $\Z^n \oplus (\mathfrak{1}{2}+\Z)^n$, in such a way that $(\lambda_1,\ldots, \lambda_n)$ is dominant if and only if $\lambda_1\ge \cdots \lambda_n$ and $\lambda_{n-1} \ge |\lambda_n|$. The positive coroots are $e_i\pm e_j$, $i\ne j$, and the simple coroots are $e_i-e_{i+1}$,$ i=1\,\ldots, n-1$ together with $e_{n-1} + e_n$. 

In this setting, $\sigma$ acts as
\begin{equation*}
    (\lambda_1, \lambda_2,\ldots, \lambda_n)
    \mapsto
    (\lambda_1, - \lambda_2, \ldots, -\lambda_n).
\end{equation*}
The real coweights are those of the form $(\mu_1,0,\ldots, 0)$, $\mu_1\in \Z$. We see immediately that the coroots $e_1\pm e_j$, $j>1$ are complex, with all other coroots being imaginary, so $B$ fulfils the requirements. Moreover, the only nonzero real minuscule coweight is $\check{\varpi}_1 = (1,0,\ldots, 0)$. We also find $2\check{\rho}_M = \sum_{i=2}^{n-1}2(n-i)e_i$. Lastly, the Weyl group of $G$ is generated by permutations and changes of an even number of signs. Thus, $a=1$ if $n$ is odd. If $n$ is even, then $a$ corresponds to the diagram automorphism exchanging the last two simple roots.

\subsubsection{\texorpdfstring{$E_{6,\mathrm{ad}}$}{E6ad}}
\label{subsubsec:case4}
\begin{figure}[ht]
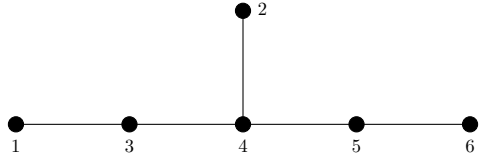

    \centering
        \dynkin[edge length = 1.5cm, rootradius = 0.1cm, labels = {1,2,3,4,5,6}] E6
    \caption{Indexed Dynkin diagram of $E_6$}
    \label{fig:E6_indexing}
\end{figure}

Here we work more abstractly, starting with the Cartan involution. We use the established labelling of the Dynkin diagram displayed in Figure \ref{fig:E6_indexing}, in which the nontrivial diagram automorphism is the one fixing nodes $2$ and $4$ and swapping $1\leftrightarrow 6$, $2\leftrightarrow 5$.
This diagram automorphism lifts to a pinning-preserving involution $\theta$ of $E_{6,\mathrm{ad}}$, and $\sigma$ is a real structure with Cartan involution $\theta$. We begin with the (real form of the) maximal torus $T$ and Borel subgroup $B$ obtained from the preserved pinning. By Proposition \ref{prop:qa_max_tori}, the maximal torus is automatically maximally split, but $B$ need not fulfil the requirements. We will therefore have to change the initial system of simple roots $\{\alpha_1,\ldots, \alpha_6\}$.

The action of $\sigma$ on the roots is minus that of $\theta$, so it can be inferred from the diagram automorphism. We claim that the following roots form a simple system which fulfils the requirements\footnote{It was found by defining positivity with respect to a suitable ordered basis of the Cartan subalgebra, followed by direct computation of the resulting simple roots.}: 
\begin{align*}
    &\alpha^\mathrm{new}_1 = -\alpha_2 - \alpha_3 - \alpha_4, \\
    &\alpha^\mathrm{new}_2 = \alpha_1 + \alpha_2 + 2\alpha_3 + 2\alpha_4 + 2\alpha_5 + \alpha_6,\\
    &\alpha^\mathrm{new}_3 = \alpha_2 + \alpha_3 + 2\alpha_4 + \alpha_5,\\
    &\alpha^\mathrm{new}_4 = - \alpha_1 - \alpha_2 - 2\alpha_3 - 3\alpha_4  -2\alpha_5 - \alpha_6,\\
    &\alpha^\mathrm{new}_5 = \alpha_1 + \alpha_2 +\alpha_3 + 2\alpha_4 + \alpha_5 + \alpha_6,\\
    &\alpha^\mathrm{new}_6 = -\alpha_2 - \alpha_4  - \alpha_5 - \alpha_6.
\end{align*}
Indeed, one may check that they are linearly independent and their pairwise inner products match the Dynkin diagram as labelled in Figure \ref{fig:E6_indexing}. Moreover, $\sigma$ acts on these roots as follows: 
\begin{align*}
    &\sigma(\alpha_1^\mathrm{new}) = \alpha^\mathrm{new}_1 +\alpha^\mathrm{new}_2 + 2\alpha^\mathrm{new}_3 + 2\alpha^\mathrm{new}_4 +\alpha^\mathrm{new}_5,\\
    &\sigma(\alpha^\mathrm{new}_2) = -\alpha^\mathrm{new}_2,\\
    &\sigma(\alpha^\mathrm{new}_3) = -\alpha^\mathrm{new}_3,\\
    &\sigma(\alpha^\mathrm{new}_4) = -\alpha^\mathrm{new}_4,\\
    &\sigma(\alpha^\mathrm{new}_5) = -\alpha^\mathrm{new}_5,\\
    &\sigma(\alpha^\mathrm{new}_6) = \alpha^\mathrm{new}_2 + \alpha^\mathrm{new}_3 +2\alpha^\mathrm{new}_4 +2\alpha^\mathrm{new}_5 +\alpha^\mathrm{new}_6.
\end{align*}
A general positive root $\alpha$ for this new simple system is of the form $\alpha = \sum_{i=1}^6n_i\alpha^\mathrm{new}_i$ for $n_1,\ldots, n_6\in \N$. If $n_1=n_6=0$, then $\alpha$ is imaginary. Otherwise, $\sigma(\alpha)$ is also positive because $\sigma$ preserves the coefficients of $\alpha^\mathrm{new}_1$ and $\alpha^\mathrm{new}_6$, so the new simple system (or its corresponding Borel subgroup) fulfils the requirements.

To ease the notation, we now drop the superscript $\mathrm{new}$ and redefine $\alpha_i\coloneqq \alpha_i^\mathrm{new}$ -- this effectively only changes the action of $\sigma$ on simple (co)roots to be the one computed above. By expressing the fundamental coweights $\check{\varpi}_1,\ldots \check{\varpi}_6$ in terms of simple coroots (according to the Cartan matrix), one finds that $\Lambda^\vee_\R$ is spanned by $\check{\varpi}_1$ and $\check{\varpi}_6$. Since these are precisely the nonzero minuscule coweights, we see that all components of $\Gr_{E_{6,\mathrm{ad}}}$ are admissible. Finally, the coefficient argument above also shows that the positive imaginary coroots are precisely those spanned by the subsystem $\{\check{\alpha}_2,\ldots, \check{\alpha}_4\}$ of type $D_4$. It follows that $2\check{\rho}_M = 6\check{\alpha}_2 + 6\check{\alpha}_3 + 10\check{\alpha}_4 + 6\check{\alpha}_5$. For later use, we also record that $\check{rho} = 8\check{\alpha}_1 +  11\check{\alpha}_2 + 15\check{\alpha}_3 + 21\check{\alpha}_4 + 15\check{\alpha}_5 + 8\check{\alpha}_6.$

The class $a=1$ because $\sigma$ is inner to a split real structure. Indeed, it is not inner to an anisotropic real structure (by construction) and neither is the split real structure. To see this, note that a split real structure has a Cartan involution acting on the roots by $-1$. If this were an inner automorphism, then the Weyl group would contain $-1$, so all its exponents would have to be odd, which is not the case for $E_6$. Since there are exactly two inner classes, $\sigma$ is inner to a split real structure as claimed.

\subsection{Proofs}
We first prove two lemmas from \S\ref{sec:coh_context}, and then a proposition from \S\ref{subsec:asV_fixed_loci}.
\begin{lem}[= Lemma \ref{lem:qa_NG}]
\label{lem:appx_A_qa_NG_lemma}
    If $G$ is a quasi-anisotropic real form of a connected simply connected semisimple complex group. Then, with notation as in \S\ref{subsubsec:NG}:
    \begin{itemize}
        \item[(i)] $\Check{L}_1 = \Check{L}_0$, and
        \item[(ii)] $2\rho_M\in 2\Lambda$, so $\rho_M$ is a well-defined element of $\Lambda$.
    \end{itemize}
\end{lem}
\begin{proof}
    Since $G$ is assumed simply connected, $\Lambda$ is the dual lattice of the system of coroots. Both statements then only concern the root systems, so we may work with the Lie algebra $\g_\C$ of $G_\C$ and assume it to be either simple or the direct sum of two simple Lie algebras interchanged by the real structure $\sigma$ corresponding to $G$. In the latter case, there are no imaginary roots, so that $\Check{L}_0 = \Check{G}$ and the first claim is easy to see. The second claim is then also clear because $\rho_M = 0$ in this case. If $\g_\C$ is simple, then either $\sigma$ is an anisotropic real structure or $\g_\C$ is the Lie algebra of one of the four groups considered in Table \ref{tab:simple_qa}. In the anisotropic case, all roots are imaginary, so $\Check{L}_1 = \Check{L}_0 = \Check{T}$ and $2\rho_M$ is the sum of positive roots, which is well known to lie in $2\Lambda$.

    We are left with the four cases considered in the setup part of this appendix. In each of them, we have to check that $2\rho_M$ has even pairing with all (simple) coroots and that the diagram involution of its centraliser induced by $\sigma$ does not interchange two adjacent nodes. Note that those cases were discussed in adjoint type, so we should read all mentions of coroots there to be about roots. 
    
    In case \ref{subsubsec:case1}, $2\rho_M = 0$ and the diagram automorphism is trivial. 
    
    In case \ref{subsubsec:case2}, $2\rho_M$ is represented by $\sum_{i=1}^{n-1}e_{2i-1}-e_{2i}$, and the parity claim is easy to see. The derived subalgebra of its centraliser is of type $A_{n-1}\oplus A_{n-1}$ and the diagram involution swaps the two summands, so it does not interchange adjacent nodes.

    In case \ref{subsubsec:case3}, we have $2\rho_M  = \sum_{i=2}^{n-1}2(n-i)e_i$ which is visibly even. The roots of $\Check{L}_0$ are the coroots pairing trivially with $2\rho_M$, namely those with roots $\pm e_1 \pm e_n$. The positive roots $e_1+e_n$ and $e_1-e_n$ are interchanged by $\sigma$ but they are orthogonal. 
    
    Lastly, in case \ref{subsubsec:case4}, the formula for $2\rho_M$ again makes it clearly even. Writing a general positive root as $\alpha = \sum_{i=1}^6 n_i\alpha_i$, the pairing $(\rho_M,\alpha)$ is seen to equal $-3(n_1+n_6) + n_2 + n_3 + n_4 + n_5$. One then checks (by exhaustion) that the roots orthogonal to $\rho_M$ are those with $(n_1,\ldots, n_6)$ in the list $(1,1,1,1,0,0)$, $(1,0,1,1,1,0)$, $(0,1,0,1,1,1)$, $(0,0,1,1,1,1)$, $(1,1,2,2,1,1)$ and $(1,1,1,2,2,1)$. The action of $\sigma$ on these matches the pairs 
    \begin{equation*}
        (1,1,1,1,0,0)\leftrightarrow (1,0,1,1,1,0), 
        \,
        (0,1,0,1,1,1) \leftrightarrow (0,0,1,1,1,1),
        \,
        (1,1,2,2,1,1) \leftrightarrow (1,1,1,2,2,1),
    \end{equation*}
    and one checks that they are all orthogonal.
\end{proof}

\begin{lem}[= Lemma \ref{lem:NG_lemma}]
\label{lem:appx_A_sr_NG_lemma}
    Let $G$ be a quasi-anisotropic real form of a simply connected almost simple complex group, and assume that $\Lie(G)$ does not contain a copy of $\esl_{2n+1}(\R)$, $n\in \Z_{>0}$. By Lemma \ref{lem:appx_A_qa_NG_lemma}, we have $\check{\g}_\mathrm{Nad} = \Lie(\Check{L}_0)^\eta$ for an involution $\eta$, which can be chosen so that the following hold:
    \begin{itemize}
        \item[(i)] $\theta_a\circ \Ad_{\rho(-1)}$ and $\Ad_{\rho_M(-1)}$ fix the same elements of $\check{\g}_e$.
        \item[(ii)] $\check{\g}_e$ contains a regular nilpotent element $e'$ of $\check{\g}_\mathrm{Nad}$ such that $\check{\g}_e\cap \check{\g}_\mathrm{Nad} = (\check{\g}_\mathrm{Nad})_{e'}$.
    \end{itemize}
\end{lem}
\begin{proof}
    We can work on the level of Lie algebras and assume $\g_\C$ to be simple, with $\g$ not equal to $\esl_{2n+1}(\R)$, $n\in \Z_{>0}$. If $\g$ is anisotropic, then $\check{\g}_\mathrm{Nad} = 0$ and the claims hold trivially. Otherwise, $\g_\C$ and the real structure $\sigma$ are as in one of the last three cases discussed in the setup.
    
    \textbf{Case \ref{subsubsec:case2}:} Here $\theta_a = \id$. Moreover, $\rho-\rho_M = \sum_{i=1}^{n}(n-i)(e_{2i-1} + e_{2i})$, is an even weight, so that the adjoint actions of $\rho_M(-1)$ and $\rho(-1)$ agree. The rest of this case is covered by Example \ref{expl:A2n-1_coh}.

    \textbf{Case \ref{subsubsec:case3}:} Here $\check{\g}_\mathrm{Nad}$ is simple of rank one, and its only positive root space is spanned by a sum $e'$ of root vectors for the roots $e_1+e_n$ and $e_1-e_n$ of height $n-1$. This is then automatically a regular nilpotent of $\check{\g}_\mathrm{Nad}$. We moreover have complete freedom in choosing the two root vectors summing to $e'$ (among nonzero elements of the given root spaces). 
    
    The principal grading (induced by the cocharacter $\rho$) decomposes $\check{\g}_e$ into graded pieces $(\check{\g}_e)^d$, $d > 0$, and $(\check{\g}_e)^d$ is contained in the span of positive root spaces of $\check{\g}$ for roots of height $d$. The degrees $d$ that appear are one less than the degrees of the invariant ring $\C[\so_{2n}(\C)]^{\SO_{2n}}$ (\cf \cite{Kostant1959_TDS}), so they are, with multiplicities, the odd integers $1,3,\ldots, 2n-3$, as well as $n-1$. Moreover, we can arrange for $e'$ to be contained in $(\check{\g}_e)^{n-1}$. We claim that the fixed subspaces of $\theta_a\circ\Ad_{\rho(-1)}$ and $\Ad_{\rho_M(-1)}$ on $\check{\g}_e$ are both spanned by $e'$.
    
    We have $\rho-\rho_M = \sum_{i=1}^{n-1}(n-i)e_i - \sum_{i=2}^{n-1}(n-i)e_i = (n-1)e_1$. If $n$ is odd, then this shows that $\rho(-1)$ and $\rho_M(-1)$ have the same adjoint actions. Moreover, for odd $n$ we also have $\theta_a = \id$. This immediately yields part (i) in this case, and for part (ii) it suffices to recall that $e'$ spans the only positive root space of $\check{\g}_\mathrm{Nad}$ and is contained in $\check{\g}_e$.
    
    If $n$ is even, then $\theta_a$ comes from the diagram involution exchanging the last two simple roots. We then realise $\so_{2n}(\C)$ as in the setup \ref{subsubsec:case3} and find that -- for the standard choice of $e = \begin{psmallmatrix}
        X & Y\\
        0 & - X^t
    \end{psmallmatrix}$, where $X=\sum_{i=1}^{n-1}E_{i,i+1}$ and $Y = E_{n-1,n} - E_{n,n-1}$, $\theta_a$ acts on it by $A\mapsto -CA^tC\inv$ where $C=\begin{psmallmatrix}
        I_n - E_{n,n} & E_{n,n}\\
        E_{n,n} & I_n - E_{n,n}
    \end{psmallmatrix}$. Either directly or from Lemma \ref{lem:qa_real_and_nci_roots}, we can see that $\theta_a$ acts trivially on all root spaces it preserves. For $d\ne n-1$, we have $\dim (\check{\g}_e)^d\le 1$, and if this equals one then a root computation shows that a nonzero element has nonzero component in a root space preserved by $\theta_a$. The two-dimensional space $(\check{\g}_e)^{n-1}$ is spanned by $e' = \begin{psmallmatrix}
        E_{1,n} & E_{n,1} - E_{1,n}\\
        0 & - E_{n,1}
    \end{psmallmatrix}$ and one other element $e''$ with nonzero component in a root space whose root is preserved by $\theta_a$. Thus, $\theta_a(e'') = e''$ and a direct computation shows $\theta_a(e') = -e'$. We conclude that $\check{\g}_e$ is spanned by eigenvectors of $\theta_a$, all but one of which have eigenvalue $1$, and only $e'$ has eigenvalue $-1$. Since all principal degrees are odd, the fixed subspace of $\theta_a\circ \Ad_{\rho(-1)}$ is spanned by $e'$. Using the identity $\rho - \rho_M = (n-1)e_1$, we similarly find that the fixed subspace of $\Ad_{\rho_M(-1)}$ is spanned by $e'$, which finishes this case.

    \textbf{Case \ref{subsubsec:case4}:} Here we again have $a=1$ and $\rho-\rho_M$ is the even weight $8(\alpha_1+\alpha_2 + \alpha_6) + 12(\alpha_3 + \alpha_5) + 16\alpha_4$, so $\theta_a\circ \Ad_{\rho(-1)}$ and $\Ad_{\rho_M(-1)}$ agree entirely. We have $\check{\g}_\mathrm{Nad} \cong \esl_3(\C)$, and upon choosing linearly independent simple root vectors $v_1,v_2$ of $\check{\g}_\mathrm{Nad}$, their sum $e'$ will be a regular nilpotent and span $(\check{\g}_{\mathrm{Nad}})_{e'}$ together with $[v_1,v_2]$. It thus remains to choose $e,$ $v_1$ and $v_2$ such that $e$ centralises $v_1+v_2$ and $[v_1,v_2]$.

    We begin by arbitrarily choosing nonzero simple root vectors $x_1,\ldots, x_6$ for $\check{\g} \cong \mathrm{e}_6$ and putting $e\coloneqq x_1 + \cdots + x_6$. The centraliser of $e$ then intersects the span of height $d$ root spaces whenever there are fewer height $d+1$ roots than height $d$ roots, and one finds that these values of $d$ are $1$, $4$, $5$, $7$, $8$ and $11$. In terms of coefficient vectors with respect to simple roots, the five height four roots are $\beta_1 = (1,1,1,1,0,0)$, $\beta_2 = (1,0,1,1,1,0)$, $\beta_3 = (0,1,0,1,1,1)$, $\beta_4 = (0,0,1,1,1,1)$ and $\beta_5 = (0,1,1,1,1,0)$. By \cite[p.\,197]{Gradl1996_exponents}, $(\check{\g}_e)^4$ has nonzero components in only the root spaces for $\beta_1,\ldots, \beta_4$, and we have seen that these are precisely the simple roots of the Levi subalgebra $\Lie(\Check{L}_0)$ relevant to $\check{\g}_\mathrm{Nad}$.\footnote{Note that the cited paper uses a different labelling of $E_6$.} We may therefore choose root vectors $y_i$ of $\beta_i$, $i=1,\ldots, 4$, so that $e$ centralises $y_1+\cdots + y_4$, and obtain simple root vectors $v_1 = y_1 + y_2$ and $v_2 = y_3+y_4$ of $\check{\g}_\mathrm{Nad}$. This entails $[e,v_1+v_2]=0$, and it remains to check $[e,[v_1,v_2]]= 0$. 

    To see this, consider the pinning-preserving involution $\xi$ which swaps $x_1$ with $x_6$ and $x_2$ with $x_5$. It then fixes $e$, so it must also fix $v_1 + v_2$. On the other hand, it interchanges $\beta_1$ with $\beta_3$ and $\beta_2$ with $\beta_4$. Therefore, it must map $y_1$ to $y_3$, $y_2$ to $y_4$ and $v_1$ to $v_2$. It follows that 
    \begin{equation*}
        \xi([e,[v_1,v_2]]) = 
        [e,[\xi(v_1),\xi(v_2)]] =
        [e, [v_2, v_1]] = 
        - [e,[v_1,v_2]].
    \end{equation*}
    Now, $[e,[v_1,v_2]]$ is a (possibly zero) root vector for the single root $\gamma$ of height $9$. We certainly have $\xi(\gamma) = \gamma$, so $\xi$ restricts to the root space of $\gamma$ as either $1$ or $-1$. If it was $-1$, then $\gamma$ would be a noncompact imaginary root for a quasi-anisotropic real structure with Cartan involution $\xi$, which is impossible by Lemma \ref{lem:qa_real_and_nci_roots}. It follows that $[e,[v_1,v_2]] = 0$ as required. This finishes the last case, hence the proof.
\end{proof}

\begin{prop}[= Proposition \ref{prop:unique_max_fixed_cwt}]
\label{prop:appx_A_main}
    Let $\sigma$ be a quasi-anisotropic real structure of $G$, and $\lambda\in \Lambda^\vee_+$ with $\sigma(\lambda) \in W\cdot \lambda$. Write $(\Lambda^\vee_+)_{\le \lambda}$ for the set of dominant coweights less than or equal to $\lambda$. Then the partially ordered set $(\Lambda^\vee_+)_{\le \lambda}\cap \Lambda^\vee_\R$ is either empty or contains a unique maximal element $\mu_\mathrm{max}$. It follows that $\Gr_G^{\le\lambda}(\C)^\sigma$ is either empty or equal to $\Gr_{G^\sigma}^{\le \mu_\mathrm{max}}(\R)$. 
\end{prop}
\begin{proof}
    The coweight lattice of $G$ projects to the coweight lattice of $G_\mathrm{ad}$, and this projection restricts to isomorphisms on given cosets (or $\sigma$-orbits of cosets) of the coroot lattice. Replacing $G$ by $G_\mathrm{ad}$ if necessary, we may therefore assume that $G$ is of adjoint type. Then $G$ splits into simple factors, and we may further restrict to a single $\sigma$-orbit of simple factors. 
    
    If there are two factors, then $G^\sigma$ is the underlying real reductive group of a complex reductive group $H$. We may then arrange for $\Lambda^\vee = \Lambda^\vee_H\oplus \Lambda^\vee_H$ where $\Lambda^\vee_H$ is the coweight lattice of $H$, with $\sigma$ acting on $\Lambda^\vee_H\oplus \Lambda^\vee_H$ as the swap involution. Then the assumption implies $\lambda = (\mu,\mu)$ for some dominant weight $\mu$ of $\Lambda_H$, and so $\lambda$ itself is the sought-after unique maximal element $\mu_\mathrm{max}$. If $G$ is simple and $\sigma$ is an anisotropic real structure, then $\Lambda^\vee_\R = 0$ and the claim is trivially true. It remains to prove the proposition under the assumption that $G$ is simple of adjoint type and $\sigma$ a non-anisotropic quasi-anisotropic real structure. These are the four cases set up in \ref{sec:appx_setup}, and we treat each of them separately.

    \textbf{Case \ref{subsubsec:case1}:} This is a split real structure, acting trivially on the coweight lattice, so the claim holds trivially.

    \textbf{Case \ref{subsubsec:case2}:}  
    Let $\lambda = [(\lambda_1, \ldots, \lambda_{2n}]$ be a dominant coweight. We first observe the following admissibility criterion:
    \begin{equation}
    \label{eq:appx_case_2_criterion}
        \text{If}
        \quad
        (\Lambda^\vee_+)_{\le \lambda}\cap \Lambda^\vee_\R\ne \emptyset
        \quad
        \text{then}
        \quad
        \sum_{i=1}^{2n} \lambda_i
        \text{ is even.}
    \end{equation}
    Indeed, subtracting a coroot from $\lambda$ leaves the parity of the sum in question invariant, but the sum is certainly even for a $\sigma$-fixed coweight. Therefore, we now assume that $\sum_{i=1}^{2n}\lambda_i$ is even, and construct the desired maximal $\sigma$-fixed coweight in two steps, illustrated by the example $\lambda = [(12,9,8,7,7,4,4,2,1,0)]$ for $n=5$:
    \begin{enumerate}
        \item Define $\lambda'\coloneqq \lambda - \sum_{i=1}^n\lfloor (\lambda_{2i-1}-\lambda_{2i})/2\rfloor \check{\alpha}_{2i-1} \in \Lambda^\vee_+$. Observe that any $\sigma$-fixed coweight $\mu\le \lambda$ must also fulfil $\mu\le \lambda'$ because the only way of lowering of the difference $\lambda_{2i-1}-\lambda_{2i}$ by subtraction of simple coroots is via $\check{\alpha}_i$. By construction, we have $\lambda'_{2i-1}-\lambda'_{2i}\in \{0,1\}$ for every $i\in \{1,\ldots, n\}$, and the difference 1 occurs an even number of times. In the example, we get $\lambda' = [(11,10,8,7,6,5,3,3,1,0)]$.
        \item Let $i_1, \ldots, i_{2m}$ denote the odd indices $i$, in ascending order, for which $\lambda'_i-\lambda'_i = 1$. We define 
        \begin{equation*}
            \mu_\mathrm{max}\coloneqq \lambda' - \sum_{j=1}^m\sum_{i=i_{2j-1}}^{i_{2j}}\check{\alpha}_i  
        \end{equation*}
        and claim that this has the desired properties. In the example, we have $i_1=1$, $i_2=3$, $i_3=5$ and $i_4=9$ and subtract $[e_1-e_4 + e_5 - e_{10}]$, obtaining $\mu_\mathrm{max} = [(10,10,8,8,5,5,3,3,1,1)]$. 
        
        In general, it is clear that $\mu_\mathrm{max}\le \lambda$, and it is easy to check that $\mu_\mathrm{max}\in \Lambda^\vee_{\R,+}$. It remains to show that, if $\mu\le \lambda'$ is any $\sigma$-fixed dominant coweight, then $\mu\le \mu_\mathrm{max}$. 
        The coweights below $\lambda'$ in the partial order $\le$ are of the general form $\lambda' - \sum_{i=1}^{2n-1}c_i\check{\alpha}_i$, with $c_i\in \N$. For each $i$, we define $k_i$ to be the minimum of the coefficients $c_i$ occurring in any such expression for which $\mu$ is dominant and fixed by $\sigma$.  
        
        Now, to reach such $\mu$ from $\lambda'$ by subtracting simple coroots, we have to subtract each $\check{\alpha}_{i_j}$ for $j\in\{1,\ldots, m\}$ at least once, so $k_{i_j}\ge 1$. (Indeed, this is the only way of decreasing $\lambda'_{i_j}-\lambda'_{i_j+1} = 1$.) Doing so yields $(\lambda'-\check{\alpha}_{i_j})_{i_j} - (\lambda'-\check{\alpha}_{i_j})_{i_j+1} = -1$, so to obtain a dominant coweight we must also subtract at least one simple coroot $\check{\alpha}_\ell$ with $|\ell-i_j| = 1$. Though this seems to leave a choice between $\ell = i_j-1$ and $\ell = i_j+1$, for $j=1$ we are in fact forced to use $i_1+1$ eventually. Indeed, if $i_1 = 1$, then $i_1-1$ is not available, and otherwise the difference $\lambda'_{i_1-2}-\lambda'_{i_1-1}$ is already even, which would be destroyed by using $i_1-1$. Repeating this observation shows (inductively) that $k_{i_1+1}\ge 1$. 
        
        After subtracting $\check{\alpha}_{i_1+1}$, we will have changed the parity of $\lambda'_{i_1+2}-\lambda'_{i_1+3}$. If $i_1+2 = i_2$, this is compatible with the necessary subtraction of $\check{\alpha}_{i_2}$; in general a similar inductive argument shows $k_j\ge 1$ for $j=i_1+1,\ldots, i_2$. Repeating this line of reasoning at each $i_{2\ell-1}$ yields $k_{i_{2j-1}}, \ldots, k_{i_{2j}}\ge 1$. This shows that $\mu_\mathrm{max}$ is indeed the unique maximal element of the $\sigma$-fixed dominant coweights below $\lambda$.
    \end{enumerate}

    \textbf{Case \ref{subsubsec:case3}:} 
    Here all elements of $\Lambda^\vee_\R$ are of the form $(\mu_1,0,\ldots, 0)$, $\mu_1\in \Z$. Since the sum of all coordinates of a coroot is even, all elements of $(\Lambda^\vee_+)_{\le \lambda}\cap \Lambda^\vee_\R$ have first coordinates of the same parity. Since we also have $(\mu_1,0,\ldots,0) \ge (\mu_1-2,0,\ldots, 0)$ via subtracting the simple coroots $e_1-e_2$ and $e_1+e_2$, if $(\Lambda^\vee_+)_{\le \lambda}\cap \Lambda^\vee_\R \ne \emptyset$ then required maximal $\mu^\mathrm{max}$ is the one with largest first coordinate. 
    
    \textbf{Case \ref{subsubsec:case4}:} We will need the following expressions of simple coroots in terms of fundamental coweights which can be read off from the Cartan matrix:
    \begin{align*}
        &\check{\alpha}_1 = 2\check{\varpi}_1 - \check{\varpi}_3,\\ &\check{\alpha}_2 = 2\check{\varpi}_2 - \check{\varpi}_4,\\
        &\check{\alpha}_3 = 2\check{\varpi}_3 - \check{\varpi}_1 - \check{\varpi}_4,\\
        &\check{\alpha}_4 = 2\check{\varpi}_4 - \check{\varpi}_2 - \check{\varpi}_3 - \check{\varpi}_5\\
        &\check{\alpha}_5 = 2\check{\varpi}_5 - \check{\varpi}_4 - \check{\varpi}_6,\\
        &\check{\alpha}_6 = 2\check{\varpi}_6 - \check{\varpi}_5.
    \end{align*}
    As we have seen, $\Lambda^\vee_\R = \Z\cdot\{\check{\varpi}_1, \check{\varpi}_6\}$.

    Now let $\lambda = \sum_{i=1}^6 a_i \check{\varpi}_i$, $a_i\in \N$. We study coweights below $\lambda$ from the point of view of the simple coroots subtracted to reach them. For a general $\lambda'\in \Lambda^\vee$, we call a subtraction (of a simple coroot) \emph{necessary}, or more precisely \emph{$\lambda'$-necessary}, if it has to be made in order to eventually reach a $\sigma$-fixed coweight below $\lambda'$. Of course, this notion depends on $\lambda'$, although if $\lambda' \le \lambda$ and there does not exist any $\sigma$-fixed coweight in the interval $[\lambda',\lambda]$, then $\lambda'$-necessity implies $\lambda$-necessity. Let us list some kinds of $\lambda'$-necessary subtractions in terms of $\lambda' = \sum_{i=1}^6 a_i' \check{\varpi}_i$:
\begin{itemize}
    \item[(a)] If $a_i' \ge 1$ for $i\in \{2,3,4,5\}$, then subtracting $\check{\alpha}_i$ is necessary because it is the only simple coroot with positive $\check{\varpi}_i$-coefficient.
    \item[(b)] If $a_2'$ is odd, then subtracting $\check{\alpha}_4$ is necessary. Indeed, the only other way of influencing the coefficient $a_2'$ by simple coroots is via $\check{\alpha}_2$ which has even $\check{\varpi}_2$-coefficient.
    \item[(c)] If $a_3'$ does not have the same parity as $a_2'$, then subtracting $\check{\alpha}_1$ is necessary. Indeed, this is the only way of altering the parity $a_2'+a_3'$ by subtracting a simple root, and we claim that this parity has to be even in order to reach a fixed coweight using only simple coroots different from $\check{\alpha}_1$. Indeed, the other simple coroots correspond to a Levi subalgebra of Dynkin type $D_5$, and this claim then follows from the previous case (by comparing the real root lattices). Similarly, if $a_5'$ does not have the same parity as $a_2'$, then subtracting $\check{\alpha}_6$ is necessary.
\end{itemize}

Consider now the following iterative process: \emph{starting from $\lambda' \coloneqq \lambda$, replace $\lambda'$ by $\lambda'-\check{\alpha}_i$ for some $i$ such that subtracting $\check{\alpha}_i$ is necessary with respect to $\lambda'$.} To prove the proposition in the current case, it is sufficient to show that this process terminates at an element of $\Lambda^\vee_{\R, +}$.\footnote{The process always terminates because all cosets are admissible, and this happens once it reaches any element of $\Lambda^\vee_\R$. Until that point, it will not have seen any elements of $\Lambda^\vee_\R$, so the notion of necessary coweight subtraction does not change until that point. What has to be checked is that the termination point is at a \emph{dominant} real coweight.}

To simplify checking this, we observe that it is clearly possible to achieve $a_i' \in \{0,1\}$ for $2\le i \le 5$ using only necessary subtractions. Below, we check the claim for each of the sixteen possible tuples $(a_2',a_3',a_4',a_5')\in \{0,1\}^4$. To simplify this, we use the following conventions:
\begin{itemize}
    \item All type (c) subtractions are done (simultaneously) at the start.
    \item If several type (a) subtractions are necessary, they are performed simultaneously.
    \item If $a_4' \ge 1$ but $a_2' = 0$, then the necessary subtraction of $\alpha_4$ leads to $a_2' = 1$, hence another necessary subtraction of $\check{\alpha}_4$ and a necessary subtraction of $\check{\alpha}_2$. We perform these at once, thus subtracting $\check{\alpha}_2 + 2\check{\alpha}_4 = -2\check{\varpi}_3 + 3\check{\varpi}_4 - 2 \check{\varpi}_5$ and denote this by (ab)
\end{itemize}

We now come to checking the sixteen cases, using wherever possible the symmetry exchanging the labels $1\leftrightarrow 6$, $3\leftrightarrow 5$. The reader may prefer to check these on the Dynkin diagram (where the process is reminiscent of the classic arcade game whack-a-mole). The coordinates used for $\lambda'$ in the list below refer to fundamental coweights.
\begin{enumerate}
    \item $\lambda' = (a,0,0,0,0,b)$: trivial.
    
    \item\label{item:1000} $\lambda' = (a,1,0,0,0,b)$: perform subtractions
    \begin{align*}
        &(a,1,0,0,0,b) \overset{(c)}{\rightarrow} 
        (a-2,1,1,0,1,b-2) \overset{(b)}{\rightarrow}
        (a-2,2,2,-2,2,b-2) \overset{(a)}{\rightarrow}\\
        \overset{(a)}{\rightarrow}& (a-1,0,0,1,0,b-1)
        \overset{(ab)}{\rightarrow}
        (a-1,0,2,-2,2,b-1) \overset{(a)}{\rightarrow}
        (a,0,0,0,0,b).
    \end{align*}

    \item\label{item:0100} $\lambda' = (a,0,1,0,0,b)$: perform subtractions 
    \begin{align*}
        &(a,0,1,0,0,b) \overset{(c)}{\rightarrow}
        (a-2,0,2,0,0,b) \overset{(a)}{\rightarrow}
        (a-1,0,0,1,0,b) \overset{(ab)}{\rightarrow}\\
        \overset{(ab)}{\rightarrow}&
        (a-1,0,2,-2,2,b) \overset{(a)}{\rightarrow}
        (a,0,0,0,0,b+1).
    \end{align*}

    \item\label{item:0010} $\lambda'=(a,0,0,1,0,b)$: perform subtractions
    \begin{align*}
        &(a,0,0,1,0,b) \overset{(ab)}{\rightarrow}
        (a,0,2,-2,2,b) \overset{(a)}{\rightarrow}
        (a+1,0,0,0,0,b+1).
    \end{align*}
    This case turns out to be very useful as we will see below.

    \item\label{item:0001} $\lambda'=(a,0,0,0,1,b)$: analogous to case \ref{item:0100}.

    \item\label{item:1100} $\lambda'=(a,1,1,0,0,b)$: perform subtractions
    \begin{align*}
        &(a,1,1,0,0,b) \overset{(c)}{\rightarrow}
        (a,1,1,0,1,b-2) \overset{(a)}{\rightarrow}
        (a+1,-1,-1,3,-1,b-1) \overset{(a)}{\rightarrow}\\
        &\overset{(a)}{\rightarrow}
        (a+1,0,0,1,0,b-1),
    \end{align*}
    then go to case \ref{item:0010}.
    
    \item\label{item:1010} $\lambda' = (a,1,0,1,0,b)$: perform subtractions
    \begin{align*}
        &(a,1,0,1,0,b) \overset{(c)}{\rightarrow}
        (a-2,1,1,1,1,b-2) \overset{(a)}{\rightarrow}
        (a-1,0,0,2,0,b-1) \overset{(ab)}{\rightarrow}\\
        &\overset{(a)}{\rightarrow}
        (a-1,0,2,-1,2,b-1) \overset{(a)}{\rightarrow}
        (a,0,0,1,0,b),
    \end{align*}
    then go to case \ref{item:0010}.

    \item\label{item:1001} $\lambda' = (a,1,0,0,1,b)$: analogous to case \ref{item:1100}.

    \item\label{item:0110} $\lambda' = (a,0,1,1,0,b)$: perform subtractions
    \begin{align*}
        &(a,0,1,1,0,b) \overset{(c)}{\rightarrow}
        (a-2,0,2,1,0,b) \overset{(a)}{\rightarrow}
        (a-1,0,0,2,0,b) \overset{(ab)}{\rightarrow}\\
        \overset{(ab)}{\rightarrow}&
        (a-1,0,2,-1,2,b) \overset{(a)}{\rightarrow}
        (a,0,0,1,0,b+1),
    \end{align*}
    then go to case \ref{item:0010}.

    \item\label{item:0101} $\lambda' = (a,0,1,0,1,b)$: perform subtractions
    \begin{align*}
        &(a,0,1,0,1,b) \overset{(c)}{\rightarrow}
        (a-2,0,2,0,2,b-2) \overset{(a)}{\rightarrow}
        (a-1,0,0,2,0,b-1) \overset{(ab)}{\rightarrow}\\
        \overset{(ab)}{\rightarrow}&
        (a-1,0,2,-1,2,b-1) \overset{(a)}{\rightarrow}
        (a,0,0,1,0,b),
    \end{align*}
    then go to case \ref{item:0010}.

    \item\label{item:0011} $\lambda' = (a,0,0,1,1,b)$: analogous to case \ref{item:0110}.

    \item\label{item:1110} $\lambda' = (a,1,1,1,0,b)$: perform subtractions
    \begin{align*}
        &(a,1,1,1,0,b) \overset{(c)}{\rightarrow}
        (a,1,1,1,1,b-2) \overset{(a)}{\rightarrow}
        (a+1,0,0,2,0,b-1) \overset{(ab)}{\rightarrow}\\
        \overset{(ab)}{\rightarrow}&
        (a+1,0,2,-1,2,b-1) \overset{(a)}{\rightarrow}
        (a+2,0,0,1,0,b) \overset{(a)}{\rightarrow},
    \end{align*}
    then go to case \ref{item:0010}.

    \item\label{item:1101} $\lambda'=(a,1,1,0,1,b)$: perform subtractions
    \begin{align*}
        &(a,1,1,0,1,b) \overset{(a)}{\rightarrow}
        (a+1,-1,-1,3,-1,b+1) \overset{(a)}{\rightarrow}
        (a+1,0,0,1,0,b+1),
    \end{align*}
    then go to case \ref{item:0010}.

    \item\label{item:1011} $\lambda' = (a,1,0,1,1,b)$: analogous to case \ref{item:1110}.

    \item\label{item:0111} $\lambda' = (a,0,1,1,1,b)$: perform subtractions
    \begin{align*}
        &(a,0,1,1,1,b) \overset{(c)}{\rightarrow}
        (a-2,0,2,1,2,b-2) \overset{(a)}{\rightarrow}
        (a-1,0,0,3,0,b-1) \overset{(ab)}{\rightarrow}\\
        \overset{(ab)}{\rightarrow}&
        (a-1,0,2,0,2,b-1) \overset{(a)}{\rightarrow}
        (a,0,0,2,0,b) \overset{(ab)}{\rightarrow}
        (a,0,2,-1,2,b) \overset{(a)}{\rightarrow}\\
        \overset{(a)}{\rightarrow}&
        (a+1,0,0,1,0,b+1),
    \end{align*}
    then go to case \ref{item:0010}.

    \item\label{item:1111} $\lambda' = (a,1,1,1,1,b)$: perform subtractions
    \begin{align*}
        &(a,1,1,1,1,b) \overset{(a)}{\rightarrow}
        (a+1,0,0,2,0,b+1) \overset{(ab)}{\rightarrow}
        (a+1,0,2,-1,2,b+1) \overset{(a)}{\rightarrow}\\
        \overset{(a)}{\rightarrow}&
        (a+2,0,0,1,0,b+2),
    \end{align*}
    then go to case \ref{item:0010}.   
\end{enumerate}
This concludes the proof of for the case of $E_{6,\mathrm{ad}}$, which was the last case that remained to be checked.
\end{proof}

\clearpage
\printbibliography

\end{document}